\documentclass{amsproc}

\usepackage{epsfig,graphicx}

\newtheorem{theorem}{Theorem}[section]
\newtheorem{lemma}[theorem]{Lemma}
\newtheorem{corollary}[theorem]{Corollary}
\newtheorem{proposition}[theorem]{Proposition}
\newtheorem{definition}[theorem]{Definition}

\newtheorem{remark}[theorem]{Remark}

\numberwithin{equation}{section}

\newcommand{\bz}{{\mathbb B}}
\newcommand{\cz}{{\mathbb C}}
\newcommand{\dz}{{\mathbb D}}
\newcommand{\gz}{{\mathbb Z}}

\newcommand{\nz}{{\mathbb N}}
\newcommand{\rz}{{\mathbb R}}

\newcommand{\calA}{\mathcal{A}}
\newcommand{\calB}{\mathcal{B}}
\newcommand{\calC}{\mathcal{C}}
\newcommand{\calD}{\mathcal{D}}
\newcommand{\calE}{\mathcal{E}}
\newcommand{\calF}{\mathcal{F}}
\newcommand{\calG}{\mathcal{G}}
\newcommand{\calH}{\mathcal{H}}
\newcommand{\calK}{\mathcal{K}}
\newcommand{\calL}{\mathcal{L}}
\newcommand{\calM}{\mathcal{M}}

\newcommand{\calP}{\mathcal{P}}

\newcommand{\calS}{\mathcal{S}}

\newcommand{\ulC}{\underline{C}}
\newcommand{\ulS}{\underline{S}}
\newcommand{\ulU}{\underline{U}}
\newcommand{\ulV}{\underline{V}}

\newcommand{\as}{\text{\rm As}}
\newcommand{\bsgp}{\mathcal{B}^{s,\gamma}_p}

\newcommand{\bsgpb}{\mathcal{B}^{s,\gamma}_p({\mathbb B})}
\newcommand{\ci}{\mathcal{C}^\infty}
\newcommand{\cicomp}{\mathcal{C}^\infty_{\text{\rm comp}}}
\newcommand{\cig}{\mathcal{C}^{\infty,\gamma}}

\newcommand{\cii}{\mathcal{C}^{\infty,\infty}}

\newcommand{\dbar}{d\hspace*{-0.08em}\bar{}\hspace*{0.1em}}
\newcommand{\diag}{\text{\rm diag}} 
\newcommand{\eps}{\varepsilon}
\newcommand{\hsgp}{\mathcal{H}^{s,\gamma}_p}

\newcommand{\hsgpd}{\mathcal{H}^{s,\gamma}_p({\mathbb D})}

\newcommand{\intb}{\text{\rm int}\,{\mathbb B}}
\newcommand{\intd}{\text{\rm int}\,{\mathbb D}}
\newcommand{\kringel}[1]{\mbox{\tiny$\overset{\circ}{\mbox{\normalsize$#1$}}$}{}}
\newcommand{\op}{\text{\rm op}}
\newcommand{\pit}{\,{\widehat{\otimes}}_\pi\,}
\newcommand{\re}{\text{\rm Re}\,}
\newcommand{\rpbar}{\overline{{\mathbb R}}{}_+}
\newcommand{\schnitt}{\mathop{\mbox{\Large$\cap$}}}
\newcommand{\skp}[2]{(#1,#2)}
\newcommand{\spk}[1]{\left<#1\right>}

\newcommand{\smsum}{\mathop{\mbox{\large$\sum$}}}
\newcommand{\st}{\mbox{\boldmath$\;|\;$\unboldmath}}

\newcommand{\wt}{\widetilde}

\begin{document}

\title[Realizations of Differential Operators on Conic Manifolds with Boundary]{Realizations 
of Differential Operators\\ on Conic Manifolds with Boundary}

\author{S.\ Coriasco}
\address{Universit\'a di Torino, Dipartimento di Matematica,
         Via Carlo Alberto 10, 10123 Torino, Italy}
\email{coriasco@dm.unito.it}

\author{E.\ Schrohe}
\address{Universit\"at Hannover, Institut f\"ur Mathematik, 
         Welfengarten 1, 30167 Hannover, Germany}
\email{schrohe@math.uni-hannover.de}
\author{J.\ Seiler}
\address{Universit\"at Hannover, Institut f\"ur Angewandte Mathematik, 
         Welfengarten 1, 30167 Hannover, Germany}
\email{seiler@ifam.uni-hannover.de}



\begin{abstract}
We study the closed extensions (realizations) of differential operators subject to homogeneous 
boundary conditions on weighted $L_p$-Sobolev spaces over a manifold with boundary and conical 
singularities. Under natural ellipticity conditions we determine the domains of the minimal and 
the maximal extension. We show that both are Fredholm operators and give a formula for the
relative index.
\end{abstract}

\maketitle


\tableofcontents


\section{Introduction}\label{section1}

Operator semigroups are one of the most efficient tools for the analysis of parabolic 
differential equations. In fact, these problems can often be reformulated as abstract 
evolution equations of the form $\dot{u}+Au=f$, $u(0)=u_0$, where $A$ is a closed 
unbounded operator in a Banach space $E$, induced by the original partial differential 
equation. Also semihomogeneous boundary value problems can be treated in this way, namely 
by incorporating the boundary condition into the choice of the domain: Given a boundary 
condition $T$, one studies the operator $A$ on a domain contained in the kernel of $T$. 
In general, one can think of the operator $A$ to be defined initially on a small space of 
functions in $E$. In order to apply the machinery of evolution equations, as a first step one  
has to determine those closed extensions of $A$ in $E$ which reflect the original problem.

In the case of a closed manifold, there are at least two distinguished choices: The 
{\em minimal} extension is the closure of the operator $A$ acting on all smooth functions, 
and the {\em maximal} extension has as domain all those elements of $E$ which are mapped into 
$E$ by $A$. 

For a manifold with boundary, the latter extension is no longer relevant, since it does not 
involve any boundary condition. Instead one considers the closed extensions -- in this case 
also called {\em realizations} -- with respect to a given boundary condition $T$, i.e.\ one 
only considers domains $\calD$ for which $Tu=0$ whenever $u\in\calD$. Of course, boundary 
conditions can only be imposed if one has a certain a priori regularity. As the minimal 
(respectively maximal) extension we therefore choose the closure of the operator acting on 
all smooth (respectively all sufficiently regular) functions which vanish under the boundary 
condition. 

A classical example is the heat equation in a bounded open set $\Omega$ with smooth boundary 
and Dirichlet boundary condition: In order to solve the problem in $L_p(\Omega)$, for example, 
one studies $A=-\Delta$, initially defined on all smooth functions on $\Omega$ vanishing at 
the boundary. Then the closure can be shown to coincide with the maximal extension, given by 
the action of $-\Delta$ on the domain $\{u\in H^2(\Omega): u|_{\partial\Omega}=0\}$. 
Hence there is only one closed realization for the Dirichlet boundary condition, namely the 
one just described. 

For Fuchs type differential operators on manifolds with conical singularities, the situation 
is less simple. The boundaryless case (for $p=2)$ has been analyzed by Lesch \cite{Lesc}. He 
showed that, under a natural ellipticity assumption (corresponding to $\dz$-ellipticity in 
this paper) both the minimal and the maximal extension are Fredholm operators and the 
quotient ${\calD}_{\rm max}/\calD_{\rm min}$ of their domains is  finite-dimensional; its 
dimension can be computed from symbol data (more precisely from the meromorphic structure of 
the conormal symbol) of $A$. Gil and Mendoza \cite{GiMe}, Proposition 3.6, obtained an 
improved  description of $\calD_{\rm min}$, while the structure of the maximal domain was 
studied in Schrohe and Seiler \cite{ScSe2}, Theorem 2.8.

In the present paper, we extend this work to elliptic boundary value problems on manifolds 
with boundary and conical singularities. While our final results -- given in detail in 
a)--d) at the end of this section -- look similar to those in \cite{Lesc}, \cite{GiMe} and 
\cite{ScSe2}, the analysis becomes much more difficult. 

This starts with simple facts: The domains in general are not invariant under multiplication 
by cut-off functions, since the boundary condition then need not remain fulfilled. This limits 
localization techniques. Fortunately, there are suitable projections to handle these problems. 
Also the adjoint of an elliptic boundary value problem is harder to analyze than that of a 
differential operator, a fact which complicates the theory already in the case of smooth 
manifolds. We adapt here a construction of Grubb \cite{Grub1}, Section 1.6, to the conic 
situation. Moreover, a basic reduction in the analysis of the boundaryless situation 
consists in switching to an operator whose coefficients are constant near the cone singularity 
(frozen at the tip of the cone). In the case of boundary value problems, one has to freeze 
both the operator and the boundary condition, and the corresponding changes of the domain 
are more difficult to handle. Finally, there is the fact that instead of the 
pseudodifferential techniques used before for the construction of parametrices we now have 
to employ Boutet de Monvel's calculus. Indeed, our main tool here will be an $L_p$-version 
of the cone calculus for boundary value problems developed by Schrohe and Schulze \cite{ScSc1}, 
\cite{ScSc2}, which we review in Section 2.

This paper lays the foundations for the study of the resolvents of elliptic boundary value 
problems on manifolds with conical singularities and questions of maximal regularity or 
solvability of certain nonlinear equations in the spirit of \cite{ScSe2}, \cite{CSS}.

Let us now explain the contents of the present paper in more detail. The intuitive picture of 
the underlying manifold is given in Figure 1.  
\begin{figure}[!ht]
 \begin{center}
  \includegraphics[width=12cm]{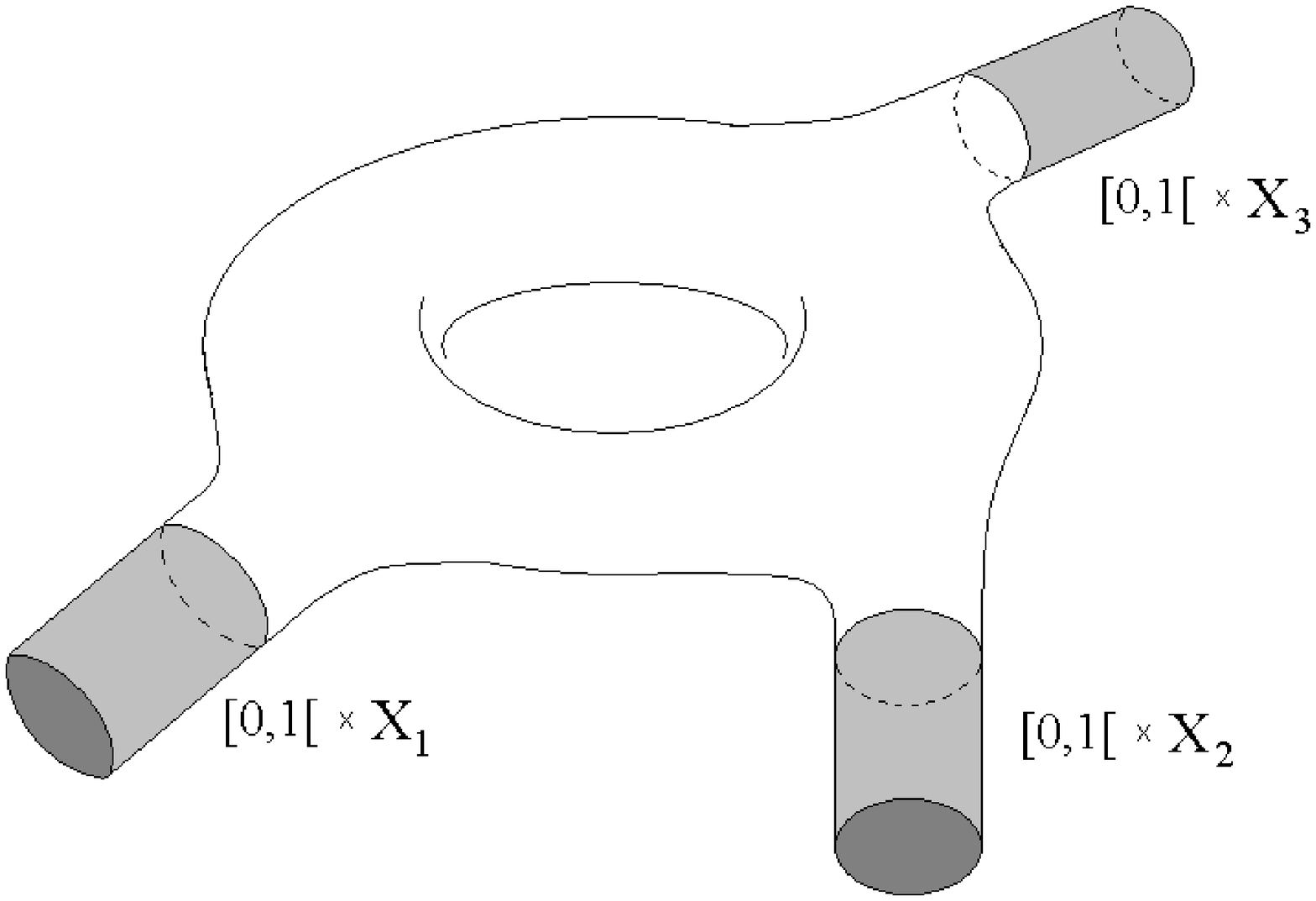}
 \end{center}
\caption{}
\end{figure}
Technically, we denote by $\intd$ an $n$-dimensional riemannian manifold with boundary having 
 finitely many cylindrical ends, i.e.\ there exists a compact set $C$ such that 
$\intd\setminus C$ is the disjoint union $U_1\cup\ldots\cup U_N$, where each $U_i$ is 
isometric to the product ${]0,1[}\times X_i$ for a smooth riemannian manifold $X_i$ with 
boundary; the $X_i$ need not be connected. We fix the coordinate in ${]0,1[}$ in such a way 
that every neighborhood of $C$ in $\intd$ has nonempty intersection with 
${]\frac{1}{2},1[}\times X_i$. We complete $\intd$ with the help of the riemannian metric, 
so that $\dz\setminus C$ can be identified with the disjoint union of all ${[0,1[}\times X_i$. 
We next denote by $\intb$ the boundary of $\intd$ and by $\bz$ its completion in the metric 
inherited from $\intd$. Then $\bz\setminus C$ can be identified with the disjoint union of 
all ${[0,1[}\times\partial X_i$; $\bz$ itself is a smooth manifold with boundary.  

To make the exposition simpler we shall assume in the following that $\dz$ has only one 
cylindrical end, denoted by ${[0,1[}\times X$. Here we use the canonical coordinates $(t,x)$, 
$0\le t<1$, $x\in X$. The subset $\{0\}\times X$ of $\dz$, where $t=0$, is  occasionally 
called the {\em singularity} of $\dz$. In the sequel, a vector bundle over $\dz$ will be a 
smooth hermitian vector bundle over $\intd$ such that $E|_{{]0,1[}\times X}$ is isometric to 
the pull-back (under the canonical projection ${]0,1[}\times X\to X$) of a hermitian bundle 
$E_0$ over $X$. 

That we call $\dz$ a manifold with {\em conical} singularity is due to the class of 
differential operators we consider on it, the so-called {\em Fuchs type} operators. A 
$\mu$-th order differential operator $A$ on $\intd$ with smooth coefficients, acting 
between sections of vector bundles $E$ and $\wt E$ over $\dz$, is of Fuchs type, if, near 
the singularity, 
\begin{equation}\label{real4A.5}
 A=t^{-\mu}\smsum_{j=0}^\mu a_j(t)(-t\partial_t)^j,\qquad
 a_j\in\ci([0,1[,\text{\rm Diff}^{\,\mu-j}(X;E_0,{\wt E}_0)),
\end{equation}
where $E_0$, $\wt E_0$ are the restrictions of $E$, $\wt E$ to $X\cong\{0\}\times X$. The 
three characteristic properties of a Fuchs type operator are the singular factor $t^{-\mu}$ 
corresponding to the order of the operator, coefficients which are smooth up to $t=0$, and 
the totally characteristic differentiation $t\partial_t$ in $t$-direction. As we are treating
elliptic operators, we assume without loss of generality that $E=\wt E$. 

One particular example of such an operator (of order $\mu=2$ and with $E=\intd\times\cz$) is 
the Laplace-Beltrami operator on $\intd$ for a so-called {\em conical} metric, i.e.\ a metric 
which is of the form $dt^2+t^2g(t)$ on ${]0,1[}\times X$ with a smooth (up to $t=0$) family 
$g(t)$ of metrics on $X$. In this case 
 $$\Delta=t^{-2}\,\Big\{(-t\partial_t)^2-(n-1+a(t))(-t\partial_t)+\Delta_X(t)\Big\}$$ 
with $a(t)=\frac{1}{2}t\partial_t\log G(t)$ for $G=\det g_{ij}$, and $\Delta_X(t)$ is the 
Laplacian on $X$ for the metric $g(t)$. 

Fuchs type operators naturally act in a scale of weighted cone Sobolev spaces 
$\calH^{s,\gamma}_p(\dz,E)$ introduced in Definition \ref{real31}, below. Given 
$s,\gamma\in\rz$ and $1<p<\infty$, the operator $A$ defines a continuous map 
$\calH^{s,\gamma}_p(\dz,E)\to\calH^{s-\mu,\gamma-\mu}_p(\dz,E)$.  For $s=0$ the space 
$\calH^{0,\gamma}_p(\dz,E)$ is isomorphic to the weighted space  
$t^{\gamma_p+\frac{n}{p}-\gamma}L_p(\dz,E)$ for $\gamma_p=(n+1)(\frac{1}{2}-\frac{1}{p})$. 

The main objective of this paper is the description of the closed realizations of the operator
\begin{equation}\label{real5B.5}
 A:\cii(\dz,E)_T\subset\calH^{0,\gamma}_p(\dz,E)\longrightarrow\calH^{0,\gamma}_p(\dz,E),
\end{equation} 
where the initial domain consists of all smooth functions that vanish to infinite order in 
the singularity and that vanish under the {\em boundary condition} $T$. We assume that $T$ is 
of the form $T=(T_0,\ldots,T_{\mu-1})$ with $T_k=\gamma_0\circ B_k$, where each $B_k$ is a 
Fuchs type differential operator of order $k$ on $\dz$, acting from sections of $E$ to 
sections of bundles $F_k$; as usual, $\gamma_0$ denotes the operator of restriction to the 
boundary, i.e. $\gamma_0:\cii(\dz,F_k) \to\cii(\bz,{F_k\!|}_\bz)$. We also allow $F_k$ to be 
zero dimensional; in this case the $k$-th condition $T_k$ is void. For more details see 
Section \ref{section4.1}. 

Our main assumption is that  $A$ is an elliptic differential operator on $\intd$, that the 
boundary condition is {\em normal} (in the sense of Grubb \cite{Grub1}, Definition 1.4.3), 
and that $A$ together with $T$ is an elliptic boundary value problem, i.e.\ that its boundary 
symbol satisfies the Shapiro-Lopatinskij condition. From \eqref{real4A.5} it is obvious that 
the principal symbol of $A$ degenerates in a controlled way as $t$ tends to 0; the same is 
true for the boundary symbol. After suitable rescaling, we can pass to limit symbols at $t=0$ 
which we require to be invertible. We call this {\em $\dz$-ellipticity}.

We shall denote by $A_{T,\min}$ the closure of \eqref{real5B.5} in the Banach space  
$\calH^{0,\gamma}_p(\dz,E)$, and by $A_{T,\max}$ the unbounded operator  acting as $A$ on the 
domain  $\{u\in\calH^{\mu,\gamma}_p(\dz,E)\st Tu=0\text{ and }Au\in\calH^{0,\gamma}_p(\dz,E)\}$. 
Assuming $\dz$-ellipticity we prove the following: 
\begin{itemize}
 \item[a)] The domain of $A_{T,\min}$ is 
  $\{u\in\schnitt\limits_{\eps>0}\calH^{\mu,\gamma+\mu-\eps}_p(\dz,E)\st Tu=0\text{ and } 
  Au\in\calH^{0,\gamma}_p(\dz,E)\}$; for $\gamma$ outside a discrete set this coincides with 
  $\{u\in\calH^{\mu,\gamma+\mu}_p(\dz,E)\st Tu=0\}$,
 \item[b)] $\calD(A_{T,\max})=\calD(A_{T,\min})+\calE$ for a finite dimensional space $\calE$ 
  of smooth functions that is determined by the (conormal) symbolic structure of $A$ and $T$, 
 \item[c)] both $A_{T,\min}$ and $A_{T,\max}$ are Fredholm operators; the difference of the
  indices can be computed,
 \item[d)] $A_{T,\min}^*$, the adjoint of the closure, can be described explicitly as 
  $(A^t)_{\wt{T},\max}$ for the formal adjoint $A^t$ of $A$ and a resulting boundary condition 
  $\wt T$. Correspondingly, $A_{T,\max}^*=(A^t)_{\wt{T},\min}$. 
\end{itemize}
As mentioned above, Lesch \cite{Lesc} treated the case of conical manifolds without boundary. 
Gil and Mendoza \cite{GiMe} determined the domain of adjoints of closed extensions of Fuchs 
type differential operators; in joint work \cite{GLM} with Loya they proved an index theorem. 
The articles \cite{ScSe2} and \cite{CSS} focus on the descripton of the resolvent of closed 
extensions and show existence of bounded imaginary powers. The analysis of boundary value 
problems on manifolds with conical singularities began with Kondratievs fundamental paper 
\cite{Kond} and has evolved in many directions. A topic related to this investigation is the 
work on the asymptotics of the solutions to elliptic equations on singular domains by Mazya, 
Nazarov, and Plamenevskij \cite{MNP}, \cite{NP}.

{\bf Notation:} Throughout this paper, a {\em cut-off function} is a non-negative, 
monotonically decreasing function in $\cicomp([0,1[)$ that is identically 1 in a neighborhood 
of 0. In the sequel, $\omega,\omega_0,\omega_1,\ldots$, and 
$\wt\omega,\wt\omega_0,\wt\omega_1,\ldots$ will denote such cut-off functions. We shall write 
$\omega_0\prec\omega$ if $\omega\equiv1$ in a neighborhood of the support of $\omega_0$. 
The operator of multiplication (in various spaces) with $\omega$ shall be denoted by $\omega$, 
again. 

By $t$ we denote a real variable as well as a smooth positive function on $\intd$ that 
coincides on the cylindrical part ${]0,1[}\times X$ with $(t,x)\mapsto t$.


\section{Boutet de Monvel's algebra on manifolds with conical singularities}\label{section3}

We review Boutet de Monvel's algebra for manifolds with conical singularities. Main references 
are \cite{ScSc1} and \cite{ScSc2}, where this algebra is studied in detail. For other 
expositions we refer to \cite{KaSc} and \cite{Schu2}. We assume some familiarity with Boutet's 
calculus on smooth manifolds. However, for completeness, we give an introduction to this 
calculus in the appendix, cf.\ Section \ref{section2}.  

In a slight extension of the $L_2$-based previous work, we shall consider the operators on 
$L_p$-spaces with $1<p<\infty$. This requires two continuity results which we deduce from 
\cite{GrKo}. 

\subsection{Weighted distribution spaces on $\dz$ and $\bz$}\label{section3.1}
Let $X$ be embedded in $\Omega$, a smooth compact manifold without boundary. Using the product 
structure of $\rz\times\Omega$, it is straightforward to define the Sobolev spaces 
 $$H^s_p(\rz\times\Omega),\qquad s\in\rz,\;1<p<\infty.$$
Extending the map $S_\gamma:\cicomp(\rz\times\Omega)\to\cicomp(\rz_+\times\Omega)$, 
\begin{equation}\label{change}
 (S_\gamma u)(t,x):=t^{-\frac{n+1}{2}+\gamma}u(\log t,x),\qquad \gamma\in\rz,
\end{equation}
from functions to distributions, we obtain the weighted Sobolev spaces 
 $$\hsgp(\rz_+\times\Omega),\qquad s,\gamma\in\rz,\;1<p<\infty$$
on the half-cylinder as the image of  $H^s_p(\rz\times\Omega)$ under $S_\gamma$ with the 
canonically induced norm. Restricting to $\rz_+\times X$ we obtain $\hsgp(\rz_+\times X)$ with 
the quotient norm 
 $$\|u\|_{\calH^{s,\gamma}_p(\rz_+\times X)}=\inf\Big\{
   \|v\|_{\calH^{s,\gamma}_p(\rz_+\times\Omega)}\st v|_{\rz_+\times\textrm{int}\,X}=u\Big\}.$$

Gluing together two copies of $\dz$ along the singularity $\{0\}\times X$ yields a smooth 
manifold with boundary $2\dz$ and the standard scale of Sobolev spaces $H^s_p(2\dz)$. 

\begin{definition}\label{real31}
For $s,\gamma\in\rz$ and $1<p<\infty$ let 
 $$\hsgpd=\left\{u\in\calD^\prime(\kringel{\dz})\st
   \omega u\in\hsgp(\rz_+\times X)\text{ and }
   (1-\omega)u\in H^s_p(2\dz)\right\},$$
where $\omega$ is an arbitrary fixed cut-off function and $\kringel{\dz}$ denotes the interior 
of $\intd$. The norm is given by 
 $$\|u\|_{\hsgpd}=\|\omega u\|_{\hsgp(\rz_+\times X)}+
   \|(1-\omega)u\|_{H^s_p(2\dz)}.$$
We  denote the closure of $\cicomp(\kringel{\dz})$ in $\hsgpd$ by 
$\kringel{\calH}^{s,\gamma}_p(\dz)$. 
\end{definition} 

The spaces $\hsgpd$ and $\kringel{\calH}^{s,\gamma}_p(\dz)$ are Banach spaces and, in case 
$p=2$, even Hilbert spaces. While the index $s$ indicates the smoothness of a distribution, 
the index $\gamma$ indicates its {\em flatness} towards the singularity, since
\begin{equation}\label{weight}
 u\in\calH^{s,\gamma}_p(\dz)\iff t^{-\gamma}u\in\calH^{s,0}_p(\dz).
\end{equation}
The standard embedding and duality properties of Sobolev spaces have corresponding analogues 
for cone Sobolev spaces: 

\begin{remark}\label{real32}
We have continuous embeddings 
 $$\calH^{s^\prime,\gamma^\prime}_p(\dz)\hookrightarrow\calH^{s,\gamma}_p(\dz),\qquad 
   \calH^{r,\gamma}_p(\dz)\hookrightarrow\calH^{t,\gamma}_q(\dz),$$
provided $s^\prime\ge s$, $\gamma^\prime\ge\gamma$, and $q\ge p$, 
$r-\frac{n+1}{p}\ge t-\frac{n+1}{q}$. The first embedding is compact in case of 
$s^\prime>s$ and $\gamma^\prime>\gamma$. Via the scalar product of $\calH^{0,0}_2(\dz)$ the 
spaces $\calH^{s,\gamma}_p(\dz)$ and $\kringel{\calH}^{-s,-\gamma}_{p^\prime}(\dz)$ are dual 
to each other if $\frac{1}{p}+\frac{1}{p^\prime}=1$. 
\end{remark}

Similarly, using the Besov spaces $B^s_p(\rz\times\partial X):=B^s_{p,p}(\rz\times\partial X)$ 
on the cylinder and the transformation 
\begin{equation}\label{change2}
 (S_\gamma^\prime u)(t,x^\prime):=t^{-\frac{n}{2}+\gamma}u(\log t,x^\prime),\qquad\gamma\in\rz,  
\end{equation}
we introduce the Banach spaces (Hilbert spaces in case $p=2$)
\begin{equation}\label{real3A}
 \bsgpb=\left\{u\in\calD^\prime(\intb)\st
   \omega u\in\bsgp(\rz_+\times\partial X)\text{ and }
   (1-\omega)u\in B^s_p(2\bz)\right\}. 
\end{equation}

In a neighborhood in $\intd$ of the boundary $\intb$ we fix a normal coordinate (using the 
product structure on the cylindrical part) and a normal derivative $\partial_\nu$. By 
$\gamma_j$ we denote the usual boundary operator $\gamma_0\circ\partial_\nu^j$, where 
$\gamma_0$ denotes restriction to the boundary. From standard theorems on restriction of 
Sobolev spaces, the following lemma follows immediately.

\begin{lemma}\label{real33}
For any $1<p<\infty$ and $s>\frac{1}{p}+j$ the boundary operator $\gamma_j$ induces continuous 
maps  
 $$\gamma_j:\hsgpd\longrightarrow\calB^{s-j-\frac{1}{p},\gamma-\frac{1}{2}}_p(\bz).$$
\end{lemma}

Note that the weight $\gamma$ is shifted to the weight $\gamma-\frac{1}{2}$, since the 
definition of the weighted spaces involves the dimension of the underlying manifold. 

For later purpose we also define subspaces with asymptotics. 

\begin{definition}\label{real34}
For $\gamma\in\rz$ and $\theta>0$ let 
$\as(X;\gamma,\theta)$ denote the set of all finite sets 
 $$P=\Big\{(p,m,M)\st-\gamma-\theta<\re p-\frac{n+1}{2}<-\gamma,\;m\in\nz_0,\;M\subset
   \ci(X)\Big\},$$
where $M$ is a finite dimensional vector spaces. Any such $P$ is called an 
{\em asymptotic type}. We assume that to any $p\in\cz$ there exists at most one element 
$(p,m,M)$ in $P$. Similarly, $\as(\partial X;\gamma,\theta)$ consists of all finite sets   
 $$Q=\Big\{(q,l,L)\st-\gamma-\theta<\re q-\frac{n}{2}<-\gamma,\;l\in\nz_0,\;L\subset
     \ci(\partial X)\Big\}$$
with finite dimensional spaces $L$. We set 
 $$\as(X,\partial X;\gamma,\theta)=\Big\{(P,Q)\st P\in\as(X;\gamma,\theta)\textrm{ and }
   Q\in\as(\partial X;\gamma-\frac{1}{2},\theta)\Big\}.$$
\end{definition}

With $P\in\as(X;\gamma,\theta)$ we associate a finite dimensional space 
$\calE_P(\rz_+\times X)$ of smooth functions: Its elements are all functions $u$ with
 $$u(t,x)=\omega(t)\smsum_{(p,m,M)\in P}\smsum_{j=0}^m c_{pj}(x)t^{-p}\log^jt,$$
where $c_{pj}\in M$ and $\omega$ is a fixed cut-off function. We consider $u$ as a function 
both on $\intd$ and $\rz_+\times X$. In the same way we obtain the spaces 
$\calE_Q(\rz_+\times\partial X)$ for $Q\in\as(\partial X;\gamma,\theta)$. 

\begin{definition}\label{real35}
For $s,\gamma\in\rz$, $1<p<\infty$, and $\theta>0$ set 
 $$\calH^{s,\gamma}_{p,\theta}(\dz)=\schnitt_{\varepsilon>0}
   \calH^{s,\gamma+\theta-\varepsilon}_{p}(\dz),\qquad
   \calB^{s,\gamma}_{p,\theta}(\bz)=\schnitt_{\varepsilon>0}
   \calB^{s,\gamma+\theta-\varepsilon}_{p}(\bz).$$
For $P\in\as(X;\gamma,\theta)$ and $Q\in\as(\partial X;\gamma,\theta)$ we set 
 $$\calH^{s,\gamma}_{p,P}(\dz)=\calH^{s,\gamma}_{p,\theta}(\dz)\oplus\calE_P(\rz_+\times X),
   \qquad \calB^{s,\gamma}_{p,Q}(\bz)=
   \calB^{s,\gamma}_{p,\theta}(\dz)\oplus\calE_Q(\rz_+\times\partial X).$$
\end{definition}

All these space carry a natural Fr\'{e}chet topology as a projective limit of Banach spaces. 

\begin{definition}\label{real35.5}
For $\gamma\in\rz$ set 
 $$\cig(\dz)=\Big\{u\in\ci(\intd)\st(t\partial_t)^k\log^lt(\omega u)\in
            \schnitt_{s\in\rz}\calH^{s,\gamma}_p(\dz)\quad\forall\;k,l\in\nz_0\Big\},$$
where $\omega$ is some cut-off function (the definition is independent of the involved  
$1<p<\infty$: if one transports $\omega u$ to $\rz\times X$ via $S_\gamma$ from \eqref{change}, 
the resulting smooth function is rapidly decreasing in $t$). For $\theta>0$ and 
$P\in\as(X;\gamma,\theta)$ we then set 
 $$\calC^{\infty,\gamma}_\theta(\dz):=\schnitt_{\eps>0}
   \calC^{\infty,\gamma+\theta-\eps}(\dz),\qquad
   \calC^{\infty,\gamma}_P(\dz):=\calC^{\infty,\gamma}_\theta(\dz)\oplus
   \calE_P(\rz_+\times X).$$
Analogously we introduce $\calC^{\infty,\gamma}(\bz)$, $\calC^{\infty,\gamma}_\theta(\bz)$, 
and $\calC^{\infty,\gamma}_Q(\bz)$ for $Q\in\as(\partial X;\gamma,\theta)$. 
\end{definition}
 
By the previous definitions and Lemma \ref{real33} we obtain: 

\begin{lemma}\label{real36}
Let $P\in\as(X;\gamma,\theta)$. For any $1<p<\infty$ and $s>\frac{1}{p}+j$ we have 
 $$\gamma_j:\calH^{s,\gamma}_{p,P}(\dz)\longrightarrow
   \calB^{s-j-\frac{1}{p},\gamma-\frac{1}{2}}_{p,Q}(\bz)$$
with $Q=\{(p,m,\gamma_j M)\st (p,m,M)\in P\}\in\as(\partial X;\gamma-\frac{1}{2},\theta)$, 
where $\gamma_jM$ is the trace space of $M$ under $\gamma_j:\ci(X)\to\ci(\partial X)$. 
\end{lemma}

\subsection{Integral operators}
In the sequel we shall say that 
$\calG=\begin{pmatrix}\calG_{11}&\calG_{12}\\ \calG_{21}&\calG_{22}\end{pmatrix}$ is an 
{\em integral operator with kernel $k$ with respect to the scalar product in} 
$\calH^{0,0}_2(\dz)\oplus\calB^{-\frac{1}{2},-\frac{1}{2}}_2(\bz)$, if $\calG$ acts on 
$\cicomp(\intd)\oplus\cicomp(\intb)$ (and then possibly extends by continuity to other spaces) 
by 
\begin{equation}\label{real3C.5}
 (\calG_{j1}u_1+\calG_{j2}u_2)(y_j)=
   \skp{k_{j1}(y_j,\cdot)}{\overline{u}_1}_{\calH^{0,0}_2(\dz)}+
   \skp{k_{j2}(y_j,\cdot)}{\overline{u}_2}_{\calB^{-1/2,-1/2}_2(\bz)},
   \qquad j=1,2.
\end{equation}
In \eqref{real3C.5}, $y_1$ denotes the variable of $\dz$ and $y_2$ that of $\bz$. The kernels 
$k_{11}$, $k_{12}$, $k_{21}$, and $k_{22}$ are locally integrable on $\intd\times\intd$, 
$\intd\times\intb$, $\intb\times\intd$, and $\intb\times\intb$, respectively. In view of the 
use of the scalar product of $\calB^{-\frac{1}{2},-\frac{1}{2}}_2(\bz)$ in the second 
component, this notion differs slightly from the standard notion of integral operators with a 
kernel $k$. However, this form makes the exposition easier. We shall use the short-hand notation 
 $$k=\begin{pmatrix}k_{11}&k_{12}\\k_{21}&k_{22}\end{pmatrix}\in
   \begin{pmatrix}V_1\\ \oplus\\V_2\end{pmatrix}\pit
   (W_1\oplus W_2),$$
if $k_{ij}\in V_i\pit W_j$, $1\le i,j\le2$, for certain subspaces $V_1$, $W_1$ of smooth 
functions on $\intd$ and subspaces $V_2$, $W_2$ of smooth functions on $\intb$. Moreover, 
$\pit$ denotes the completed projective tensor product. If the spaces $V_i$, $W_j$ are 
Fr\'{e}chet spaces, $V_i\pit W_j$ consists of all functions $f$ that have a representation 
 $$f(y_i,y_j)=\smsum_{k=1}^\infty\,\lambda_k\,v_i^k(y_i)\,w_j^k(y_j)$$
with $(\lambda_k)_{k\in\nz}$ being an absolutely summable sequence of complex numbers, and 
$(v^k_i)_{k\in\nz}$, $(w^k_j)_{k\in\nz}$ being zero sequences in $V_i$ and $W_j$, respectively. 

\subsection{The flat cone algebra}\label{section3.2}
A {\em flat Green operator} $\calG$ of type 0 is an integral operator 
with respect to the scalar product in 
$\calH^{0,0}_2(\dz)\oplus\calB^{-\frac{1}{2},-\frac{1}{2}}_2(\bz)$ with kernel in 
 $$\begin{pmatrix}\cii(\dz)\\ \oplus\\ \cii(\bz)\end{pmatrix}\pit  
   \left(\cii(\dz)\\ \oplus\\ \cii(\bz)\right),$$
where 
\begin{equation}\label{real3C.7}
 \cii(\dz)=\schnitt\limits_{\gamma\in\rz}\calC^{\infty,\gamma}(\dz),\qquad
 \cii(\bz)=\schnitt\limits_{\gamma\in\rz}\calC^{\infty,\gamma}(\bz)
\end{equation}
(these functions vanish to infinite order in the singularity). In particular, $\calG$ induces 
an operator 
\begin{equation}\label{real3D}
 \calG:\begin{matrix}\cii(\dz)\\ \oplus\\ \cii(\bz)\end{matrix}\longrightarrow
 \begin{matrix}\cii(\dz)\\ \oplus\\ \cii(\bz)\end{matrix}. 
\end{equation}

\begin{definition}\label{real37}
An operator $\calG$, acting as in \eqref{real3D}, is called a {\em flat Green operator} of 
type $d\in\nz_0$ if it is of the form 
 $$\calG=\smsum_{j=1}^d \calG_j \begin{pmatrix}D^j&0\\0&1\end{pmatrix}+\calG_0$$
with flat Green operators $\calG_0,\ldots,\calG_d$ of type 0, and a differential operator $D$ 
on $\dz$ with coefficients smooth up to $t=0$, that coincides with $-i\partial_\nu$ near the 
boundary $\bz$ of $\dz$. The space of all such operators is denoted by $C^d_G(\dz)_\infty$. 
\end{definition}

\begin{definition}\label{real38}
A {\em holomorphic Mellin symbol} of order $\mu\in\gz$ and type $d\in\nz_0$ is a holomorphic 
function $h:\cz\to B^{\mu,d}(X)$ such that 
 $$\delta\mapsto h(\delta+i\tau):\rz\longrightarrow B^{\mu,d}(X;\rz_\tau)$$
is a continuous function. We denote the space of all such symbols by $M^{\mu,d}_O(X)$ and set 
 $$M^{\mu,d}_O(\rpbar\times X)=\ci(\rpbar)\pit M^{\mu,d}_O(X).$$
\end{definition}

Each $h\in M^{\mu,d}_O(\rpbar\times X)$ defines an operator 
 $$\op_M(h):
   \begin{matrix}\cicomp(\rz_+\times X)\\ \oplus\\ \cicomp(\rz_+\times\partial X)\end{matrix}
   \longrightarrow
   \begin{matrix}\ci(\rz_+\times X)\\ \oplus\\ \ci(\rz_+\times\partial X)\end{matrix}$$
by 
\begin{equation}\label{real3E}
(\op_M(h)u)(t)=\int_\Gamma t^{-z}h(t,z)(\calM u)(z)\,\dbar z,
\end{equation}
where $\calM$ denotes the Mellin transform, $\Gamma$ is an arbitrary vertical line in the 
complex plane, and we identify 
$\ci(\rz_+\times X)\oplus\ci(\rz_+\times\partial X)$ with 
$\ci(\rz_+,\ci(X)\oplus\ci(\partial X))$. Observe the following useful relation: For any 
$\sigma\in\rz$,  
\begin{equation}\label{210}
 \op_M(h)\,t^{-\sigma}=t^{-\sigma}\,\op_M(T^\sigma h),
\end{equation}
where $T^\sigma$ is the operator of shifting $z$ by $\sigma$, i.e. 
\begin{equation}\label{211}
 (T^\sigma h)(t,z)=h(t,z+\sigma). 
\end{equation}

\begin{definition}\label{real39}
With $\mu\in\gz$, $d\in\nz_0$ let $C^{\mu,d}_O(\dz)$ denote the space of all operators 
 $$\calA:\begin{matrix}\cii(\dz)\\ \oplus\\ \cii(\bz)\end{matrix}\longrightarrow
   \begin{matrix}\cii(\dz)\\ \oplus\\ \cii(\bz)\end{matrix}$$
of the form 
 $$\calA=\omega\,t^{-\mu}\,\op_M(h)\,\omega_0+(1-\omega)\calP(1-\omega_1)+\calG,$$
where $\omega_1\prec\omega\prec\omega_0$ are cut-off functions and 
 \begin{itemize}
  \item[i)] $h(t,z)\in M^{\mu,d}_O(\rpbar\times X)$ is a holomorphic Mellin symbol,
  \item[ii)] $\calP\in B^{\mu,d}(2\dz)$ is an element of Boutet's algebra on $2\dz$, 
  \item[iii)] $\calG\in C^d_G(\dz)_\infty$ is a flat Green operator. 
 \end{itemize}
\end{definition}

\subsection{The full cone algebra}\label{section3.3}
In order to allow the construction of a parametrix to elliptic boundary value problems it is 
necessary to enlarge the flat cone algebra by a more general class of smoothing remainders. 

Flat Green operators of type 0 have, roughly speaking, integral kernels that vanish to 
infinite order in the conical singularity. General Green operators instead have a specific 
kind of asymptotic behaviour near the singularity. 

A {\em Green operator} of type 0 associated with asymptotic types 
$(P,Q)\in\as(X,\partial X;-\gamma,\theta)$ and 
$(P^\prime,Q^\prime)\in\as(X,\partial X;\gamma^\prime,\theta^\prime)$ 
is an integral operator with respect to the scalar product in 
$\calH^{0,0}_2(\dz)\oplus\calB^{-\frac{1}{2},-\frac{1}{2}}_2(\bz)$ with a kernel belonging to 
the intersection of the two spaces 
\begin{equation*} 
 \begin{pmatrix}
  \calC^{\infty,\gamma^\prime}_{P^\prime}(\dz)\\ \oplus\\ 
  \calC^{\infty,\gamma^\prime-\frac{1}{2}}_{Q^\prime}(\bz)   
 \end{pmatrix}\pit  
   \left(\calC^{\infty,-\gamma}(\dz)\\ \oplus\\ 
   \calC^{\infty,-\gamma-\frac{1}{2}}(\bz)\right)
\end{equation*}
and
\begin{equation*} 
 \begin{pmatrix}
  \calC^{\infty,\gamma^\prime}(\dz)\\ \oplus\\ 
  \calC^{\infty,\gamma^\prime-\frac{1}{2}}(\bz)   
 \end{pmatrix}\pit  
   \left(\calC^{\infty,-\gamma}_{\overline{P}}(\dz)\\ \oplus\\ 
   \calC^{\infty,-\gamma-\frac{1}{2}}_{\overline{Q}}(\bz)\right),
\end{equation*} 
where $\overline{P}=\{(\overline{p},m,\overline{M})\st(p,m,M)\in P\}$ and analogously 
$\overline{Q}$ is given; $\overline{M}$ denotes the space of all complex conjugates of 
elements of $M$. 

Using the technique from \cite{Seil}, Section 4, it can be seen that Green operators 
equivalently can be characterized by mapping properties in Sobolev (Besov) spaces: 

\begin{theorem}\label{real310}
For an operator $\calG$ to have a kernel as above, it is necessary and sufficient that, for 
some fixed $1<p<\infty$, 
 $$\calG:\begin{matrix}
      \hsgpd\\ \oplus\\ \calB^{r,\gamma-\frac{1}{2}}_p(\bz)
     \end{matrix}
     \longrightarrow
     \begin{matrix}
      \calC^{\infty,\gamma^\prime}_{P^\prime}(\dz)\\ \oplus\\ 
      \calC^{\infty,\gamma^\prime-\frac{1}{2}}_{Q^\prime}(\bz)
     \end{matrix}$$
for all $s>-1+\frac{1}{p}$ and all $r\in\rz$, and that
 $$\calG^t:
     \begin{matrix}
      \calH^{s,-\gamma^\prime}_{p^\prime}(\dz)\\ \oplus\\ 
      \calB^{r,-\gamma^\prime-\frac{1}{2}}_{p^\prime}(\bz)
     \end{matrix}
     \longrightarrow
     \begin{matrix}
      \calC^{\infty,-\gamma}_{P}(\dz)\\ \oplus\\ 
      \calC^{\infty,-\gamma-\frac{1}{2}}_{Q}(\bz)
     \end{matrix}$$
for all $s>-1+\frac{1}{p^\prime}$ and all $r\in\rz$. Here, the exponent $t$ refers to the 
formal adjoint of $\calG$ with respect to the scalar product in 
$\calH^{0,0}_2(\dz)\oplus\calB^{-\frac{1}{2},-\frac{1}{2}}_2(\bz)$ and $p^\prime$ is the dual 
number to $p$. 
\end{theorem}

\begin{definition}\label{real311}
Let $(P,Q)\in\as(X,\partial X;-\gamma,\theta)$, 
$(P^\prime,Q^\prime)\in\as(X,\partial X;\gamma^\prime,\theta^\prime)$, and $d\in\nz_0$. Then 
$C_G^d(\dz;\gamma,\theta;\gamma^\prime,\theta^\prime)$ denotes the space of all operators 
$\calG$ that are of the form 
 $$\calG=\smsum_{j=1}^d \calG_j \begin{pmatrix}D^j&0\\0&1\end{pmatrix}+\calG_0$$
with operators $\calG_0,\ldots,\calG_d$ having a kernel as described above (with certain 
asymptotic types depending on $\calG_j$), and $D$ as in Definition \mbox{\rm\ref{real37}}. We 
write  $C_G^d(\dz;\gamma,\gamma^\prime,\theta)= C_G^d(\dz;\gamma,\theta;\gamma^\prime,\theta)$. 
\end{definition}

Besides Green operators we need another kind of smoothing remainders, built on 
{\em meromorphic Mellin symbols}. 

\begin{definition}\label{real311.5}
An {\em asymptotic type} for Mellin symbols $P$ of type $d\in\nz_0$ is a set of triples 
$(p,n,N)$ with $p\in\cz$, $n\in\nz_0$, and $N$ a finite dimensional subspace of finite
rank operators from $B^{-\infty,d}(X)$. Moreover, we require that, for each $\delta>0$,
  $$\{p\in\cz\st |\re p|\le\delta\text{ and }(p,n,N)\in P\text{ for some }n,\,N\}$$
is a finite set. A meromorphic Mellin symbol with asymptotic type $P$ of type $d\in\nz_0$ is 
a meromorphic function $f:\cz\to B^{-\infty,d}(X)$ with poles at most in points $p\in\cz$ with 
$(p,n,N)\in P$ for some $n,N$. Moreover $f$ satisfies: If $(p,n,N)\in P$, then the principal 
part of the Laurent series of $f$ in $p$ is of the form
$\sum\limits_{k=0}^nR_k(z-p)^{-k-1}$ with $R_k\in N$; if $\chi\in\ci(\cz)$ is identically zero 
in a small neighborhood around each pole and identically 1 outside another neighborhood, then
$(\chi f)(\delta+i\tau)$ is a continuous function of $\delta\in\rz$ with values in 
$B^{-\infty,d}(X;\rz_\tau)=\calS(\rz_\tau,B^{-\infty,d}(X))$. We shall then write 
$f\in M_P^{-\infty,d}(X)$. 
\end{definition}
 
As in \eqref{real3E} we can associate with meromorphic Mellin symbols a pseudodifferential 
operator. Now, however, the operator will depend on the choice of the line $\Gamma$. Denoting 
the vertical line $\re z = \frac{1}{2}-\delta$ by $\Gamma_{\frac{1}{2}-\delta}$, we define 
$\op_M^{\delta}(f)$ by 
  \begin{equation*}
   (\op_M^{\delta}(f)u)(t)=\int_{\Gamma_{\frac{1}{2}-\delta}}t^{-z}
   f(t,z)(\calM u)(z)\,\dbar z.
  \end{equation*}
Using this notation, we always require that none of the poles of $f$ lies on the chosen line. 

\begin{definition}\label{real312}
Let $\gamma\in\rz$, $\mu\in\gz$, $d\in\nz_0$, and $k\in\nz$. Then 
$C^{\mu,d}_{M+G}(\dz;\gamma,\gamma-\mu,k)$ is the space of all operators 
 $$\calA:\begin{matrix}
      \calC^{\infty,\gamma}(\dz)\\ \oplus\\ 
      \calC^{\infty,\gamma-\frac{1}{2}}(\bz)
     \end{matrix}
     \longrightarrow 
     \begin{matrix}
      \calC^{\infty,\gamma-\mu}(\dz)\\ \oplus\\ 
      \calC^{\infty,\gamma-\mu-\frac{1}{2}}(\bz)
     \end{matrix}$$
which are of the form 
\begin{equation}\label{mellin}
 \calA=\omega\,\smsum_{j=0}^{k-1}\,
 t^{-\mu+j}\,\op_M^{\gamma_j-\frac{n}{2}}(f_j)\,\omega_0+\calG,
\end{equation}
where $\omega,\omega_0$ are some cut-off functions, $\calG\in C^d_G(\dz;\gamma,\gamma-\mu,k)$, 
the $f_j$ are meromorphic Mellin symbols with asymptotic types $P_j$ of type $d$, and 
$\gamma-j\le\gamma_j\le\gamma$. 
\end{definition}

Now we are in the position to define the full cone algebra.

\begin{definition}\label{real313}
Let $\gamma\in\rz$, $\mu\in\gz$, $d\in\nz_0$, and $k\in\nz$. We set
 $$C^{\mu,d}(\dz;\gamma,\gamma-\mu,k)=C^{\mu,d}_O(\dz)+
   C^{\mu,d}_{M+G}(\dz;\gamma,\gamma-\mu,k),$$ 
acting as in Definition \mbox{\rm\ref{real312}}.  
\end{definition}

The terminology `algebra' is due to the following fact: 

\begin{theorem}\label{real314}
Composition of operators, $(\calA_0,\calA_1)\mapsto \calA_0\calA_1$, induces a map
 $$C^{\mu_0,d_0}(\dz;\gamma-\mu_1,\gamma-\mu_0-\mu_1,k)\times
   C^{\mu_1,d_1}(\dz;\gamma,\gamma-\mu_1,k)
   \longrightarrow C^{\mu,d}(\dz;\gamma,\gamma-\mu,k),$$
with $\mu=\mu_0+\mu_1$ and $d=\max(\mu_1+d_0,d_1)$. If one of the factors $\calA_0$ or 
$\calA_1$ belongs to the $C_{M+G}$- or $C_G$-classes, the same holds true for the product.   
\end{theorem}

\begin{theorem}\label{real315}
Each $\calA\in C^{\mu,d}(\dz;\gamma,\gamma-\mu,k)$ extends by continuity to mappings 
\begin{equation}\label{map1} 
 \calA:\begin{matrix}
      \hsgpd\\ \oplus\\ \calB^{s-\frac{1}{p},\gamma-\frac{1}{2}}_p(\bz)
     \end{matrix}
     \longrightarrow
     \begin{matrix}
      \calH^{s-\mu,\gamma-\mu}_p(\dz)\\ \oplus\\ 
      \calB^{s-\mu-\frac{1}{p},\gamma-\mu-\frac{1}{2}}_p(\bz)
     \end{matrix},\qquad s>d-1+\frac{1}{p}.
\end{equation}
For each asymptotic type $(P,Q)\in\as(X,\partial X;\gamma,k)$ there exists a  
$(P^\prime,Q^\prime)\in\as(X,\partial X;\gamma-\mu,k)$ such that 
\begin{equation}\label{map2} 
 \calA:\begin{matrix}
      \calH^{s,\gamma}_{p,P}(\dz)\\ \oplus\\
      \calB^{s-\frac{1}{p},\gamma-\frac{1}{2}}_{p,Q}(\bz)
     \end{matrix}
     \longrightarrow
     \begin{matrix}
      \calH^{s-\mu,\gamma-\mu}_{p,P^\prime}(\dz)\\ \oplus\\ 
      \calB^{s-\mu-\frac{1}{p},\gamma-\mu-\frac{1}{2}}_{p,Q^\prime}(\bz)
     \end{matrix},\qquad s>d-1+\frac{1}{p}.
\end{equation}
\end{theorem}
\begin{proof}
We may assume type $d=0$. Let $\calA$ be as described in Definitions \ref{real39} and 
\ref{real313}. Then $(1-\omega)\calP(1-\omega_1)$ has the required mapping properties by the 
results of \cite{GrKo}, since away from the tip cone Sobolev and cone Besov spaces coincide 
with the usual ones. Also $\calG$ behaves correctly due to Theorem \ref{real310}. Thus we may 
assume that 
 $$\calA=\omega\,t^{-\mu}\,\Big(\op_M(h)+\smsum_{j=0}^{k-1}t^j\,
   \op_M^{\gamma_j-\frac{n}{2}}(f_j)\Big)\,\omega_0,$$
cf.\ \eqref{mellin}. In order to obtain \eqref{map1} it is then enough to consider 
 $$\wt\calA=\omega\,t^{-\mu}\,\op_M^{\gamma-\frac{n}{2}}(f)\,\omega_0$$
for some $f(t,\gamma+i\tau)\in\ci(\rpbar,B^{\mu,0}(X;\rz_\tau))$. We pull back $t^\mu\wt\calA$ 
to an operator on $\rz\times X$ using the maps $S_\gamma$ and $S_{\gamma-\frac{1}{2}}$ in 
\eqref{change} and \eqref{change2}. Since
 $$\begin{pmatrix}S_\gamma^{-1}&0\\0&S_\gamma^{\prime -1}\end{pmatrix}
   \op_M^{\gamma-\frac{n}{2}}(f)
   \begin{pmatrix}S_\gamma&0\\0&S_\gamma^\prime\end{pmatrix}=
   \op(f_\gamma),\qquad
   f_\gamma(r,\varrho)=f(e^r,\frac{n+1}{2}-\gamma+i\varrho),$$
this yields a `standard' operator in Boutet's calculus on $\rz\times X$ and the assertion 
follows from \cite{GrKo}. 

Let us turn to the proof of \eqref{map2}. From the known case $p=2$, cf.\ Theorem 4.1.14 in 
\cite{ScSc1} and Lemma 3.1.8 in \cite{ScSc2}, we have that 
 $$\calA:\begin{matrix}
      \calC^{\infty,\gamma}_{P}(\dz)\\ \oplus\\ \calC^{\infty,\gamma-\frac{1}{2}}_{Q}(\bz)
     \end{matrix}
     \longrightarrow
     \begin{matrix}
      \calC^{\infty,\gamma-\mu}_{P^\prime}(\dz)\\ \oplus\\ 
      \calC^{\infty,\gamma-\mu-\frac{1}{2}}_{Q^\prime}(\bz)
     \end{matrix}.$$
By the definition of the spaces with asymptotics, cf. Definition \ref{real35}, we thus may 
assume that $P$ and $Q$ are the empty asymptotic types. Since $h$ is holomorphic, 
$\op_M^\sigma(h)$ on $\cicomp(\intd)\oplus\cicomp(\intb)$ is independent of the choice of 
$\sigma$. By density and \eqref{map1} this yields that 
 $$\omega\,t^{-\mu}\,\op_M(h)\,\omega_0:\begin{matrix}
      \calH^{s,\gamma}_{p,k}(\dz)\\ \oplus\\
      \calB^{s-\frac{1}{p},\gamma-\frac{1}{2}}_{p,k}(\bz)
     \end{matrix}
     \longrightarrow
     \begin{matrix}
      \calH^{s-\mu,\gamma-\mu}_{p,k}(\dz)\\ \oplus\\ 
      \calB^{s-\mu-\frac{1}{p},\gamma-\mu-\frac{1}{2}}_{p,k}(\bz)
     \end{matrix}$$
is continuous. It remains to consider $\wt\calA$ as above, but now with a meromorphic Mellin 
symbol $f\in M^{-\infty,0}_R(X)$. For small $\eps>0$ we write 
 $$\wt\calA=\omega\,t^{-\mu}\,\op_M^{\gamma+k-\eps-\frac{n}{2}}(f)\,\omega_0+
   t^{-\mu}\,\omega\,
   \big(\op_M^{\gamma-\frac{n}{2}}(f)-\op_M^{\gamma+k-\eps-\frac{n}{2}}(f)\big)\,\omega_0.$$
If $\eps$ is sufficiently small $f$ will be holomorphic in the strip 
$-\gamma-k<\re z-\frac{n+1}{2}<-\gamma-k+\eps$, and thus the action of 
$\op_M^{\gamma+k-\eps-\frac{n}{2}}(f)$ will be independent of  $\eps>0$. Hence the first 
summand maps 
$\calH^{s,\gamma}_{p,k}(\dz)\oplus\calB^{s-\frac{1}{p},\gamma-\frac{1}{2}}_{p,k}(\bz)$ into 
$\calC^{\infty,\gamma-\mu}_k(\dz)\oplus\calC^{\infty,\gamma-\mu-\frac{1}{2}}_k(\bz)$ 
The second summand can be calculated explicitly by means of Cauchy's integral formula, see 
Proposition \ref{real60}, giving the desired result.  
\end{proof}

\begin{theorem}\label{real315.5}
Let $\calA\in C^{\mu,0}(\dz;\gamma,\gamma-\mu,k)$ with $\mu\le0$. Then the formal adjoint of 
$\calA$ (with respect to the scalar product of 
$\calH^{0,0}_2(\dz)\oplus\calB^{-\frac{1}{2},-\frac{1}{2}}_2(\bz)$) belongs to 
$C^{\mu,0}(\dz;-\gamma+\mu,-\gamma,k)$. If $\calA$ belongs to the $C_{M+G}$- or $C_G$-classes, 
the same is true for the formal adjoint.   
\end{theorem}

\subsection{Symbolic structure, ellipticity, and parametrix}\label{section3.4}

Let $\calA=\begin{pmatrix}A&K\\T&Q\end{pmatrix}\in C^{\mu,d}(\dz;\gamma,\gamma-\mu,k)$ be 
given. On $\intd$ (i.e.\ away from the conical singularity) the operator $\calA$ is a usual 
element of Boutet de Monvel's algebra. Hence we can associate with $\calA$ the usual 
(homogeneous) {\em principal symbol}
 $$\sigma_\psi^\mu(\calA)=\sigma_\psi^\mu(A)\in\ci(T^*\intd\setminus 0)$$
and the usual (principal) {\em boundary symbol} 
\begin{equation}\label{bdrysymb}
 \sigma_\partial^\mu(\calA)\in\ci(T^*\intb\setminus 0,
 \calL(H^s_p(\rz_+)\oplus\cz,H^{s-\mu}_p(\rz_+)\oplus\cz)).
\end{equation}
Near the conical singularity these symbols are `cone degenerate', i.e.\ 
the limits 
 $$\wt\sigma_\psi^\mu(\calA)(x,\tau,\xi)=
   \lim_{t\to0+}t^\mu\,\sigma_\psi^\mu(\calA)(t,x,t^{-1}\tau,\xi),\qquad 
 \wt\sigma_\partial^\mu(\calA)(x^\prime,\tau,\xi^\prime)=
   \lim_{t\to0+}t^\mu\,\sigma_\partial^\mu(\calA)(t,x^\prime,t^{-1}\tau,\xi^\prime)$$
exist. We call 
 $$\wt\sigma_\psi^\mu(\calA)\in\ci((T^*X\times\rz)\setminus 0)$$
the {\em rescaled principal symbol}. The {\em rescaled boundary symbol} is
\begin{equation}\label{rescbdrysymb} 
 \wt\sigma_\partial^\mu(\calA)\in\ci((T^*\partial X\times\rz)\setminus 0,
 \calL(H^s_p(\rz_+)\oplus\cz,H^{s-\mu}_p(\rz_+)\oplus\cz)).
\end{equation}
According to Definition \ref{real313} $\calA$ is of the form  
 $$\calA=\omega\,t^{-\mu}\,\op_M(h)\,\omega_0+(1-\omega)\calP(1-\omega_1)+\calM+\calG,$$
where $h$ and $\calP$ are as in Definition \ref{real39}.i,ii), and $\calM+\calG$ are as in 
definition \ref{real312}. Then, by definition, the {\em conormal symbol} of $\calA$ is 
 $$\sigma_M^\mu(\calA)(z)=h(0,z)+f_0(z).$$
It is a meromorphic function with values in $B^{\mu,d}(X)$ or, alternatively, 
\begin{equation}\label{consymb}
 \sigma_M^\mu(\calA)(z)\in\calL(H^s_p(X)\oplus B^{s-\frac{1}{p}}_p(\partial X),
 H^{s-\mu}_p(X)\oplus B^{s-\mu-\frac{1}{p}}_p(\partial X))
\end{equation}
for any $s>d-1+\frac{1}{p}$. Recall also that 
$\sigma_M^\mu(\calA)(\delta+i\tau)\in B^{\mu,d}(X;\rz_\tau)$ whenever $f_0$ has no pole on 
the line $\re z=\delta$. We then have the following compatibility relations: 
 $$\sigma_\psi^\mu(\sigma_M^\mu(\calA))(x,\xi,\tau)=
   \wt\sigma_\psi^\mu(\calA)(x,-\tau,\xi),\qquad
   \sigma_\partial^\mu(\sigma_M^\mu(\calA))(x^\prime,\xi^\prime,\tau)=
   \wt\sigma_\partial^\mu(\calA)(x^\prime,-\tau,\xi^\prime)$$
independently of $\delta\in\rz$ (note that $f_0$ is of order $-\infty$ and thus does not 
contribute to the principal symbolic levels). 

\begin{remark}\label{real316}
Let $\calA\in C^{\mu,d}(\dz;\gamma,\gamma-\mu,k)$ with $d=\max(0,\mu)$ and assume that both 
$\wt\sigma_\psi^\mu(\calA)$ and $\wt\sigma_\partial^\mu(\calA)$ are pointwise everywhere 
invertible. Then the before mentioned compatibility relations ensure that 
$h(0,\delta+i\tau)\in B^{\mu,d}(X;\rz_\tau)$ is (parameter-dependent) elliptic for each 
$\delta$, hence $\sigma_M^\mu(\calA)(\delta+i\tau)$ is invertible for large $|\tau|$. By a 
classical theorem on the invertibility of Fredholm families one then can deduce that 
$\sigma_M^\mu(\calA)$ is meromorphically invertible with 
 $$\sigma_M^\mu(\calA)^{-1}(z)=h^\prime(z)+f_0^\prime(z),$$
where $h^\prime\in M^{\mu,d^\prime}_O(X)$ and 
$f_0^\prime\in M^{-\infty,d^\prime}_{P^\prime}(X)$ with $d^\prime=\max(0,-\mu)$. For details 
see \cite{ScSc2}. 
\end{remark}

Note that the invertibilty of the (rescaled) boundary symbol \eqref{bdrysymb}, 
\eqref{rescbdrysymb}, and the conormal symbol \eqref{consymb} is independent of $s$ and $p$. 

\begin{definition}\label{real317}
We call $\calA\in C^{\mu,d}(\dz;\gamma,\gamma-\mu,k)$ {\em elliptic} if the following three 
conditions are satisfied: 
 \begin{itemize}
  \item[i)] Both $\sigma_\psi^\mu(\calA)$ and $\wt\sigma_\psi^\mu(\calA)$ are pointwise
   everywhere invertible, 
  \item[ii)] both $\sigma_\partial^\mu(\calA)$ and $\wt\sigma_\partial^\mu(\calA)$ are 
   pointwise everywhere invertible, 
  \item[iii)] $\sigma_M^\mu(\calA)^{-1}$ has no pole on the line $\re z=\frac{n+1}{2}-\gamma$. 
 \end{itemize}
\end{definition}

\begin{theorem}\label{real318}
Let $\calA\in C^{\mu,d}(\dz;\gamma,\gamma-\mu,k)$, $d=\max(0,\mu)$, be elliptic. Then there 
exists a $\calB \in C^{-\mu,d^\prime}(\dz;\gamma-\mu,\gamma,k)$, $d^\prime=\max(0,-\mu)$, such 
that 
 $$\calB\calA-1\in C^{d}_G(\dz;\gamma,\gamma,k),\qquad
   \calA\calB-1\in C^{d^\prime}_G(\dz;\gamma-\mu,\gamma-\mu,k).$$
\end{theorem}

Any such $\calB$ with these properties we call a {\em parametrix} of $\calA$. Using the 
compactness of Green operators and the mapping properties of elements of the cone algebra, one 
can easily deduce the following results on Fredholm property and elliptic regularity: 

\begin{theorem}\label{real319}
Let $\calA\in C^{\mu,d}(\dz;\gamma,\gamma-\mu,k)$, $d=\max(0,\mu)$, be elliptic. Then 
 $$\calA:\begin{matrix}
      \hsgpd\\ \oplus\\ \calB^{s-\frac{1}{p},\gamma-\frac{1}{2}}_p(\bz)
     \end{matrix}
     \longrightarrow
     \begin{matrix}
      \calH^{s-\mu,\gamma-\mu}_p(\dz)\\ \oplus\\ 
      \calB^{s-\mu-\frac{1}{p},\gamma-\mu-\frac{1}{2}}_p(\bz)
     \end{matrix}$$
is a Fredholm operator for any $1<p<\infty$ and $s>d-1+\frac{1}{p}$. 
\end{theorem}

\begin{theorem}\label{real320}
Let $\calA\in C^{\mu,d}(\dz;\gamma,\gamma-\mu,k)$, $d=\max(0,\mu)$, be elliptic and consider 
the equation $Au=f$ Then: 
 \begin{itemize}
  \item[a)] If $f\in\calH^{s-\mu,\gamma-\mu}_p(\dz)\oplus 
   \calB^{s-\mu-\frac{1}{p},\gamma-\mu-\frac{1}{2}}_p(\bz)$, $s>d-1+\frac{1}{p}$, and 
   $u\in\calH^{t,\gamma}_p(\dz)\oplus 
   \calB^{t-\frac{1}{p},\gamma-\frac{1}{2}}_p(\bz)$ for some $t>\max(0,-\mu)-1+\frac{1}{p}$, 
   then $u\in\calH^{s,\gamma}_p(\dz)\oplus 
   \calB^{s-\frac{1}{p},\gamma-\frac{1}{2}}_p(\bz)$. 
  \item[b)] If $f\in\calH^{s-\mu,\gamma-\mu}_{p,P^\prime}(\dz)\oplus 
   \calB^{s-\mu-\frac{1}{p},\gamma-\mu-\frac{1}{2}}_{p,Q^\prime}(\bz)$, $s>d-1+\frac{1}{p}$, 
   for some asymptotic type $(P^\prime,Q^\prime)\in\as(X,\partial X;\gamma-\mu,k)$ and 
   $u\in\calH^{t,\gamma}_p(\dz)\oplus 
   \calB^{t-\frac{1}{p},\gamma-\frac{1}{2}}_p(\bz)$ for some $t>\max(0,-\mu)-1+\frac{1}{p}$, 
   then there exists a $(P,Q)\in\as(X,\partial X;\gamma,k)$ such that 
   $u\in\calH^{s,\gamma}_{p,P}(\dz)\oplus 
   \calB^{s-\frac{1}{p},\gamma-\frac{1}{2}}_{p,Q}(\bz)$. 
 \end{itemize}
\end{theorem}

\subsection{Operators acting on sections into vector bundles}\label{section3.5}
For simplicity of the presentation, we so far restricted ourselves to operators acting on 
scalar valued functions. However, the results in Sections \ref{section3.1} to \ref{section3.4} 
readily extend to operators acting on sections of vector bundles. 

If $E,\wt E$ and $F,\wt F$ are vector bundles over $\dz$ and $\bz$ (see the introduction), 
respectively, we can speak of the cone algebra 
$C^{\mu,d}(\dz;\gamma,\gamma-\mu,k;E,F;\wt E,\wt F)$. Elements $\calA$ of this space then 
induce operators 
 $$\calA:\begin{matrix}
      \calH^{s,\gamma}_p(\dz,E)\\ \oplus\\ 
      \calB^{s-\frac{1}{p},\gamma-\frac{1}{2}}_p(\bz,F)
     \end{matrix}
     \longrightarrow
     \begin{matrix}
      \calH^{s-\mu,\gamma-\mu}_p(\dz,\wt E)\\ \oplus\\ 
      \calB^{s-\mu-\frac{1}{p},\gamma-\mu-\frac{1}{2}}_p(\bz,\wt F)
     \end{matrix}$$
for any $1<p<\infty$ and $s>d-1+\frac{1}{p}$. It is straightforward to adapt the previous 
definitions to the vector bundle situation. The holomorphic Mellin symbol, for example, takes 
values in $B^{\mu,d}(X;E_0,F_0;\wt E_0,\wt F_0)$, where subscript $0$ indicates restriction of 
the bundles to $t=0$. 

The principal and the rescaled principal symbols are then homomorphisms acting between the 
pull-backs of $E$ and $\wt E$ to the cotangent bundle. Similar statements for the (rescaled) 
boundary symbol. More details will be given in Section \ref{section4.2}.   

If the bundles $F$ or $\wt F$ are zero dimensional, we shall omit them from the notation, 
e.g.\ we write $C^{\mu,d}(\dz;\gamma,\gamma-\mu,k;E;\wt E)$ if $F=\wt F=0$. 


\section{Differential boundary value problems on 
         manifolds with conical singularities}\label{section4}

In the following let $A$ be a Fuchs type differential operator of order $\mu\in\nz$ on $\dz$, 
acting between sections of a vector bundle $E$ over $\dz$. 
Using the notation from \eqref{real4A.5}, the principal symbol of $A$ near $t=0$ is 
 $$\sigma_\psi^\mu(A)(t,x,\tau,\xi)=t^{-\mu}\smsum_{j=0}^\mu
   \sigma_\psi^{\mu-j}(a_j)(t,x,\xi)(-it\tau)^j,$$
and thus the rescaled principal symbol is 
 $$\wt\sigma_\psi^\mu(A)(x,\tau,\xi)=\smsum_{j=0}^\mu
   \sigma_\psi^{\mu-j}(a_j)(0,x,\xi)(-i\tau)^j.$$
We shall say that $A$ has {\em $t$-independent coefficients} (near the singularity $t=0$), if 
all the $a_j$ are constant in $t$ near $t=0$. 

\subsection{Differential boundary conditions}\label{section4.1}
A typical boundary condition for $A$ is of the form $\gamma_0 B$, where $\gamma_0$ denotes 
the operator of restriction to $\bz$, and $B$ is a Fuchs type differential operator on $\dz$, 
acting from sections of $E$ to sections of a bundle $\wt F$ over $\dz$. Near $t=0$ let  
 $$B=t^{-\nu}\smsum_{j=0}^\nu b_j(t)(-t\partial_t)^j,\qquad
   b_j\in\ci([0,1[,\text{\rm Diff}^{\nu-j}(X;E_0,\wt F_0)).$$
Using the splitting $x=(x^\prime,x_n)$ near the boundary $\partial X$, we write 
 $$b_j(t)=\smsum_{k=0}^{\nu-j} b_{jk}(t)D_{x_n}^k,\qquad
   b_{jk}\in\ci([0,1[,\text{\rm Diff}^{\nu-j-k}(\partial X;E_0^\prime,\wt F_0^\prime)),$$
where $^\prime$ indicates restriction of the bundles to $\partial X$. Then   
\begin{equation}\label{real4A}
 \gamma_0 B=\smsum_{k=0}^{\nu}t^{-k}\Big(t^{-\nu+k}\smsum_{j=0}^{\nu-k}
 b_{jk}(t)(-t\partial_t)^j\Big)\gamma_k=
 \smsum_{k=0}^{\nu}\ulS_k\,t^{-k}\,\gamma_k,
\end{equation}
with $\gamma_k=\gamma_0 D_n^k$ and Fuchs type differential operators $\ulS_k$ of order $\nu-k$ 
on $\bz$ (acting from $E_\bz$ to $\wt F_\bz$, the restrictions of the bundles to $\bz$). 

Given $\mu$ such boundary conditions
 $$T_k=\gamma_0 B_k,\qquad k=0,\ldots,\mu-1,$$
with Fuchs type differential operators $B_k$ of order $k$ acting from sections of $E$ to 
sections of $\wt F_k$  and if we set 
\begin{equation}\label{real4B}
 T=\begin{pmatrix}T_0&T_1&\ldots&T_{\mu-1}\end{pmatrix}^t,\qquad
 \varrho= \begin{pmatrix}\gamma_0&\gamma_1&\ldots&\gamma_{\mu-1}\end{pmatrix}^t,
\end{equation}
we obtain $T=S\varrho$ with a left-lower triangular matrix 
\begin{equation}\label{real4C}
 S=(\ulS_{jk})_{0\le j,k<\mu}\;\text{\rm diag}(1,t^{-1},\ldots,t^{-\mu+1}).
\end{equation}
Here the $\ulS_{jk}$ are Fuchs type operators of order $j-k$ on $\bz$, and $\ulS_{jk}=0$ if 
$k>j$. Each $\ulS_{jk}$ acts from sections of $E_\bz$ to sections of $F_k:=\wt F_{k\bz}$. 

For the remaining part of the paper $T=S\varrho$ will be the standard form of a boundary 
condition for $A$. Note that we allow $F_k$ to be zero dimensional. In this case the $k$-th 
boundary condition $T_k$ is void. Nevertheless, for systematic reasons, it is convenient to 
work with the full matrix $S$. 

We shall say that $T$ has {\em $t$-independent coefficients} (near the singularity $t=0$), 
if all the $\ulS_{jk}$ have $t$-independent coefficients. 

\subsection{Relation with the cone algebra}\label{section4.2}
Let $T=S\varrho$ be as above and consider the boundary value problem 
\begin{equation}\label{bvp}
 \calA=\begin{pmatrix}A\\T\end{pmatrix}:\calH^{s,\gamma}_p(\dz,E)\longrightarrow
 \begin{matrix}
  \calH^{s-\mu,\gamma-\mu}_p(\dz,E)\\ \oplus\\
  \mathop{\oplus}\limits_{k=0}^{\mu-1}\calB^{s-k-\frac{1}{p},\gamma-k-\frac{1}{2}}(\bz,F_k)
 \end{matrix}
\end{equation}
for $s>\mu-1+\frac{1}{p}$. We choose order reductions on the boundary
 $$\Lambda_k\in C^{\mu-k}(\bz;\gamma-k,\gamma-\mu,\theta;F_k,F_k)$$
where $\theta\in\nz$ can be chosen arbitrarily large, and let 
$\Lambda:=\diag(\Lambda_0,\ldots,\Lambda_{\mu-1})$. In particular,  
 $$\Lambda:\mathop{\oplus}\limits_{k=0}^{\mu-1}
   \calB^{s-k-\frac{1}{p},\gamma-k-\frac{1}{2}}(\bz,F_k)
   \longrightarrow
   \mathop{\oplus}\limits_{k=0}^{\mu-1}
   \calB^{s-\mu-\frac{1}{p},\gamma-\mu-\frac{1}{2}}(\bz,F_k)=
   \calB^{s-\mu-\frac{1}{p},\gamma-\mu-\frac{1}{2}}(\bz,F)$$
for $F:=F_0\oplus\ldots\oplus F_{\mu-1}$ is an isomorphism. Then 
 $$\calA_\Lambda:=\begin{pmatrix}1&0\\0&\Lambda\end{pmatrix}
   \begin{pmatrix}A\\T\end{pmatrix}=
   \begin{pmatrix}A\\\Lambda T\end{pmatrix}:
   \calH^{s,\gamma}_p(\dz,E)\longrightarrow
   \begin{matrix}
    \calH^{s-\mu,\gamma-\mu}_p(\dz,E)\\ \oplus\\
    \calB^{s-\mu-\frac{1}{p},\gamma-\mu-\frac{1}{2}}(\bz,F)
   \end{matrix}$$
is an element of the cone algebra as it has been described in Section \ref{section3}, i.e.\ 
$\calA_\Lambda\in C^{\mu,\mu}(\dz;\gamma,\gamma-\mu,\theta;E;E,F)$. 

\begin{definition}\label{real41}
The boundary value problem $\calA$ is called {\em elliptic with respect to} $\gamma$ if the 
resulting $\calA_\Lambda$ is elliptic in the sense of Definition \mbox{\rm\ref{real317}}. 
\end{definition}

For later application it is convenient to describe ellipticity of $\calA$ intrinsically and 
not refering to $\calA_\Lambda$. To this end we introduce a symbolic structure for $\calA$, 
using the representation of $T=S\varrho$ as in \eqref{real4C}. 

The principal and rescaled principal symbol of $\calA$ are 
 $$\sigma_\psi^\mu(\calA):=\sigma_\psi^\mu(A),\qquad
   \wt\sigma_\psi^\mu(\calA):=\wt\sigma_\psi^\mu(A).$$
The boundary symbol 
 $$\sigma_\partial^\mu(\calA):\pi_\bz^*E_\bz\otimes H^s_p(\rz_+)\longrightarrow\pi_\bz^*
   \begin{pmatrix}
    E_\bz\otimes H^{s-\mu}_p(\rz_+)\\ \oplus\\ F
   \end{pmatrix}$$
is given by 
 $$\sigma_\partial^\mu(\calA)(y,\eta)=
   \begin{pmatrix}
    \sigma_\psi^\mu(A)(y,0,\eta,D_s)\\
    \left( \sigma_\psi^{j-k}(\ulS_{jk}(y,\eta)\right)_{0\le j,k<\mu}\,\varrho
   \end{pmatrix}.$$
Here, $\pi_\bz:T^*\intb\setminus0\to\intb$ is the canonical projection, $(y,s)$ refers to a 
splitting of coordinates in a collar neighborhood of $\bz$, and 
$\varrho=(\gamma_0,\ldots,\gamma_{\mu-1})_s$ acts on the half-axis $\rz_{+}$. Similarly, with 
the projection $\pi_{\partial X}:T^*\partial X\setminus0\to\partial X$, we define 
 $$\wt\sigma_\partial^\mu(\calA):\pi_{\partial X}^*E_0^\prime\otimes H^s_p(\rz_+)
   \longrightarrow\pi_{\partial X}^*
   \begin{pmatrix}
    E_0^\prime\otimes H^{s-\mu}_p(\rz_+)\\ \oplus\\ F_0
   \end{pmatrix}$$
by  
 $$\wt\sigma_\partial^\mu(\calA)(x^\prime,\xi^\prime,\tau)=
   \begin{pmatrix}
    \wt\sigma_\psi^\mu(A)(x^\prime,0,\xi^\prime,\tau,D_{x_n})\\
    \left(\wt\sigma_\psi^{j-k}(\ulS_{jk}(x^\prime,\xi^\prime,\tau)\right)_{0\le j,k<\mu}\,
    \varrho
   \end{pmatrix}.$$
Finally, the conormal symbol of $\calA$,
 $$\sigma_M^\mu(\calA)(z):
   H^s_p(X,E_0)\longrightarrow
   \begin{matrix}
     H^{s-\mu}_p(X,E_0)\\ \oplus\\ 
     \mathop{\oplus}\limits_{k=0}^{\mu-1}B^{s-k-\frac{1}{p}}_p(\partial X,F_{k0})
   \end{matrix}$$
for $s>\mu-1+\frac{1}{p}$ is defined by 
 $$\sigma_M^\mu(\calA)(z)=
   \begin{pmatrix}
    \sigma_M^\mu(A)(z)\\
    \left(T^k\sigma_M^{j-k}(\ulS_{jk})(z)\right)_{0\le j,k<\mu}\,\varrho_{x_n}
   \end{pmatrix}.$$
Here, $T^k$ acts on functions as operator of translation by $k$, i.e.\ $(T^kf)(z)=f(z+k)$. 
Note that this shift appears here, since we commuted in \eqref{real4A} the factor $t^{-k}$ 
to the right. 

A straightforward calculation then shows that 
 $$\sigma^\mu(\calA_\Lambda)=
   \begin{pmatrix}1&0\\0&\sigma(\Lambda)\end{pmatrix}\sigma^\mu(\calA),$$
where $\sigma$ is any of the symbols $\sigma_\psi$, $\sigma^{}_\partial$, $\wt\sigma_\psi$, 
$\wt\sigma^{}_\partial$, or $\sigma^{}_M$ and 
 $$\sigma(\Lambda):=\diag(\sigma^\mu(\Lambda_0),\sigma^{\mu-1}(\Lambda_1),\ldots,
   \sigma^1(\Lambda_{\mu-1}))$$
in all of the first four cases, while for $\sigma=\sigma^{}_M$ we have 
 $$\sigma^{}_M(\Lambda):=\diag(T^0\sigma_M^\mu(\Lambda_0),
   T^1\sigma_M^{\mu-1}(\Lambda_1),\ldots,T^{\mu-1}\sigma_M^1(\Lambda_{\mu-1})).$$
From these relations the following statement is immediately clear. 

\begin{proposition}\label{real42}
$\calA$ of \textrm{\eqref{bvp}} is elliptic with respect to $\gamma$ if and only if 
 \begin{itemize}
  \item[i)] Both $\sigma_\psi^\mu(\calA)$ and $\wt\sigma_\psi^\mu(\calA)$ are invertible, 
  \item[ii)] both $\sigma_\partial^\mu(\calA)$ and $\wt\sigma_\partial^\mu(\calA)$ are 
   invertible, 
  \item[iii)] $\sigma_M^\mu(\calA)(z)$ is invertible for each $z$ with 
   $\re z=\frac{n+1}{2}-\gamma$. 
 \end{itemize}
\end{proposition}

\subsection{A parametrix construction}

Let $\calA=\begin{pmatrix}A\\T\end{pmatrix}$ be elliptic with respect to $\gamma$. 
Then there exists a parametrix $\calB_\Lambda$ to $\calA_\Lambda$ in the sense of Theorem 
\ref{real318} for any arbitrarily large prescribed $k\in\nz$. 
We call 
 $$\begin{pmatrix}R&K\end{pmatrix}:=\calB_\Lambda
   \begin{pmatrix}1&0\\0&\Lambda^{-1}\end{pmatrix}$$
a parametrix of the original problem $\calA$. We then have that 
\begin{equation}\label{real4D}
 R\in C^{-\mu,0}(\dz;\gamma-\mu,\gamma,k;E;E)
\end{equation}
and, modulo Green remainders, 
 $$RA+KT\equiv 1,\qquad 
   \begin{pmatrix}AR&AK\\TR&TK\end{pmatrix}\equiv
   \begin{pmatrix}1&0\\0&1\end{pmatrix}.$$
This construction can be improved. For notational convenience let us introduce the following 
abbreviation: Whenever $\calF$ is a function or distribution space on $\dz$ on which the 
boundary operator $T$ acts continuously, we set 
 $$\calF_T=\{u\in\calF\st Tu=0\}.$$
This is then a closed subspace of $\calF$. 

\begin{proposition}\label{real43}
Let $\begin{pmatrix}A\\T\end{pmatrix}$ be elliptic with respect to the
weight $\gamma$. Then there exists a parametrix 
$\begin{pmatrix}R_0&K_0\end{pmatrix}$ such that 
\begin{itemize}
 \item[a)] $TR_0=0$, i.e.\ $R_0$ maps $\calH^{s-\mu,\gamma-\mu}_p(\dz,E)$, $s>-1+\frac{1}{p}$, 
  into the kernel of $T$ in $\calH_p^{s,\gamma}(\dz,E)$. We denote this space by
  $\calH_p^{s,\gamma}(\dz,E)_T$.
 \item[b)]  $AR_0-1\in C^\mu_G(\dz;\gamma-\mu,\gamma-\mu,k;E;E)$ is a finite rank operator 
  on $\calH_p^{s,\gamma}(\dz,E)$ with image in 
  $\calC_P^{\infty,\gamma-\mu}(\dz,E)$ for a suitable asymptotic type $P$.
 \item [c)] On $\calH_p^{s,\gamma}(\dz,E)_T$,
  the operator $R_0A-1$ has finite rank; 
  it coincides with an element in $C^0_G(\dz;\gamma,\gamma,k;E;E)$ and has its image in 
  $\calC_{P^\prime}^{\infty,\gamma}(\dz,E) $ for a suitable asymptotic type  $P^\prime$.
\end{itemize}
\end{proposition}
\begin{proof}
Our construction is based on an argument by Grubb \cite{Grub1}, Proposition 1.4.2.
We first reduce order and weight to zero: Let $\Lambda$ be the operator constructed in Section 
\ref{section4.2}, and denote by $\Lambda_-^{-\mu}$  an invertible cone differential 
operator of order $-\mu$ , cf.\ Theorem 2.10 in \cite{HaSc03}. 
For $s>-1+1/p$ we consider the operator 
\begin{eqnarray}\label{0}
\widetilde \calA= \begin{pmatrix}A\\ \Lambda T\end{pmatrix}\Lambda_-^{-\mu} :
\calH_p^{s,\gamma-\mu}(\dz,E)
\longrightarrow
\begin{matrix} 
\calH_p^{s,\gamma-\mu}(\dz,E)\\ \oplus \\ \calB^{s-1/p, \gamma-\mu-1/2}_p(\bz,F)
\end{matrix}.
\end{eqnarray}
By conjugating with $t^{\gamma-\mu}$ we may also assume that $\gamma=\mu$. 

As a composition of an elliptic operator with an invertible operator, 
$\widetilde\calA\in C^{0,0}(\dz,0,0,k;E;E,F)$ is elliptic of order and type zero.
Hence  there exists a parametrix $ \calC$ 
of order and type zero such that 
\begin{equation}\label{-1}
 \wt\calA{\calC}=1+\calS\quad\text{and}\quad  \calC\widetilde\calA=1+\calS^\prime
\end{equation}
with Green operators $\calS\in  C_G^0(\dz;0,0,k;E,F;E,F)$ and $\calS^\prime\in 
C_G^0(\dz;0,0,k;E;E)$, respectively.

We denote by $\widetilde{\calA}^*$ the formal adjoint of $\widetilde\calA$ 
with respect to the scalar products in  
$\calH^{0,0}_2(\dz,E)$ and $\calH^{0,0}_2(\dz,E) \oplus  \calB^{-1/2, -1/2}_2(\bz,F)$, 
respectively. 
$\widetilde{\calA}^*$ also is elliptic of order and type zero; for $1<q<\infty$ and 
$s^\prime>-1+1/q$ it extends to a Fredholm operator 
 $$\widetilde{\calA}^*: \begin{matrix} 
   \calH_q^{s^\prime,0}(\dz,E)\\ \oplus\\ \calB^{s^\prime-1/q, -1/2}_q(\bz,F)
   \end{matrix}\longrightarrow \; \calH_q^{s^\prime,0}(\dz,E).$$
By elliptic regularity, the kernel is finite dimensional and independent of $q$ and $s^\prime$; 
it consists of functions in 
$\calC_P^{\infty,0}(\dz,E)\oplus \calC_Q^{\infty,-1/2}(\bz,F)$ 
for a suitable asymptotic type $(P,Q)$. We choose a basis $\{v_1,\ldots,v_m\}$ which is 
orthonormal with respect to the scalar product in 
$\calH^{0,0}_2(\dz,E)\oplus \calB^{-1/2,-1/2}_2(\bz,F)$.

Let us show that ${\rm ker}\,\wt{\calA}^*={\rm span}\,\{v_1,\ldots,v_m\}$ is a complement of 
the image of $\widetilde\calA $ in $\calH_p^{s,0}(\dz,E) \oplus \calB^{s-1/p, -1/2}_p(\bz,F)$. 
In fact, for $u\in \calH^{s,0}_p(\dz,E)$,  
\begin{equation}\label{1} 
 \langle \widetilde \calA u,v_j\rangle _{\calH^{s,0}_p\oplus  \calB^{s-1/p, -1/2}_p, 
 \overset{\circ}{\calH}{}^{-s,0}_{p^\prime}\oplus  \calB^{-s-1/p^\prime, -1/2}_{p^\prime}}=
 \langle u, \wt{\calA}^*v_j
 \rangle_{\calH^{s,0}_p,\overset{\circ}{\calH}{}^{-s,0}_{p^\prime}}=0,\qquad j=1,\ldots,m.
\end{equation}
For $v$ in ${\rm ker}\,\wt{\calA}^*\cap {\rm im}\, \widetilde \calA$, \eqref{1} implies that 
 $$0=\langle v,v\rangle_{\calH^{s,0}_p\oplus  \calB^{s-1/p, -1/2}_p, 
   \overset{\circ}{\calH}{}^{-s,0}_{p^\prime}\oplus  \calB^{-s-1/p^\prime, -1/2}_{p^\prime}} 
   = \|v\|^2_{\calH^{0,0}_2\oplus  \calB^{-1/2, -1/2}_2}, $$
so that $v=0$. On the other hand, for an element $v$ in 
$\kringel{\calH}^{-s,0}_{p^\prime}(\dz,E)\oplus  \calB^{-s-1/p^\prime, -1/2}_{p^\prime}(\bz,F)$
we have 
 $$\langle v,\widetilde\calA u\rangle_{\overset{\circ}{\calH}{}^{-s,0}_{p^\prime}\oplus 
   \calB^{-s-1/p^\prime, -1/2}_{p^\prime}, 
   \calH^{s,0}_{p}\oplus  \calB^{s-1/p, -1/2}_{p}}=0 \qquad 
   u\in\calH^{s,0}_p(\dz,E)$$
if and only if $v$ is in the kernel of 
$\widetilde{\calA}^*$ on $\kringel{\calH}^{-s,0}_{p^\prime}(\dz,E)\oplus 
\calB^{-s-1/p^\prime, -1/2}_{p^\prime}$. Since the image of $\widetilde\calA$  is closed, the 
quotient $\calH^{s,0}_p(\dz,E)\oplus  \calB^{s-1/p, -1/2}_p(\bz,F)/
{\rm im}\, \widetilde\calA$ is isomorphic to the space of all functionals in 
the dual of $\calH^{s,0}_p(\dz,E)\oplus  \calB^{s-1/p, -1/2}_p(\bz,F)$ 
which vanish on ${\rm im }\,\widetilde\calA$, thus to 
${\ker\widetilde \calA^*}$. Hence ${\rm ker}\,\wt{\calA}^*$ is a complement to 
${\rm  im }\,\widetilde\calA$, for it has the right dimension.

As a consequence, we can write the projection $\pi_{{\rm ker}\,\wt{\calA}^*}$ onto 
${\rm ker}\,\wt{\calA}^*$  along the image of $\widetilde \calA$, 
 $$\pi_{{\rm ker}\,\wt{\calA}^*}:
   \begin{matrix} 
    \calH_p^{s,0}(\dz,E)\\ \oplus \\ \calB^{s-1/p, -1/2}_p(\bz,F)
   \end{matrix}
   \longrightarrow {\rm ker}\,\wt{\calA}^*\subset
   \begin{matrix}
    \calH_p^{s,0}(\dz,E)\\ \oplus \\ \calB^{s-1/p, -1/2}_p(\bz,F)
   \end{matrix},$$ 
in the form 
\begin{equation}\label{1a} 
 \pi_{{\rm ker}\,\wt{\calA}^*}(v)= \sum_{j=1}^m 
 \langle v,v_j\rangle_{\calH^{s,0}_p\oplus  \calB^{s-1/p, -1/2}_p, 
 \calH^{-s,0}_{p^\prime}\oplus  \calB^{-s-1/p^\prime, -1/2}_{p^\prime}}\ v_j.
\end{equation}
This shows that it is an element of  $C^0_G(\dz;0,0,k;E,F;E,F)$ with integral kernel 
$\sum_{j=1}^m v_j\otimes \overline v_j$, cf. Section \ref{section3.3}. 

In an analogous way we can choose a basis $\{u_1,\ldots,u_k\}$ of the kernel of $\wt\calA$ 
in $\calH_p^{s,0}(\dz,E)$. It is a subset of
$\calC^{\infty,0}_{P^\prime}(\dz,E)$ for some asymptotic type $P^\prime$, 
thus independent of $s$ and $p$; we choose it orthonormal with respect to 
the scalar product of $\calH^{0,0}_2(\dz,E)$. The projection $\pi_{{\rm ker}\,\wt{\calA}}$ 
onto ${\rm ker}\,\wt{\calA}$ along  the range of $\widetilde{\calA}^*:
\calH_p^{s,0}(\dz,E) \oplus \ \calB^{s-1/p, -1/2}_p(\bz,F)\longrightarrow 
\calH^{s,0}_p(\dz,E)$ is given by
$\pi_{{\rm ker}\,\wt{\calA}}(u)=\sum_{j=1}^k\langle u,u_j\rangle_{\calH^{s,0}_p,
\overset{\circ}{\calH}{}^{-s,0}_{p^\prime}}u_j $ 
and therefore an element of $C_G^0(\dz;0,0,k;E;E)$.

The kernel of $\pi_{{\rm ker}\,\wt{\calA}}$ is the range of $\wt{\calA}^*$ 
and thus a complement of the kernel of $\wt\calA$. 
Let $\calC^\prime$ denote the operator which acts like an inverse of 
 $$\wt{\calA}:{\rm ker}\,\pi_{{\rm ker}\,\wt{\calA}} \longrightarrow {\rm im}\, 
   \wt\calA={\rm ker}\,\pi_{{\rm ker}\,\wt{\calA}^*}$$
and is zero on ${\rm ker}\,\wt{\calA}^*$.  In fact, the operator $\calC^\prime$ is defined for 
all $1<p<\infty$ and  $s>1/p-1$; its action is independent of $s$ and $p$; 
since this is the case for $\widetilde \calA$. We have    
\begin{align}
 \calC^\prime\wt\calA&=1-\pi_{{\rm ker}\,\wt{\calA}}\quad\text{on } 
 \calH^{s,0}_p(\dz,E)\label{2}\\
 \widetilde \calA\calC^\prime&=1-\pi_{{\rm ker}\,\wt{\calA}^*}\quad\text{on }
 \calH^{s,0}_p(\dz,E)\oplus\calB^{s-1/p, -1/2}_p(\bz,F).\label{3}
\end{align}
For the difference of  $\calC^\prime$ and the above parametrix $\calC$ we have by \eqref{-1}:
 $$\calC^\prime-\calC
   =(\calC\wt\calA-\calS^\prime)(\calC^\prime-\calC)= 
    \calC(1-\pi_{{\rm ker}\,\wt{\calA}^*}) -\calC(1+\calS)-\calS^\prime(\calC^\prime-\calC)
   =-\calC(\pi_{{\rm ker}\,\wt{\calA}^*}+\calS)-\calS^\prime(\calC^\prime-\calC).$$
This operator maps $\calH^{s,0}_p(\dz,E)\oplus\calB^{s-1/p, -1/2}_p(\bz,F)$ to 
$\calC_{P_1}^{\infty,0}(\dz,E)\oplus\calC_{Q_1}^{\infty,-1/2}(\bz,F)$ 
with suitable asymptotic types. For the adjoint we have 
\begin{align*}
 (\calC^\prime-\calC)^*&=(\calC^*\widetilde\calA^*-\calS^*)({\calC^\prime}^{*}-\calC^*)= 
 \calC^*(1-\pi^*_{{\rm ker}\,\wt{\calA}})-
 \calC^*(1+{\calS^\prime}^{*})-\calS^*(\calC^\prime-\calC)^*\\
 &=-\calC^*(\pi^{*}_{{\rm ker}\,\wt{\calA}}+{\calS^\prime}^{*})-\calS^*(\calC^\prime-\calC)^*.
\end{align*}
Again, this operator maps $\calH^{s,0}_p(\dz,E)\oplus\calB^{s-1/p, -1/2}_p(\bz,F)$ to 
$\calC_{P_2}^{\infty,0}(\dz,E)\oplus\calC_{Q_2}^{\infty,-1/2}(\bz,F)$ 
with suitable asymptotic types and thus is a Green operator of type zero. We write 
$\calC=\begin{pmatrix}R&  K\end{pmatrix}$,  
$\calC^\prime=\begin{pmatrix}R^\prime&K^\prime\end{pmatrix}$, and  
$\widetilde\calA=\begin{pmatrix}\widetilde A\\\widetilde T \end{pmatrix}$ 
with $\widetilde A=A\Lambda_-^{-\mu}$ and $\widetilde T=\Lambda T \Lambda_-^{-\mu}$. 
Then there are Green operators $S_1$ and $S_2$ in $C_G^0(\dz;0,0,k;E;E)$  
such that 
 $$R^\prime\widetilde A=1+S_2\quad \text{on }\calH^{s,0}_p(\dz,E)_{\widetilde T},\qquad
    \widetilde A R^\prime=1+S_1\quad\text{on }\calH^{s,0}_p(\dz,E).$$
Indeed, the first identity follows from \eqref2, since, on 
$\calH^{s,0}_p(\dz,E)_{\wt{T}}$, we have 
 $$R^\prime\widetilde A= R^\prime\widetilde A+K^\prime\widetilde T 
   =\calC^\prime\widetilde\calA= 1-\pi_{{\rm ker}\,\wt{\calA}},$$  
while the second follows from \eqref{3}. Next we write 
 $$\pi_{{\rm ker}\,\wt{\calA}^*}=\begin{pmatrix}S_{11}&S_{12}\\S_{21}&S_{22}\end{pmatrix}: 
   \begin{matrix} 
   \calH_p^{s,0}(\dz,E)\\ \oplus \\ \calB^{s-1/p, -1/2}_p(\bz,F)
   \end{matrix}
   \longrightarrow
   \begin{matrix} 
   \calH_p^{s,0}(\dz,E)\\ \oplus \\ \calB^{s-1/p, -1/2}_p(\bz,F)
   \end{matrix}.$$
We know that $S_{21}$ is of finite rank. The adjoint, $S_{21}^*$, also is of finite rank. 
By \eqref{1a}, its range lies in $\calC_{Q^\prime}^{\infty,-1/2}(\dz,E)$, for a suitable 
asymptotic type and is independent of $p$. 
We denote this range by $Z$ and note that the projection $\pi_{Z}$ onto 
$Z$ along the kernel of $S_{21}$ also is an element of $C_G^0(\dz;0,0,k;E,E)$. 
In fact, choosing a basis $\{e_1,\ldots e_r\}$ of ${Z}$, which is orthonormal with 
respect to the $\calH^{0,0}_2(\dz,E)$ scalar product, we can write $\pi_{{Z}}$ 
as the integral operator with the kernel 
$\sum_ {j=1}^re_j\otimes \overline e_j$ with respect 
to the pairing between $\calH_p^{s,0}(\dz,E)$ and $\kringel{\calH}_{p^\prime}^{-s,0}(\dz,E)$. 
We now let
 $$R^{\prime\prime}= R^\prime(1-\pi_{Z}).$$
We infer from \eqref3 that 
 $$\widetilde \calA R^{\prime\prime}=
   \begin{pmatrix}\wt AR^{\prime\prime}\\ \wt T R^{\prime\prime}\end{pmatrix} 
   = \begin{pmatrix}1-S_{11}\\-S_{21}\end{pmatrix}(1-\pi_{Z})
   =\begin{pmatrix}(1-S_{11})(1-\pi_{Z})\\0\end{pmatrix},$$
so that $R^{\prime\prime}$ maps $\calH^{0,0}_p(\dz,E)$ to 
$\calH^{0,0}_p(\dz,E)_{\wt T}$, while $\wt AR^{\prime\prime}$ differs from the identity by a 
finite rank operator, say $S^{\prime\prime}$, in $C_G^0(\dz;0,0,k;E;E)$. 

We can now conclude the proof: We let 
$R_0=\Lambda_-^{-\mu}R^\prime(1-\pi_{Z})=\Lambda_-^{-\mu}R^{\prime\prime}$. 
Since $TR_0=\Lambda^{-1}\widetilde TR^{\prime\prime}=0$, this operator 
maps $\calH_p^{s,0}(\dz,E)$ into $\calH^{s+\mu,0}_p(\dz,E)_T$, and 
we obtain assertion a). Moreover,  
$AR_0-1=\wt AR^{\prime\prime}-1=S^{\prime\prime}\in C_G^0(\dz;0,0,k;E;E)$, showing b). 
Finally, on $\calH^{s,0}(\dz,E)_T$, 
\begin{align*}
 R_0A-1&=\Lambda^{-\mu}_- R^\prime(1-\pi_{Z})\widetilde A\Lambda^{\mu}_- -1
   =\Lambda^{-\mu}_- (R^\prime\widetilde A+K^\prime\widetilde T)\Lambda^{\mu}_- -1-
   \Lambda^{-\mu}_-R^\prime\pi_{Z}\Lambda^{\mu}_-\\
 &=\Lambda^{-\mu}_- \calC^\prime\widetilde\calA\Lambda^{\mu}_- -1-
   \Lambda^{-\mu}_-R^\prime\pi_{Z}\Lambda^{\mu}_-=
   -\Lambda^{-\mu}_-(\pi_{{\rm ker}\,\wt{\calA}}+R^\prime\pi_{Z}) \Lambda^{\mu}_-. 
\end{align*}
Hence $R_0A-1\in C_G(\dz;\mu,\mu,k;E,E)$, which implies  c).
\end{proof}

\subsection{Normal boundary conditions}

\begin{lemma}\label{real44}
There exists a map $K:\mathop{\oplus}\limits_{k=0}^{\mu-1}\cii(\bz,F_k)\to\cii(\dz,E)$ that 
extends continuously to 
 $$K:\mathop{\oplus}\limits_{k=0}^{\mu-1}\calB^{s-k-\frac{1}{p},\gamma-\frac{1}{2}}_p(\bz,F_k)
   \longrightarrow\hsgp(\dz,E)$$
for any $\gamma\in\rz$, $1<p<\infty$, and $s>\mu-1+\frac{1}{p}$, such that $\varrho K=1$ and   
both $\varphi K\omega$ and $\omega K\varphi$ map into $\calC^{\infty,\infty}(\dz,E)$ whenever 
the cut-off function $\omega$ and $\varphi\in\cicomp(\intd)$ have disjoint support. 
\end{lemma}
\begin{proof}
We shall construct $K$ in the form 
$K=\omega K^\prime\omega_0+(1-\omega)K^{\prime\prime}(1-\omega_1)$ with cut-off functions 
$\omega_1\prec\omega\prec\omega_0$. For $K^{\prime\prime}$ take any right-inverse to $\varrho$ 
on the smooth manifold with boundary $2\dz$, as for example constructed in Lemma 1.6.4 of 
\cite{Grub1}. Then $(1-\omega)K^{\prime\prime}(1-\omega_1)$ has the required mapping property. 
It remains to construct $K^\prime$ on $\rz_+\times X$. The construction is of local nature, so 
we may assume all bundles to be trivial and, for notational convenience, one dimensional, and 
we may work on $\rz_+\times\rz^n_+$. We then define 
$K^\prime=\begin{pmatrix}K_0&K_1&\ldots&K_{\mu-1}\end{pmatrix}$ by 
 $$(K_ju)(t,x^\prime,x_n)=\int\int_{\Gamma}e^{ix^\prime\xi^\prime}t^{-z}k_j(x_n,z,\xi^\prime)
   (\calM_{s\to z}\calF_{y^\prime\to\xi^\prime}u)(z,\xi^\prime)\,\dbar z\dbar\xi^\prime$$
with the symbol 
\begin{align*}
 k_j(x_n,z,\xi^\prime)&=\frac{x_n^j}{j!}\int\!\!\int_0^\infty
   s^{i\sigma-z}\varphi(s)\zeta(\spk{\xi^\prime,\sigma}x_n))\,\frac{ds}{s}\dbar\sigma\\
 &=\frac{x_n^j}{j!}\int\calM(s^{-\delta}\varphi)(i\sigma)
   \zeta(\spk{\xi^\prime,\tau+\sigma}x_n)\,\dbar\sigma,\qquad z=\delta+i\tau.
\end{align*}
Here, $\varphi\in\cicomp(\rz_+)$ with $\varphi\equiv1$ near $1$ and $\zeta\in\cicomp(\rz)$ 
with $\zeta\equiv1$ near 0 (the symbols $k_j$ arise by applying a so-called kernel cut-off 
procedure, cf.\ \cite{Schu2}, to the symbols used in in the proof of \cite{Grub1}, 
Lemma 1.6.4). Moreover, $\Gamma$ is any vertical line in the complex plane; due to the 
holomorphy of $k_j$ in $z$ and Cauchy's theorem, its choice does not effect the action of 
$K_j$ on $\calC^{\infty,\infty}(\rz_+)\pit\calS(\overline{\rz}^n_+)$. Since 
$\frac{d^l}{dx_n^l}k_j(x_n,z,\xi^\prime)|_{x_n=0}=\delta_{jl}$ (Kronecker symbol) for all 
$0\le l\le j$, $K^\prime$ defines a right-inverse to $\varrho$. Also $K^\prime$ has the 
desired continuous extensions (where one takes $\Gamma$ as the line 
$\re z=\frac{n+1}{2}-\gamma$), since (up to an order reduction of order $j$) $K_j$ is the 
typical local model for a Poisson operator in the flat cone algebra. This then also induces 
the stated smoothing behaviour of $\varphi_0K\varphi_1$. 
\end{proof}

\begin{remark}\label{real44.5}
The proof of Lemma \text{\rm\ref{real44}} yields the following: Given $0<\eps<1$, we can 
construct $K$ in such a way that $K(\omega u_0,\ldots,\omega u_{\mu-1})$ is supported in 
${[0,\eps[}\times X$, whenever $\omega$ is supported sufficiently close to zero.   
\end{remark}

The previous lemma shows that $\varrho$ is surjective. This is also true for general boundary 
conditions, provided they are normal in a sense we now specify. To this end we identify zero 
order Fuchs type differential operators on $\bz$ with vector bundle morphisms.  

\begin{definition}\label{real45}
Let $T=S\varrho$ as in \eqref{real4C}. We call $T$ a {\em normal} boundary condition if 
$\ulS_{kk}:E_\bz\to F_k$ is surjective for all $0\le k<\mu$. 
\end{definition}

\begin{proposition}\label{real46}
Let $T=S\varrho$ be a normal boundary condition. For $0\le k<\mu$ let 
$\ulS^\prime_{kk}:E_\bz\to E_\bz$ be projections onto $Z_k:=\text{\rm ker}\,\ulS_{kk}$ 
$($these are subbundles of $E_\bz$, since the $\ulS_{kk}$ are surjective$)$. Then 
\begin{equation}\label{real4E}
 \begin{pmatrix}\ulS\\ \ulS^\prime\end{pmatrix}:
  \mathop{\oplus}\limits_{j=0}^{\mu-1}
  \calB^{s-j-\frac{1}{p},\gamma-j-\frac{1}{2}}_p(\bz,E_\bz)
  \longrightarrow
  \begin{matrix}
   \mathop{\oplus}\limits_{k=0}^{\mu-1}
   \calB^{s-k-\frac{1}{p},\gamma-k-\frac{1}{2}}_p(\bz,F_k)\\ \oplus\\
   \mathop{\oplus}\limits_{k=0}^{\mu-1}
   \calB^{s-k-\frac{1}{p},\gamma-k-\frac{1}{2}}_p(\bz,Z_k)
  \end{matrix}
\end{equation}
with $\ulS^\prime=\diag(\ulS^\prime_{00},\ldots,\ulS^\prime_{\mu-1,\mu-1})$ is an isomorphism 
for any $\gamma,s\in\rz$, and $1<p<\infty$. If  
 $$\begin{pmatrix}\ulC & \ulC^\prime\end{pmatrix}:=
   \begin{pmatrix}\ulS\\ \ulS^\prime\end{pmatrix}^{-1},$$
then both $\ulC$ and $\ulC^\prime$ have the same structure as $\ulS$, i.e.\ they are left 
lower triangular matrices whose entry in the $j$-th row and $k$-th column are Fuchs type 
differential operators of order $j-k$.  
\end{proposition}
\begin{proof}
The proof is parallel to Lemma 1.6.1 of \cite{Grub1}. The map
$\wt\ulS=\begin{pmatrix}\diag(\ulS_{00},\ldots,\ulS_{\mu-1,\mu-1})\\ 
\ulS^\prime\end{pmatrix}$, acting as in \eqref{real4E}, is invertible by construction of 
$\ulS^\prime$. Then $\ulU=\wt\ulS^{-1}\begin{pmatrix}\ulS\\ \ulS^\prime\end{pmatrix}$ is a 
left lower triangular matrix of Fuchs type operators acting on sections of $E_\bz$, with 
identities in the diagonal. Hence $\ulU$ is invertible, since $\ulU=1-\ulV$ with a nilpotent 
matrix $\ulV$. By Neumann series, 
 $$\begin{pmatrix}\ulC & \ulC^\prime\end{pmatrix}=\ulU^{-1}\wt\ulS^{-1}=
   (1+\ulV+\ldots+\ulV^{\mu-1})\wt\ulS^{-1}$$
is of the described structure. 
\end{proof}

As a corollary, any normal boundary condition $T=S\varrho$ is surjective. In fact, if 
\begin{equation}\label{real4F}
 C:=\diag(1,t,\ldots,t^{\mu-1})\,\ulC
\end{equation}
with $\ulC$ from the previous proposition, and $K$ is as in Lemma \text{\rm\ref{real44}}, then 
 $$TKC=S\varrho KC=
   \ulS\,\diag(1,t^{-1},\ldots,t^{-\mu+1})\,\diag(1,t,\ldots,t^{\mu-1})\,\ulC=1.$$
Thus $KC$ is a right-inverse to $T$. This together with the mapping properties of $K$ and $C$ 
immediately yields the following lemma. 

\begin{lemma}\label{real47}
Let $T=S\varrho$ be normal, $K$ as in Lemma \text{\rm\ref{real45}}, and $C$ from 
\eqref{real4F}. Then 
\begin{equation}\label{real4G}
 P:=1-KCT:\hsgp(\dz,E)\longrightarrow\hsgp(\dz,E)
\end{equation}
continuously for any $\gamma\in\rz$, $1<p<\infty$, and $s>\mu-1+\frac{1}{p}$. Moreover, $P$ is 
a projection, i.e.\ $P^2=P$, and $\text{\rm im}\,P=\hsgp(\dz,E)_T$.
\end{lemma}

For later purpose it is useful to know further properties of $P$, stated in the next lemma. 

\begin{lemma}\label{real48}
Let $T=S\varrho$ be a normal boundary condition and $P$ as in \eqref{real4G}. Let 
$\omega,\omega_0,\omega_1$ be cut-off functions with $\omega_1\prec\omega\prec\omega_0$. 
Then  
\begin{itemize}
 \item[i)] $\omega u=0$ implies $\omega_1Pu\in\calC^{\infty,\infty}(\dz,E)$, and 
  $(1-\omega)u=0$ implies $(1-\omega_0)Pu\in\calC^{\infty,\infty}(\dz,E)$, 
 \item[ii)] $\omega Tu=0$ implies $\omega_1 Pu-\omega_1u\in\calC^{\infty,\infty}(\dz,E)$.
 \item[iii)] $\omega Tu\in\mathop{\oplus}\limits_{k=0}^{\mu-1}
  \calB^{s-k-\frac{1}{p},\gamma+r-k-\frac{1}{2}}_p(\bz,F_k)$ for some $r\ge0$ implies 
  $\omega_1 Pu-\omega_1u\in\calH^{s,\gamma+r}_p(\dz,E)$. 
\end{itemize}
Here, $u\in\hsgp(\dz,E)$ with $\gamma\in\rz$, $1<p<\infty$, and $s>\mu-1+\frac{1}{p}$. 
\end{lemma}
\begin{proof}
Let us show the first part of i). We write 
 $$\omega_1Pu=\omega_1u-\omega_1KCT\omega u+\omega_1KCT(1-\omega) u.$$
The second summand vanishes by assumption, the third belongs to $\calC^{\infty,\infty}(\dz,E)$ 
due to the locality of $C$ and $T$, and due to Lemma \ref{real44}. The other claims are proved 
analogously. 
\end{proof}

\begin{corollary}\label{real49}
Let $T$ be a normal boundary condition. Then $\calC^{\infty,\infty}(\dz,E)_T$ is a dense subset 
of the closed subspace $\hsgp(\dz,E)_T$ of $\hsgp(\dz,E)$ for any $\gamma\in\rz$, $1<p<\infty$, 
and  $s>\mu-1+\frac{1}{p}$. 
\end{corollary}
\begin{proof}
Let $u\in\hsgp(\dz,E)$ with $Tu=0$ be given. There exists a sequence $(v_n)_{n\in\nz}$ in 
$\cicomp(\intd,E)$ converging to $u$ in $\hsgp(\dz,E)$. 
Now $u_n:=Pv_n$ vanishes under $T$ and belongs to $\calC^{\infty,\infty}(\dz,E)$ by Lemma 
\ref{real48}.i). By Lemma \ref{real47}, $u_n\to Pu=u$ in $\hsgp(\dz,E)$. 
\end{proof}

\subsection{Green's formula and the adjoint problem}

We denote by $A^t$ the formal adjoint of $A$. In coordinates $(t,x^\prime,x_n)$ near 
${]0,1]}\times\partial X$ let us write  
\begin{equation}\label{real5D} 
 A=t^{- \mu}\smsum_{l=0}^\mu
 \smsum_{j=0}^{\mu-l}s_{jl}(t,x^\prime,x_n,D^\prime)(-t\partial_t)^j D_n^l
\end{equation}
with differential operators $s_{jl}(t,x^\prime,x_n,D^\prime)$ of order $\mu-j-l$ on 
$\partial X$ (depending on the parameter $(t,x_n)$). 

\begin{theorem}\label{real54}
For all $u\in\calH^{\mu,\gamma+\mu}_p(\dz,E)$ and 
$v\in\calH^{\mu,-\gamma}_{p^\prime}(\dz,E)$ we have 
\begin{equation}\label{real5E} 
 \skp{Au}{v}_{\calH^{0,0}_2(\dz,E)}-\skp{u}{A^tv}_{\calH^{0,0}_2(\dz,E)}=
 \skp{\mathfrak{A}\varrho u}{\varrho v}_{\calB^{0,0}_2(\bz,E_\bz^\mu)},
\end{equation}
where $E_\bz^\mu$ is the $\mu$-fold direct sum of $E_\bz$ and 
$\mathfrak{A}=(\mathfrak{A}_{jk})_{0\le j,k<\mu}$ is a left upper triangular matrix with 
 $$\mathfrak{A}_{jk}=t^{-j}\,\underline{\mathfrak{A}}_{jk}\,t^{-k}$$
and cone differential operators 
$\underline{\mathfrak{A}}_{jk}$ of order $\mu-1-j-k$ acting on sections of $E_\bz$ and 
$\underline{\mathfrak{A}}_{jk}=0$ if $j+k\ge\mu$. Setting 
$\underline{\mathfrak{A}}=(\underline{\mathfrak{A}}_{jk})_{0\le j,k<\mu}$ we can also write 
 $$\mathfrak{A}=\text{\rm diag}(1,t^{-1},\ldots,t^{-\mu+1})\,
   \underline{\mathfrak{A}}\,\text{\rm diag}(1,t^{-1},\ldots,t^{-\mu+1}).$$
Explicitly, near the singularity, 
 $$\underline{\mathfrak{A}}_{jk}=t^{-\mu+j+k+1}\smsum_{m=0}^{\mu-j-k-1}\smsum_{l=j+k+1}^{\mu-m}
   i((-D_n)^{l-j-k-1}s_{ml})(t,x^\prime,0,D^\prime)(-t\partial_t-k)^m,$$
where we use the notation from \eqref{real5D}. In particular, for $0\le j<\mu$, 
 $$\underline{\mathfrak{A}}_{j,\mu-j-1}=\underline{\mathfrak{A}}_{0,\mu-1}=
   s_{0\mu}(t,x^\prime,0,D^\prime).$$
\end{theorem}

The proof is a standard but lengthy induction based on the relation 
 $$\skp{D_nu}{v}_{\calH^{0,0}_2(\rz_+\times\overline{\rz}^n_+)}-
   \skp{u}{D_nv}_{\calH^{0,0}_2(\rz_+\times\overline{\rz}^n_+)}=
   i\skp{\gamma_0u}{\gamma_0v}_{\calB^{0,-1/2}_2(\rz_+\times\rz^{n-1})},$$
for functions $u,v\in\cicomp(\rz_+\times\overline{\rz}^n_+)$ and the fact that 
$\skp{f}{g}_{\calB^{0,-1/2}_2(\rz_+\times\rz^{n-1})}=
\skp{tf}{g}_{\calB^{0,0}_2(\rz_+\times\rz^{n-1})}$. 
To this end note that 
$\calH^{0,0}_2(\rz_+\times\overline{\rz}^n_+)=L_2(\rz_+\times\overline{\rz}^n_+,t^{n+1}dtdx)$ 
and $\calB^{0,0}_2(\rz_+\times\rz^{n-1})=L_2(\rz_+\times\rz^{n-1},t^ndtdx^\prime)$. We shall 
omit the details. 

\begin{definition}\label{real55}
The adjoint problem of 
$\begin{pmatrix}A\\T\end{pmatrix}$ is $\begin{pmatrix}A^t\\\wt T\end{pmatrix}$ with 
 $$\wt{T}=\wt{S}\varrho=
   \wt\ulS\,\text{\rm diag}(1,t^{-1},\ldots,t^{-\mu+1})\,
   \varrho,\qquad\wt\ulS:=I^\times(\underline{\mathfrak{A}}\,\ulC^\prime)^t.$$
Here, $\underline{\mathfrak{A}}$ is from \eqref{real5E}, $\ulC^\prime$ from Proposition 
\mbox{\rm\ref{real46}}, and $I^\times=(\delta_{j+k,\mu-1})_{0\le j,k<\mu}$ is the skew unit 
matrix. The superscript $t$ refers to taking formal adjoints in each entry of the matrix and 
then passing to the transposed matrix. 
\end{definition} 

We note that 
\begin{equation}\label{real315a}
 \wt\ulS_{jk}=\smsum_{l=0}^{\mu-1}
 \ulC^{\prime t}_{l,\mu-j-1}\underline{\mathfrak{A}}_{kl}^{t}
\end{equation} 
is a cone differential operator of order $j-k$ and equal to zero if $k>j$. Letting  
$Z_k=\text{\rm ker}\,\ulS_{kk}$ for $0\le k<\mu$ the fact that $\ulC^{\prime t}_{jk}$ maps 
to sections of $Z_k$ implies that 
\begin{equation}\label{real315b}
 \wt T:\calH^{s,\gamma}_p(\dz,E)\longrightarrow
 \mathop{\oplus}\limits_{k=0}^{\mu-1}
 \calB^{s-k-\frac{1}{p},\gamma-k-\frac{1}{2}}_p(\bz,Z_{\mu-1-k}).
\end{equation} 
In general, $\wt T$ is not uniquely determined, since in Proposition \ref{real46} we may have 
some freedom in constructing $\ulC^\prime$. However, any of the possible choices for $\wt T$ 
is equally good. 

\begin{proposition}\label{real56}
If \ $\underline{\mathfrak{A}}_{0,\mu-1}$ is a  homeomorphism then $\wt T$ is normal. 
\end{proposition}
\begin{proof}
By \eqref{real315a}, 
$\wt\ulS_{kk}=
\ulC^{\prime t}_{\mu-k-1,\mu-k-1}\underline{\mathfrak{A}}_{k,\mu-k-1}^{t}=
\ulC^{\prime t}_{\mu-k-1,\mu-k-1}\underline{\mathfrak{A}}_{0,\mu-1}^{t}$, 
since the only nonzero summand is that for $l=\mu-k-1$. However, 
by construction the $\ulC^\prime_{kk}$ are injective morphisms, hence the 
$\ulC^{\prime t}_{kk}$ are surjective. The assumption yields the surjectivity of all 
$\wt\ulS_{kk}$, i.e.\ normality of $\wt T$.
\end{proof}

We shall say that $\bz$ {\em non-characteristic for} $A$ if the assumption on 
$\underline{\mathfrak{A}}_{0,\mu-1}$ in Proposition \ref{real56} is satisfied.


\section{The minimal extension and its adjoint}\label{section5}

Let us now consider $A$ as the unbounded operator 
\begin{equation}\label{real5B}
 A:\cii(\dz,E)_T\subset\calH^{0,\gamma}_p(\dz,E)\longrightarrow\calH^{0,\gamma}_p(\dz,E)
\end{equation}
with a {\em normal} boundary condition $T=S\varrho$. This operator is densely defined and 
closable, since $A$ acts continuously in the scale of Sobolev spaces. We shall assume that 
$\begin{pmatrix}A\\T\end{pmatrix}$ satisfies the ellipticity conditions i) and ii) of 
Proposition \ref{real42}. This we shall refer to as $\dz$-{\em ellipticity}. Note that then 
the conormal symbol is meromorphically invertible, cf.\ Remark \ref{real316}, and that $\bz$ 
is non-characteristic for $A$. 

\begin{definition}\label{real51}
We denote the closure (respectively minimal extension) of $A$ from \eqref{real5B} by 
$A_{T,\min}$.
\end{definition}

\begin{proposition}\label{real52}
We have the inclusions 
\begin{equation}\label{real5C}
 \calH^{\mu,\gamma+\mu}_p(\dz,E)_T\subset\calD(A_{T,\min})\subset
 \Big\{u\in\schnitt_{\eps>0}\calH^{\mu,\gamma+\mu-\eps}_p(\dz,E)_T\st 
 Au\in\calH^{0,\gamma}_p(\dz,E)\Big\}.
\end{equation}
If, additionally, $\begin{pmatrix}A\\T\end{pmatrix}$ is elliptic with respect to $\gamma+\mu$ 
then 
 $$\calD(A_{T,\min})=\calH^{\mu,\gamma+\mu}_p(\dz,E)_T.$$
\end{proposition}
\begin{proof}
Let $u\in\calH^{\mu,\gamma+\mu}_p(\dz,E)_T$. By Corollary \ref{real49} there exist 
$u_n\in\cii(\dz,E)_T$ such that $u_n\to u$ in $\calH^{\mu,\gamma+\mu}_p(\dz,E)$. Then 
$Au_n\to Au$ in $\calH^{0,\gamma}_p(\dz,E)$ and the first inclusion of \eqref{real5C} is 
verified. Next assume ellipticity with respect to $\gamma+\mu$ and let $u\in\calD(A_{T,\min})$. 
Then there exists a sequence $(u_n)_{n\in\nz}\subset\cii(\dz,E)_T$ such that $u_n\to u$ and 
$Au_n\to Au$ in $\calH^{0,\gamma}_p(\dz,E)$. Let $\begin{pmatrix}R&K\end{pmatrix}$ be the 
parametrix as constructed in Proposition \ref{real43} (with $\gamma$ replaced by $\gamma+\mu$). 
Then, in particular, 
 $$G:=RA-1:\calH^{0,\gamma+\mu}_p(\dz,E)_T\longrightarrow\calH^{\mu,\gamma+\mu}_p(\dz,E)$$
is a finite rank operator. The sequence $Gu_n=R(Au_n)-u_n$ converges in 
$\calH^{0,\gamma}_p(\dz,E)$, hence also in $\calH^{\mu,\gamma+\mu}_p(\dz,E)$ due to 
$\text{\rm dim im}\,G<\infty$. Since $TR=0$, it even converges in 
$\calH^{\mu,\gamma+\mu}_p(\dz,E)_T$. Consequently, $u_n=R(Au_n)-Gu_n$ converges in 
$\calH^{\mu,\gamma+\mu}_p(\dz,E)_T$. Thus we obtain the above stated identity. For the 
remaining inclusion in \eqref{real5C} note that in any case $\begin{pmatrix}A\\T\end{pmatrix}$ 
is elliptic with respect to $\gamma+\mu-\eps$ for all sufficiently small $\eps>0$, 
cf.\ Remark \ref{real316}. Then we argue as before. 
\end{proof}

We shall improve the second inclusion in \eqref{real5C} to an equality, cf.\ 
Theorem \ref{real514}. 

\subsection{The adjoint of the closure}

Our next aim is to study the adjoint of $A_{T,\min}$. A first simple result is the following 
lemma. 

\begin{lemma}\label{real53}
$A_{T,\min}^*$ is given by the action of $A^t$ on $\calD(A_{T,\min}^*)$. 
\end{lemma}
\begin{proof}
Since $A_{T,\min}$ is the closure of the operator in \eqref{real5B} their adjoints coincide and 
\begin{equation}\label{42a}
\begin{split}
 \calD(A_{T,\min}^*)=\Big\{u\in\calH^{0,-\gamma}_{p^\prime}(\dz,E)\st&\exists\,
   v\in\calH^{0,-\gamma}_{p^\prime}(\dz,E)\quad
   \forall\;w\in\cii(\dz,E)_T:\\
   &\skp{w}{v}_{\calH^{0,0}_2(\dz,E)}=\skp{Aw}{u}_{\calH^{0,0}_2(\dz,E)}\Big\}
\end{split}
\end{equation}
and then $A_{T,\min}^*u=v$. To given $u\in\calD(A_{T,\min}^*)$ there exists a sequence 
$u_n\in\cicomp(\kringel{\dz},E)$ with $u_n\to u$ in $\calH^{0,-\gamma}_{p^\prime}(\dz,E)$. 
Therefore, with $\skp{\cdot}{\cdot}$ denoting the scalar-product of $\calH^{0,0}_2(\dz,E)$,  
 $$\skp{w}{v}=\skp{Aw}{u}\longleftarrow\skp{Aw}{u_n}=\skp{w}{A^tu_n}
   \longrightarrow\skp{w}{A^tu}$$
due to the duality of $\calH^{s,\gamma}_p(\dz,E)$ and 
$\kringel{\calH}^{-s,-\gamma}_{p^\prime}(\dz,E)$. In particular, this is true for all 
$w\in\cicomp(\kringel{\dz},E)$ and thus, by density, $v=A^tu$.  
\end{proof}

\begin{definition}\label{real57}
The operator $A_{T,\max}$ is defined by the action of $A$ on the domain   
 $$\calD(A_{T,\max})=\{u\in\calH^{\mu,\gamma}_p(\dz,E)_T\st Au\in\calH^{0,\gamma}_p(\dz,E)\}.$$
\end{definition}

Observe that $A_{T,\max}$ is not the maximal extension of $A$ from \eqref{real5B} in the usual 
sense (which would have the domain 
$\{u\in\calH^{0,\gamma}(\dz,E)\st Au\in\calH^{0,\gamma}(\dz,E)\}$). Instead, it is the largest 
extension of $A$ in $\calH^{\mu,\gamma}_p(\dz,E)_T$. We shall call it the maximal extension of 
$A$. This is also justified by the following result. 

\begin{theorem}\label{real58}
As unbounded operators in $\calH^{0,-\gamma}_{p^\prime}(\dz,E)$ and 
$\calH^{0,\gamma}_p(\dz,E)$, respectively, we have 
\begin{equation}\label{42b}
 A_{T,\min}^*=(A^t)_{\wt{T},\max},\qquad A_{T,\min}=(A^t)_{\wt{T},\max}^*,
\end{equation}
where $\wt{T}$ is the adjoint boundary condition of Definition \text{\rm\ref{real55}}.
\end{theorem} 
\begin{proof}
It is enough to prove the first identity, for if it is valid then 
 $$(A^t)_{\wt{T},\max}^*=A_{T,\min}^{**}=\overline{A_{T,\min}}=A_{T,\min}.$$
To do so, we first shall show that 
\begin{equation}\label{real5F}
 \{u\in\calH^{\mu,-\gamma}_{p^\prime}(\dz,E)_{\wt{T}}\st 
 A^tu\in\calH^{0,-\gamma}_{p^\prime}(\dz,E)\}=
 \calD(A_{T,\min}^*)\cap\calH^{\mu,-\gamma}_{p^\prime}(\dz,E).
\end{equation}
Let $u\in\calH^{\mu,-\gamma}_{p^\prime}(\dz,E)$, $w\in\cii(\dz,E)_T$, and 
$\begin{pmatrix}\ulS\\ \ulS^\prime\end{pmatrix}$ as in Proposition \ref{real46} with 
corresponding inverse $\begin{pmatrix}\ulC&\ulC^\prime\end{pmatrix}$. If we set 
$D=\diag(1,t^{-1},\ldots,t^{-\mu+1})$, and write 
$\mathfrak{A}=D\underline{\mathfrak{A}}(\ulC\,\ulS+\ulC^\prime\ulS^\prime)D$ we obtain 
 $$\skp{\mathfrak{A}\varrho w}{\varrho u}=
   \skp{Tw}{(\underline{\mathfrak{A}}\,\ulC)^tD\varrho u}+
   \skp{\ulS^\prime D\varrho w}{I^\times\wt T u}=
   \skp{\ulS^\prime D\varrho w}{I^\times\wt T u},$$
where $\skp{\cdot}{\cdot}$ denotes the scalar product of $\calB^{0,0}_2(\bz,E^\mu_\bz)$. 
Hence, if $\wt Tu=0$, Greens formula \eqref{real5D} yields 
$\skp{Aw}{u}_{0,0}=\skp{w}{A^t u}_{0,0}$, with $\skp{\cdot}{\cdot}_{0,0}$ denoting the scalar 
product of $\calH^{0,0}_2(\dz,E)$. By \eqref{42a}, this shows that the right-hand side of 
\eqref{real5F} contains the left-hand side. Vice versa, if $u$ is an element from the 
right-hand side, 
 $$\skp{Aw}{u}_{0,0}=\skp{w}{A_{T,\min}^*u}_{0,0}=\skp{w}{A^tu}_{0,0}$$
since $u\in\calD(A_{T,\min}^*)$. On the other hand, Greens formula for 
$u\in\calH^{\mu,-\gamma}_{p^\prime}(\dz,E)$ shows that  
 $$\skp{Aw}{u}_{0,0}=\skp{w}{A^tu}_{0,0}+
   \skp{\ulS^\prime D\varrho w}{I^\times\wt Tu}.$$
Hence $\skp{\ulS^\prime D\varrho w}{I^\times\wt Tu}=0$ for all $w\in\cii(\dz,E)_T$. But now, 
by construction, 
$\ulS^\prime D\varrho:\cii(\dz,E)_T\to
\mathop{\mbox{$\oplus$}}\limits_{k=0}^{\mu-1}\cii(\bz,Z_k)$ 
with $Z_k=\textrm{ker}\,\ulS_{kk}$ is surjective. We infer from \eqref{real315b} that 
$\wt Tu$ must vanish. 

We have verified \eqref{real5F}. Next we shall show that  
$\calD(A_{T,\min}^*)\subset\calH^{\mu,-\gamma}_{p^\prime}(\dz,E)$. 

First let us assume that, additionally, $\begin{pmatrix}A\\T\end{pmatrix}$ is elliptic with 
respect to $\gamma+\mu$. Then, by Proposition \ref{real43}, there exists an 
$R\in C^{-\mu,0}(\dz;\gamma,\gamma+\mu,k;E;E)$ such that 
 $$R:\calH^{0,\gamma}_p(\dz,E)\longrightarrow\calH^{\mu,\gamma+\mu}_p(\dz,E)_T=
   \calD(A_{T,\min})$$
and that $G:=AR-1\in C^0_G(\dz;\gamma,\gamma,k;E;E)$ is a Green operator. By general facts on 
the adjoint of compositions, $R^*A_{T,\min}^*\subset(AR)^*=1+G^*$. Then 
$u=R^*A_{T,\min}^*u-G^*u$ for $u\in\calD(A_{T,\min}^*)$ yields 
$u\in\calH^{\mu,-\gamma}_{p^\prime}(\dz,E)$ by Theorem \ref{real315.5}.

In general, we have ellipticity with respect to $\gamma+\mu+\eps$ for any sufficiently small 
$\eps>0$. If we now denote by $\wt{A}$ the unbounded operator in 
$\calH^{0,\gamma+\eps}_p(\dz,E)$ acting as $A$ on the domain $\cii(\dz,E)_T$, the previous 
step shows that 
 $$\{u\in\calH^{\mu,-\gamma-\eps}_{p^\prime}(\dz,E)_{\wt{T}}\st
   A^tu\in\calH^{0,-\gamma-\eps}_{p^\prime}(\dz,E)\}=
   \calD(\wt{A}_{T,\min}^*).$$
By the definition of the domain of the adjoint, it is clear that 
$\calD(A_{T,\min}^*)\subset\calD(\wt{A}_{T,\min}^*)$. Passing to the intersection over all 
$\eps$, we obtain 
 $$\calD(A_{T,\min}^*)\subset\schnitt_{\eps>0}
   \{u\in\calH^{\mu,-\gamma-\eps}_{p^\prime}(\dz,E)_{\wt{T}}\st
   A^tu\in\calH^{0,-\gamma-\eps}_{p^\prime}(\dz,E)\}
   \,\cap\,\calH^{0,-\gamma}_{p^\prime}(\dz,E).$$ 
Now let $u$ be an element of the right-hand side. As shown in Corollary \ref{real73}, the 
adjoint problem $\begin{pmatrix}A^t\\\wt{T}\end{pmatrix}$ is elliptic with respect to 
$-\gamma-\eps$ for small $\eps$. From the inclusion 
$\calH^{0,-\gamma-\eps}_{p^\prime}(\dz,E)\subset
\calH^{0,-\gamma-\mu-\eps}_{p^\prime,\mu}(\dz,E)$, cf.\ Definition \ref{real35}, and elliptic 
regularity, cf.\ Theorem \ref{real320}.b), we obtain
\begin{equation}\label{43a}
 u\in\calH^{\mu,-\gamma-\eps}_{p^\prime,P}(\dz,E)=
 \schnitt_{\delta>0}\calH^{\mu,-\gamma-\eps+\mu-\delta}_{p^\prime}(\dz,E)\oplus
 \calE_P(\rz_+\times X)\subset
 \calH^{\mu,-\gamma}_{p^\prime}(\dz,E)+\calE_P(\rz_+\times X)
\end{equation}
for some asymptotic type $P=\{(p,m,M)\}\in\as(X;-\gamma-\eps,\mu)$. Since 
$u\in\calH^{0,-\gamma}_{p^\prime}(\dz,E)$ we must have $\re p<\frac{n+1}{2}+\gamma$ for all 
$(p,m,M)\in P$, and therefore $u\in\calH^{\mu,-\gamma}_{p^\prime}(\dz,E)$. 
\end{proof}

From \eqref{real5C} and \eqref{42b} we conclude: 

\begin{corollary}\label{real59}
$\skp{Au}{v}_{\calH^{0,0}_2(\dz,E)}=\skp{u}{A^tv}_{\calH^{0,0}_2(\dz,E)}$ for all 
$u\in\calH^{\mu,\gamma+\mu}_p(\dz,E)_T$, 
$v\in\calH^{\mu,-\gamma}_{p^\prime}(\dz,E)_{\wt T}$.
\end{corollary}

Identity \eqref{43a} applied to $A$ shows the following:  

\begin{corollary}\label{real510}
There exists an $\eps>0$ such that
$\calD(A_{T,\max})\subset\calH^{\mu,\gamma+\eps}_p(\dz,E)_T$.
\end{corollary}

\begin{theorem}\label{real510.5}
As unbounded operators in $\calH^{0,-\gamma}_{p^\prime}(\dz,E)$ and $\calH^{0,\gamma}_p(\dz,E)$, 
respectively, we have 
 $$A_{T,\max}^*=(A^t)_{\wt{T},\min},\qquad
   A_{T,\max}=(A^t)_{\wt{T},\min}^*,$$
where $\wt{T}$ is the adjoint boundary condition of Definition \text{\rm\ref{real55}}.
\end{theorem}
\begin{proof}
It is enough to prove the second identity, since the first follows from passing to the adjoints. 
To do so, note that $\calD((A^t)_{\wt{T},\min}^*)$ is given by 
\begin{equation*}
 \Big\{u\in\calH^{0,\gamma}_p(\dz,E)\st\exists\,v\in\calH^{0,\gamma}_p(\dz,E)\quad
 \forall\;w\in\cii(\dz,E)_{\wt{T}}:\quad (w,v)_{0,0}=(A^tw,u)_{0,0}\Big\}, 
\end{equation*}
and that Green's formula for the formal adjoint $A^t$ reads as 
 $$(A^t w,u)_{0,0}-(w,Au)_{0,0}=
   (-\mathfrak{A}^t\varrho w,\varrho u),$$
with $\skp{\cdot}{\cdot}_{0,0}$ and $\skp{\cdot}{\cdot}$ denoting the scalar product of 
$\calH^{0,0}_2(\dz,E)$ and $\calB^{0,0}_2(\dz,E^\mu_\bz)$, respectively. 

Now let $u\in\calH^{\mu,\gamma}_p(\dz,E)$ be given and $w\in\cii(\dz,E)_{\wt{T}}$ be arbitrary. 
 
Using $\underline{\mathfrak{A}}=\underline{\mathfrak{A}}(\ulC\ulS+\ulC^\prime\ulS^\prime)$ 
and setting 
$D=\textrm{diag}(1,t^{-1},\ldots,t^{-\mu+1})$, we can write 
 $$\mathfrak{A}^t\varrho=D\ulS^t\ulC^t\underline{\mathfrak{A}}^tD\varrho+
                         D\ulS^{\prime t}\ulC^{\prime t}\underline{\mathfrak{A}}^tD\varrho  
                        =D\ulS^t\ulC^t\underline{\mathfrak{A}}^tD\varrho+
                         D\ulS^{\prime t}I^\times \wt{T}.$$
This yields
\begin{equation}\label{flick1} 
 (-\mathfrak{A}^t\varrho w,\varrho u)
 =-(D\ulS^t\ulC^t\underline{\mathfrak{A}}^tD\varrho w,\varrho u)
 =-(\ulC^t\underline{\mathfrak{A}}^tD\varrho w,Tu).
\end{equation}
Thus if $Tu=0$, i.e.\ $u\in\calD(A_{T,\max})$, the right-hand side equals zero and together 
with the above Green's formula we obtain $u\in\calD((A^t)_{\wt{T},\min}^*)$ and 
$A_{T,\max}\subset(A^t)_{\wt{T},\min}^*$. 

By Theorem \ref{real58} we know that 
$(A^t)_{\wt{T},\min}^*=A_{\overset{\approx}{T},\max}$, where 
$\overset{\approx}{T}=\overset{\approx}{\ulS}D\varrho$ with 
$\overset{\approx}{\ulS}=-I^\times\wt{\ulC}{}^{\prime t}\underline{\mathfrak{A}}$ is the adjoint 
boundary condition to $\wt{T}$, obtained via Proposition \ref{real46}. Thus, for 
$u\in\calD((A^t)_{\wt{T},\min}^*)\subset\calH^{\mu,\gamma}_p(\dz,E)$, it follows from 
Green's formula and \eqref{flick1} that 
 $$(D\varrho w,\underline{\mathfrak{A}}\ulC Tu)=
   (\ulC^t\underline{\mathfrak{A}}^tD\varrho w,Tu)=0$$ 
for all $w\in\cii(\dz,E)_{\wt{T}}$. Using 
$\wt{\ulC}\wt{\ulS}+\wt{\ulC}{}^\prime\wt{\ulS}{}^\prime=1$, we write 
$D\varrho w=\wt{\ulC}\wt{T}w+\wt{\ulC}{}^\prime\wt{\ulS}{}^\prime D\varrho w$. Since 
$\wt{T}w=0$ and by surjectivity of $\wt{\ulS}{}^\prime D\varrho$, we conclude that 
$\wt{\ulC}{}^{\prime t}\underline{\mathfrak{A}}\ulC Tu=0$. Multiplying this equation with 
$\wt{\ulS}{}^{\prime t}$ and using $\wt{\ulC}{}^\prime\wt{\ulS}{}^\prime=1-\wt{\ulC}\wt{\ulS}$ 
and $\wt{\ulS}{}^{t}=\underline{\mathfrak{A}}\ulC^\prime I^\times$ yields 
 $$0=\underline{\mathfrak{A}}\begin{pmatrix}\ulC&\ulC^\prime\end{pmatrix}
   \begin{pmatrix}Tu\\-I^\times\wt{\ulC}{}^t\underline{\mathfrak{A}}\ulC Tu\end{pmatrix}.$$
However, $\begin{pmatrix}\ulC&\ulC^\prime\end{pmatrix}$ is bijective by construction and 
$\underline{\mathfrak{A}}$ is, since $\partial\bz$ is non-characteristic. In particular, 
it follows that $Tu=0$, and therefore $(A^t)_{\wt{T},\min}^*\subset A_{T,\max}$. 
\end{proof}

\subsection{Fredholm property and description of the minimal domain}

Let $\begin{pmatrix}r_0&k_0\end{pmatrix}$ denote the inverted conormal symbol of 
$\begin{pmatrix}A\\T\end{pmatrix}$. If $\begin{pmatrix}R&K\end{pmatrix}$ is the parametrix 
from Proposition \ref{real43} (with $\gamma$ replaced by $\gamma+\mu-\eps$) and $P$ denotes 
the projection from \eqref{real4G}, we define 
\begin{equation}\label{real5G}
 B=P\,\omega\,\op_M^{\gamma+\mu-\eps-\frac{n}{2}}(r_0)\,t^\mu\,\omega_0 + 
 P\,(1-\omega)\,R\,(1-\omega_1).
\end{equation}
Here, $\omega,\omega_0,\omega_1$ are cut-off functions with $\omega_1\prec\omega\prec\omega_0$, 
and $0<\eps<1$ is chosen so small that $r_0$ has no pole in the strip 
$\frac{n+1}{2}-\gamma-\mu<\re z\le\frac{n+1}{2}-\gamma-\mu+\eps$. This choice guarantees that 
$B$ as an operator on $\calH^{0,\gamma}_p(\dz,E)$ does not depend on the choice of $\eps$. 

\begin{lemma}\label{real512}
For given $T=S\varrho$ define
 $$T_0=\omega\,\diag(1,t^{-1},\ldots,t^{-\mu+1})\,\op_M(\sigma_M^\mu(T))+(1-\omega)T.$$
Then the cut-off function $\omega$ can be chosen in such a way that $T_0$ is a normal boundary 
condition, $T-T_0=tT_1$ for a resulting boundary condition $T_1$, and 
$\begin{pmatrix}A\\T_0\end{pmatrix}$ is $\dz$-elliptic.
\end{lemma}
\begin{proof}
Since the $\ulS_{jk}$ are Fuchs type differential operators of order $j-k$, we can write 
 $$\omega\,\ulS_{jk}=\omega\,t^{k-j}\,\op_M(\sigma_M^{j-k}(\ulS_{jk}))+
   \omega\,t\,\ulS^1_{jk}$$
with Fuchs type operators $\ulS^1_{jk}$ of order $j-k$. Therefore 
 $$\omega\,\ulS_{jk}\,t^{-k}=\omega\,t^{-j}\,\op_M(T^k\sigma_M^{j-k}(\ulS_{jk}))
   +\omega\,t\,\ulS^1_{jk}\,t^{-k}.$$
By definition of the conormal symbol of $T$ it follows that $T-T_0=tS_1\varrho$, where 
$\ulS_1=(\omega\,\ulS_{jk}^1)_{j,k}$. If we write $T_0=S_0\varrho$ with 
$\ulS_0=(\ulS_{jk}^0)_{j,k}$ then, by construction, 
 $$\ulS_{kk}-\ulS_{kk}^0=\omega\,t\,\ulS_{kk}^1,\qquad 0\le k<\mu.$$
If $\omega$ is supported sufficiently close to $t=0$ we see that with $\ulS_{kk}$ also 
$\ulS_{kk}^0$ is surjective. 

By construction the rescaled symbols of $\begin{pmatrix}A\\T\end{pmatrix}$ and 
$\begin{pmatrix}A\\T_0\end{pmatrix}$ coincide, and are invertible in case of ellipticity. Then 
also the (not rescaled) symbols are invertible for small $t$, say $t\le\eps$ (recall that one 
defines the rescaled symbols by a limit procedure). Replacing $\omega$ by an 
$\wt\omega\prec\omega$ with support $[0,\eps]$, we obtain invertibility of the symbols of 
$\begin{pmatrix}A\\T_0\end{pmatrix}$ for $t<\eps$. For $t\ge\eps$ they coincide with the 
symbols of $\begin{pmatrix}A\\T\end{pmatrix}$, hence are invertible there. 
\end{proof}

\begin{lemma}\label{real513}
Let $T$ and $T_0$ be as in the previous Lemma \text{\rm\ref{real512}} and let $P$ be any 
projection associated to $T$ by \eqref{real4G}. Then we can choose a projection $P_0$ 
associated with $T_0$ such that 
\begin{equation}\label{44a}
 P-P_0:\hsgp(\dz,E)\longrightarrow\calH^{s,\gamma+1}_p(\dz,E)
\end{equation}
for any $\gamma\in\rz$, $1<p<\infty$ and $s>\mu-1+\frac{1}{p}$.  
\end{lemma}
\begin{proof}
To prove this statement we need to have a closer look at the construction of $P$ and $P_0$, 
which relies on Proposition \ref{real46}. We write 
$P-P_0=\wt\omega(P-P_0)+(1-\wt\omega)(P-P_0)$ for a cut-off function $\wt\omega$. 
The second summand trivially has the desired mapping property. We rewrite the first as 
 $$\wt\omega(P-P_0)=\wt\omega(KCT-KC_0T_0)=\wt\omega K(C(T-T_0)+(C-C_0)T_0)
   =\wt\omega\,Kt\wt CT_1+\wt\omega\,KtC_1T_0$$
with $\wt C=t^{-1}Ct$, $T_1$ as in Lemma \ref{real512}, and $C_1=t^{-1}(C_0-C)$. We conclude 
from Lemma \ref{real44} and \eqref{210} that the first summand has the mapping property in 
\eqref{44a} and focus on the second. Write $T_0=S_0\varrho$ and $T_1=S_1\varrho$ as in the 
proof of Lemma \ref{real512}. Let 
$\begin{pmatrix}\ulC&\ulC^\prime\end{pmatrix}=
\begin{pmatrix}\ulS\\ \ulS^\prime\end{pmatrix}^{-1}$ and 
$\begin{pmatrix}\ulC_0&\ulC_0^\prime\end{pmatrix}=
\begin{pmatrix}\ulS_0\\ \ulS_0^\prime\end{pmatrix}^{-1}$ be as in Proposition \ref{real46}. 
Now assume the existence of a bundle morphism 
$\ulS_1^\prime=\diag(\ulS_{00}^{1\prime},\ldots,\ulS_{\mu-1,\mu-1}^{1\prime})$ such that 
\begin{equation}\label{real5G.5}
\begin{pmatrix}\ulS\\ \ulS^\prime\end{pmatrix}-
\begin{pmatrix}\ulS_0\\ \ulS_0^\prime\end{pmatrix}=t
\begin{pmatrix}\ulS_1\\ \ulS_1^\prime\end{pmatrix}.
\end{equation}
Multiplying this equation from the left with 
$\begin{pmatrix}\ulC&\ulC^\prime\end{pmatrix}$ and with 
$\begin{pmatrix}\ulC_0&\ulC_0^\prime\end{pmatrix}$ from the right yields 
 $$\begin{pmatrix}\ulC_0-\ulC&\ulC_0^\prime-\ulC^\prime\end{pmatrix}=
   \begin{pmatrix}\ulC t \ulS_1\ulC_0+\ulC^\prime t \ulS_1^\prime\ulC_0&
      \ulC t \ulS_1\ulC_0^\prime+\ulC^\prime t \ulS_1^\prime\ulC_0^\prime
   \end{pmatrix}$$
so that $\ulC_1=t^{-1}(\ulC_0-\ulC)$ is a left lower triangular matrix of Fuchs type 
differential operators of order $j-k$ at position $(j,k)$. Multiplying from the left with 
$\textrm{diag}(1,\ldots,t^{-\mu-1})$ we obtain $C_1$ and conclude that $\wt\omega\,KtC_1T_0$ 
has the mapping property in \eqref{44a}. 

Thus it remains to verify \eqref{real5G.5}. To this end fix a $0\le k<\mu$ and let 
$\underline{s}^\prime$ be a projector onto $E^\prime:=\text{\rm ker}\,\ulS_{kk}$. 
By construction of $\ulS_0$, the bundle $E^\prime_0:=\text{\rm ker}\,\ulS_{0,kk}$ is constant 
in $t$ near the singularity (i.e. the fibre at $(t,x)$ only depends on $x$), and $E^\prime_0$ 
coincides with $E^\prime$ in $t=0$. Thus, if 
$\underline{s}^\prime_{0,\perp}$ is the orthogonal projection onto $E^\prime_0$, also 
 $$\underline{s}^\prime_0:=\omega\underline{s}^\prime(0,\cdot)+
   (1-\omega)\underline{s}^\prime_{0,\perp}$$
is a projection onto $E^\prime_0$ (note that convex combinations of projections having the 
same image again are projections). Obviously, 
$\underline{s}^\prime-\underline{s}_0^\prime=t\underline{s}_1^\prime$ for a resulting  
$\underline{s}_1^\prime$. 
\end{proof}

\begin{theorem}\label{real511}
If $B$ is as in \eqref{real5G}, then $B:\calH^{0,\gamma}_p(\dz,E)\to\calD(A_{T,\min})$. 
\end{theorem}
\begin{proof}
Since $P(1-\omega)R(1-\omega_1)$ maps to $\calH^{\mu,\infty}_p(\dz,E)_T$ by Lemma 
\ref{real48}.i), it remains to analyze the first summand in \eqref{real5G}, which we now denote 
by $B_1$. By construction, $B_1$ maps into 
$\schnitt\limits_{\delta>0}\,\calH^{\mu,\gamma+\mu-\delta}_p(\dz,E)_T$. 
Next we shall show that $B_1$ maps into the domain of $A_{T,\max}$, i.e.\ 
$ABu\in\calH^{0,\gamma}_p(\dz,E)$ whenever $u\in\calH^{0,\gamma}_p(\dz,E)$. 

Since, near $t=0$, $A-t^{-\mu}\,\op_M(\sigma_M^\mu(A))=tA_1$ for a Fuchs type differential 
operator $A_1$ of order $\mu$, it suffices to show that 
 $$t^{-\mu}\,\op_M(\sigma_M^\mu(A))\,\wt\omega\,B_1\,u\in\calH^{0,\gamma}_p(\dz,E)$$
for some cut-off function $\wt\omega$. Since  
$T=\diag(1,t^{-1},\ldots,t^{-\mu+1})\op_M(\sigma_M^\mu(T))+tT_1$ near $t=0$ for another 
boundary condition $T_1$ and since $\sigma_M^\mu(T)r_0=0$, we obtain 
 $$\wt\omega\,T\Big(\omega\,\op_M^{\gamma+\mu-\varepsilon-\frac{n}{2}}(r_0)\,
   t^\mu\,\omega_0u\Big)\in
   \mathop{\oplus}\limits_{k=0}^{\mu-1}
   \calB^{\mu-k-\frac{1}{p},\gamma+\mu+1-\eps-k-\frac{1}{2}}_p(\bz,F_k)\qquad
   \forall\;\eps>0$$ 
for any $u\in\calH^{0,\gamma}_p(\dz,E)$ and any cut-off function $\wt\omega$ with 
$\wt\omega\prec\omega$. Thus, by Lemma \ref{real48}.iii), 
 $$\wt\omega\,P\omega\,\op_M^{\gamma+\mu-\varepsilon-\frac{n}{2}}(r_0)\,
   t^\mu\,\omega_0u\equiv 
   \wt\omega\,\op_M^{\gamma+\mu-\varepsilon-\frac{n}{2}}(r_0)\,
   t^\mu\,\omega_0u$$
modulo $\calH^{\mu,\gamma+\mu}_p(\dz,E)$. Therefore, if $\wt\omega,\omega_2$ are cut-off 
functions supported sufficiently close to zero with $\omega_2\prec\wt\omega$, 
\begin{align*}
 t^{-\mu}\,\op_M(\sigma_M^\mu(A))\,\wt\omega\,B_1\,u&\equiv
 t^{-\mu}\,\op_M(\sigma_M^\mu(A))\,
 \wt\omega\,\op_M^{\gamma+\mu-\eps-\frac{n}{2}}(r_0)\,t^\mu\,\omega_0u\\
 &\equiv\omega_2\omega_0u+(1-\omega_2)
 t^{-\mu}\,\op_M(\sigma_M^\mu(A))\,
 \wt\omega\,\op_M^{\gamma+\mu-\eps-\frac{n}{2}}(r_0)\,t^\mu\,\omega_0\,u
 \equiv0
\end{align*}
modulo $\calH^{0,\gamma}_p(\dz,E)$. We used that $\sigma_M^\mu(A)r_0=1$. This shows, as 
claimed, that $\text{\rm im}\,B_1\subset\calD(A_{T,\max})$.

Now let $u\in\calH^{0,\gamma}_p(\dz,E)$ and $f\in\calD((A^t)_{\wt T,\max})$. Since then
$f\in\calH^{\mu,-\gamma+\delta}_{p^\prime}(\dz,E)_{\wt T}$ for some $\delta>0$ by Corollary 
\ref{real510}, we have
 $$\skp{A^tf}{B_1u}_{\calH^{0,0}_2(\dz,E)}=
   \skp{f}{AB_1u}_{\calH^{0,0}_2(\dz,E)}=\skp{f}{v}_{\calH^{0,0}_2(\dz,E)}$$
with $v:=AB_1u\in\calH^{0,\gamma}_p(\dz,E)$. Here, we used Corollary \ref{real59} (with 
$\gamma$ replaced by $\gamma-\delta$). Thus  
$\text{\rm im}\,B_1\subset\calD((A^t)_{\wt T,\max}^*)=\calD(A_{T,\min})$ by Theorem 
\ref{real58}.
\end{proof}

\begin{theorem}\label{real514}
\eqref{real5C} can be improved to 
\begin{align*}
 \calD(A_{T,\min})&=
   \Big\{u\in\schnitt_{\eps>0}\calH^{\mu,\gamma+\mu-\eps}_p(\dz,E)_T\st 
   Au\in\calH^{0,\gamma}_p(\dz,E)\Big\}\\
 &=\Big\{u\in\schnitt_{\eps>0}\calH^{\mu,\gamma+\mu-\eps}_p(\dz,E)_T\st 
   t^{-\mu}\,\op_M(\sigma_M^\mu(A))\,\wt\omega\,u\in\calH^{0,\gamma}_p(\dz,E)\Big\},
\end{align*}
where $\wt\omega$ is an arbitrary cut-off function. 
\end{theorem}
\begin{proof}
That the second identity holds is already contained in the previous proof. The independence of 
the choice of $\wt\omega$ is clear. Now let $u$ belong to the right-hand side. 
If $\omega_2\prec\wt\omega$, then $\omega_2T\wt\omega u=0$ implies 
$\omega_2P\wt\omega u-\omega_2u\in\calC^{\infty,\infty}(\dz,E)$ by Lemma \ref{real48}.ii). If 
$P_0$ denotes the projection from Lemma \ref{real513}, then 
$(P-P_0)\wt\omega u\in\calH^{\mu,\gamma+\mu}_p(\dz,E)$. By Remark \ref{real44.5} we may assume 
that $P_0\wt\omega u$ is supported arbitrarily near to the tip. 
Therefore 
 $$t^{-\mu}\,\op_M(\sigma_M^\mu(A))\,P_0\,\wt\omega\,u\equiv
   t^{-\mu}\,\op_M(\sigma_M^\mu(A))(\omega_2u+(1-\omega_2)P\wt\omega u)\equiv0$$
modulo $\calH^{0,\gamma}_p(\dz,E)$ and, if $B$ is as in \eqref{real5G}, then 
 $$B\,t^{-\mu}\,\op_M(\sigma_M^\mu(A))\,P_0\,\wt\omega\,u=
   P\,\omega\,\op_M^{\gamma+\mu-\eps-\frac{n}{2}}(r_0\sigma_M^\mu(A))\,P_0\,\wt\omega\,u=
   P\,\wt\omega\,u+P(P_0-P)\,\wt\omega\,u.$$
The last equality holds, since $r_0\sigma_M^\mu(A)=1-k_0\sigma_M^\mu(T)$ and 
$\op_M(\sigma_M^\mu(T))P_0\wt\omega=T_0P_0\wt\omega=0$. Since $P$ vanishes under $T$ and by 
Lemma \ref{real47}, $P(P_0-P)\,\wt\omega\,u$ belongs to 
$\calH^{\mu,\gamma+\mu}_p(\dz,E)_T\subset\calD(A_{T,\min})$. Then 
 $$u=Pu=P\,\wt\omega\,u+P(1-\wt\omega)\,u\in\calD(A_{T,\min}),$$
since $P(1-\wt\omega)\,u\in\calH^{\mu,\infty}_p(\dz,E)_T$.
\end{proof}

\begin{theorem}\label{real514.5}
Both minimal and maximal extension are Fredholm operators. 
\end{theorem}
\begin{proof}
Let $B$ as in \eqref{real5G}. Since $B$ maps into the domain of $A_{T,\min}$, we have 
$A_{T,\text{\rm min/max}}B=AB$ on $\calH^{0,\gamma}_p(\dz,E)$. Write  
\begin{equation}\label{49a}
 AB=A\,P\,\omega\,\op_M^{\gamma+\mu-\eps-\frac{n}{2}}(r_0)\,t^\mu\,\omega_0+
 A\,P\,(1-\omega)\,R\,(1-\omega_1).
\end{equation}
Let $\omega_2\prec\omega_3\prec\omega_4\prec\omega$. Then, if $[\cdot,\cdot]$ denotes the 
commutator, 
 $$[AP,1-\omega]=\omega_3[\omega,AP]\omega_4+(1-\omega_3)[\omega,AP](1-\omega_2)+
   \omega_3[\omega,AP](1-\omega_4)+(1-\omega_3)[\omega,AP]\omega_2.$$
The first summand vanishes, the second is an operator of order $\mu-1$ located away from the 
singularity, the third and fourth map into $\calC^{\infty,\infty}(\dz,E)$ by Lemma 
\ref{real48}.i) and the locality of $A$. In particular, $[AP,\omega]\in\calK$, the compact 
operators in $\calH^{0,\gamma}_p(\dz,E)$. Using that $PR=R$ and that $R$ inverts $A$ up to a 
Green remainder, we obtain 
 $$A\,P\,(1-\omega)\,R\,(1-\omega_1)\equiv1-\omega\quad\text{modulo }\calK.$$
To treat the first summand in \eqref{49a}, note that, if $\omega$ is chosen to be supported 
sufficiently close to zero, then each of the operators $P\omega$, $P_0\omega$, $\omega P$, and 
$\omega P_0$ can be assumed to be located arbitrarily close to the singularity, cf.\ Remark 
\ref{real44.5}. In view of the compactness of the commutator above, we get 
 $$A\,P\,\omega\,\op_M^{\gamma+\mu-\eps-\frac{n}{2}}(r_0)\,t^\mu\,\omega_0\equiv 
   \omega\,A\,P_0\,\op_M^{\gamma+\mu-\eps-\frac{n}{2}}(r_0)\,t^\mu\,\omega_0+
   \omega\,A\,(P-P_0)\,\op_M^{\gamma+\mu-\eps-\frac{n}{2}}(r_0)\,t^\mu\,\omega_0$$
modulo $\calK$. In the first term, $P_0$ can be dropped, since 
$\op_M(\sigma_M^\mu(T))\op_M^{\gamma+\mu-\eps-\frac{n}{2}}(r_0)=0$. Writing 
$A=t^{-\mu}\op_M(\sigma_M^\mu(A))+tA_1$ near $t=0$ for a $\mu$-th order Fuchs type operator 
$A_1$, 
 $$AB\equiv1+\omega\,t\,A_1\op_M^{\gamma+\mu-\eps-\frac{n}{2}}(r_0)\,t^\mu\,\omega_0+
   \omega\,A\,(P-P_0)\,\op_M^{\gamma+\mu-\eps-\frac{n}{2}}(r_0)\,t^\mu\,\omega_0=:1+C$$
modulo $\calK$. Due to the presence of the $t$-factor and Lemma \ref{real513}, 
$\|C\|_{\calL(\calH^{0,\gamma}_p(\dz,E))}<\frac{1}{2}$ if $\omega$ is located sufficiently close 
to 0 (since the $\sup$-norm of $t\wt{\omega}$ tends to zero if the support of $\wt{\omega}$ 
does), and $1+C$ is invertible. Altogether, we found a $B_R$ such that 
$A_{T,\text{\rm min/max}}B_R=1+C_R$ for a $C_R\in\calK$. This yields that 
$\text{\rm im}\,A_{T,\text{\rm min/max}}$ is finite codimensional, since the image of the 
Fredholm operator $1+C_R$ is. Analogously one can construct a right-inverse modulo compact 
operators for $A^t_{\wt T,\text{\rm min/max}}$. Passing to the adjoint yields a $B_L$ such that 
$B_L A_{T,\text{\rm min/max}}=1+C_L$ on $\calD(A_{T,\text{\rm min/max}})$ with a $C_L\in\calK$. 
This shows that $\text{\rm ker}\,A_{T,\text{\rm min/max}}$ is finite dimensional, since the 
kernel of $1+C_L$ is. 
\end{proof}

We conclude this section with an observation that will have some importance in describing the 
maximal domain. 

\begin{lemma}\label{real516}
Let $T$ and $T_0$ be as in Lemma \text{\rm\ref{real512}} and $P$, $P_0$ be as in Lemma 
\text{\rm\ref{real513}}. Let $0<\delta<1$. Then $P_0$ induces an isomorphism 
 $$\calD(A_{T,\max})\cap\calH^{\mu,\gamma+\mu-\delta}_p(\dz,E)\Big/\calD(A_{T,\min})
   \longrightarrow
   \calD(A_{T_0,\max})\cap\calH^{\mu,\gamma+\mu-\delta}_p(\dz,E)\Big/\calD(A_{T_0,\min}).$$
Its inverse is induced by $P$. 
\end{lemma}
\begin{proof}
We use the description of minimal domains established in Theorem \ref{real514}. If 
$u\in\calD(A_{T,\min})$, then 
$P_0u\in\schnitt\limits_{\eps>0}\calH^{\mu,\gamma+\mu-\eps}_p(\dz,E)_{T_0}$ by 
\eqref{real4G}. Moreover, $AP_0u=Au+A(P_0-P)u\in\calH^{0,\gamma}_p(\dz,E)$, since $Pu=u$ and 
by Lemma \ref{real513}. Thus $P_0:\calD(A_{T,\min})\to\calD(A_{T_0,\min})$ and, analogously 
$P:\calD(A_{T_0,\min})\to\calD(A_{T,\min})$. 
If $u\in\calD(A_{T,\max})\cap\calH^{\mu,\gamma+\mu-\delta}_p(\dz,E)$ the same argument shows 
that 
$P_0u\in\calD(A_{T_0,\max})\cap\calH^{\mu,\gamma+\mu-\delta}_p(\dz,E)$. Similarly for $P$. 
Hence the maps $[u]\mapsto[P_0u]$ and $[v]\mapsto[Pv]$ are well-defined. 
Moreover, 
 $$[PP_0u]=[P^2u+P(P_0-P)u]=[u],$$
since $P^2u=Pu=u$  and $P(P_0-P)u\in\calH^{\mu,\gamma+\mu}_p(\dz,E)_T$ by Lemma \ref{real513}. 
Analogously, $[P_0Pv]=[v]$. 
\end{proof}


\section{The maximal extension}\label{section6}

In this section we shall discuss the maximal extension of the operator $A$ in \eqref{real5B}. 
As before we shall assume normality and $\dz$-ellipticity. 
Let $T_0$ be as in Lemma \ref{real512} and define 
\begin{equation}\label{real6A}
 A_0=\omega\,t^{-\mu}\,\op_M(\sigma^\mu_M(A))+(1-\omega)\,A.
\end{equation}
Then $A-A_0=tA_1$ with a Fuchs type differential operator $A_1$ of order $\mu$. 
Arguing as in the proof of  Lemma \ref{real512}, we see that 
the function $\omega$ can be chosen in such a way that both 
$\begin{pmatrix}A\\T_0\end{pmatrix}$ and 
$\begin{pmatrix}A_0\\T_0\end{pmatrix}$ are $\dz$-elliptic. 
Observe that both $A_0$ and $T_0$ have $t$-independent coefficients near $t=0$. 
We choose projections $P$ and $P_0$ for $T$ and $T_0$ as in Lemma \ref{real513}. 

\subsection{A class of smoothing operators}\label{section6.1}
Let us recall basic properties of operators of the form
\begin{equation*}
 G=\omega\left(\op_M^{\gamma_1-\frac{n}{2}}(r)-
 \op_M^{\gamma_2-\frac{n}{2}}(r)\right):
 \cicomp(\rz_+\times X,E_0)\longrightarrow\ci(\intd,E),
\end{equation*}
where $\gamma_1<\gamma_2$ and 
$r\in M^{\mu,0}_O(X;E_0;E_0)+M^{-\infty,0}_P(X;E_0;E_0)$ is a meromorphic Mellin symbol 
of type 0 with asymptotic type $P$ as described in Section \ref{section3}. 
An analysis of such operators for closed $X$ can be found in \cite{Lesc}; 
the results in the present case with boundary are completely analogous. 
Let 
\begin{equation*}
 \smsum_{k=0}^{n_p}R_{pk}(z-p)^{-(k+1)},\qquad 
 R_{pk}\in N_p\subset \calB^{-\infty,0}(X;E_0;E_0),
\end{equation*}
denote the principal part of $r$ around $p$ for $(p,n_p,N_p)\in P$. The $R_{pk}$ have finite 
rank by definition. Thus, also each of the mappings 
 $$C_p:\mathop{\oplus}\limits_{j=0}^{n_p}\,\ci(X,E_0)\longrightarrow
   \mathop{\oplus}\limits_{k=0}^{n_p}\,\ci(X,E_0)$$
defined by the left-upper triangular block-matrix $(c_{jk}^p)_{0\le j,k\le n_p}$ with 
$c_{jk}^p=R_{p,j+k}$ if $j+k\le n_p$ and $c_{jk}^p=0$ otherwise, is of finite rank. Let us set 
$M(r,p):=\text{\rm rank}\,C_p$ and $M(r,p)=0$ if $r$ is holomorphic in $p$.

\begin{proposition}\label{real60}
If $G$ is as above, then $G$ is of finite rank and, for each $u\in\cicomp(\rz_+\times X,E_0)$,   
 $$(Gu)(t,x)=\omega(t)
   \smsum_{\substack{(p,n_p,N_p)\in\,P\\ 
           -\gamma_2<\text{\rm Re}\,p-\frac{n+1}{2}<-\gamma_1}}
   \smsum^{n_p}_{l=0}\,\zeta_{pl}(u)(x)\,t^{-p}(\log t)^l$$
with the linear maps $\zeta_{pl}:\cicomp(\rz_+\times X,E_0)\to
\text{\rm im}\,R_{pl}+\ldots+\text{\rm im}\,R_{pn_p}\subset\ci(X,E_0)$ given by 
 $$\zeta_{pl}(u)(x)=\smsum_{k=l}^{n_p}\,\frac{(-1)^{l}}{l! (k-l)!}\,
   R_{pk}\,\frac{\partial^{k-l}}{\partial z^{k-l}}(\calM u)(p,x),$$
where $\calM=\calM_{t\to z}$ denotes the Mellin transform. The rank is 
 $$\text{\rm rank}\,G=
   \smsum_{-\gamma_2<\text{\rm Re}\,z-\frac{n+1}{2}<-\gamma_1}M(r,z).$$
\end{proposition}

Note that $G$ only depends on the meromorphic structure of $r$; 
a change by a holomorphic symbol leaves $G$ invariant. 
Changing the domain of $G$ to be 
$\calH^{0,\gamma_1}_p(\rz_+\times X,E_0)\cap\calH^{0,\gamma_2}_p(\rz_+\times X,E_0)$ does not 
alter the image.

\subsection{The maximal domain in case of $t$-independent coefficients}\label{section6.2}

Let $\begin{pmatrix}r_0&k_0\end{pmatrix}$ be the inverted conormal symbol of 
$\begin{pmatrix}A\\T\end{pmatrix}$ and thus of $\begin{pmatrix}A_0\\T_0\end{pmatrix}$. 

\begin{proposition}\label{real61}
Let $\eps>0$ be so small that $r_0$ has no pole in the strips 
$-\gamma-\mu<\re z-\frac{n+1}{2}\le-\gamma-\mu+\eps$ and 
$-\gamma-\eps\le\re z-\frac{n+1}{2}<-\gamma$. 
Define
 $$G_0=\omega\left\{\op_M^{\gamma+\eps-\frac{n}{2}}(r_0)-
        \op_M^{\gamma+\mu-\eps-\frac{n}{2}}(r_0)\right\}$$
for an arbitrary fixed cut-off function $\omega$. Then
\begin{itemize}                     
\item[a)] $\text{\rm Im}\,P_0\,G_0
   \,\cap\,\calH^{\mu,\gamma+\mu-\delta}_p(\dz,E)=\{0\}$, if $\delta>0$ is sufficiently small.
Moreover, $\text{\rm im}\,P_0\,G_0\subset\calD((A_0)_{T_0,\max})$ and
 $$\text{\rm rank}\,P_0\,G_0=\text{\rm rank}\,G_0=
   \smsum_{-\gamma-\mu<\text{\rm Re}\,z-\frac{n+1}{2}<-\gamma}M(r_0,z).$$
\item[b)] $\calD((A_0)_{T_0,\max})=\calD((A_0)_{T_0,\min})\oplus\text{\rm im}\,P_0\,G_0$. 
\end{itemize}
\end{proposition}

\begin{proof}
a) 
Since $\sigma^\mu_M(T_0)r_0=0$ and $T_0$ is a differential operator 
with $t$-independent coefficients,
$T_0G_0u=0$ near $t=0$ for $u\in \cicomp(\rz_+\times X,E_0)$. 
By Lemma \ref{real48}.i), $P_0G_0u \equiv  G_0u$ modulo $\calC^{\infty,\infty}(\dz,E)$.
Hence also $A_0P_0G_0u\in \calC^{\infty,\infty}(\dz,E)$,
since $A_0G_0u=0$ near $t=0$.
Moreover, since the image of $G_0$ has trivial intersection with  
$\calH^{\mu,\gamma+\mu-\delta}_p(\dz,E)$
for small $\delta >0$ by Proposition \ref{real60}, we conclude that
also $\text{\rm im} P_0G_0$ has trivial intersection
and that the two ranks coincide. 

b) The directness of the sum and the inclusion `$\supset$' follow from 
a). So let $u\in\calD((A_0)_{T_0,\max})$ and let $\omega_0$ 
be a cut-off function with $\omega_0\prec\omega$. 
Then $P_0(1-\omega_0)u\in\calH^{\mu,\infty}_p(\dz,E)$ 
by Corollary \ref{real510} and Lemma \ref{real48}.i). 
Hence $P_0\omega_0u=u-P_0(1-\omega_0)u\in\calD((A_0)_{T_0,\max})$, 
and thus $A_0P_0\omega_0u\in\calH^{0,\gamma}_p(\dz,E)$. 
As $A_0$ has $t$-independent  coefficients and $P_0\omega_0 u$ is 
supported close to the singularity, we can write
\begin{align*}
 P_0\,\omega_0\,u&=P_0\,\omega\,\op_M^{\gamma+\eps-\frac{n}{2}}(r_0)\,t^\mu\,
    A_0\,P_0\,\omega_0\,u\\
 &=P_0\,\omega\,\op_M^{\gamma+\mu-\eps-\frac{n}{2}}(r_0)(t^\mu\,A_0\,P_0\,\omega_0\,u)+
 P_0\,G_0(t^\mu\,A_0\,P_0\,\omega_0\,u),
\end{align*}
which is in $\calD((A_0)_{T_0,\min})\oplus\text{\rm im}\,P_0\,G_0$ by 
Theorem \ref{real511}.
\end{proof}

\subsection{On the index of the minimal and the maximal extension}\label{section6.3}

\begin{lemma}\label{real63}
For sufficiently small $\eps>0$, 
 $$\calD(A_{T,\max})\cap\calH^{\mu,\gamma+\mu-\eps}_p(\dz,E)=\calD(A_{T,\min}).$$
\end{lemma}
\begin{proof}
The identity $A-A_0=tA_1$ implies that
 $$\calD(A_{T_0,\max})\cap\calH^{\mu,\gamma+\mu-\eps}_p(\dz,E)=
   \calD((A_0)_{T_0,\max})\cap\calH^{\mu,\gamma+\mu-\eps}_p(\dz,E),\quad 0<\eps<1.$$
For small $\eps$,  Proposition \ref{real61} implies that the right hand side 
equals $\calD((A_0)_{T_0,\min})$ which in turn equals  $\calD(A_{T_0,\min})$
by Theorem \ref{real514}. 
Now the assertion follows from Lemma \ref{real516}. 
\end{proof}

The next proposition says that for considerations of the index one may assume that 
$\begin{pmatrix}A\\T\end{pmatrix}$ is elliptic with respect to $\gamma+\mu$, and thus the 
domain of $A_{T,\min}$ coincides with the space $\calH^{\mu,\gamma+\mu}_p(\dz,E)_T$. This 
observation was made in \cite{GLM} for operators on conic manifolds without boundary. 

\begin{proposition}\label{real64}
Let $\delta>0$ and let $\wt A$ denote the unbounded operator  
 $$\cii(\dz,E)_T\subset\calH^{0,\gamma-\delta}_p(\dz,E)\longrightarrow 
   \calH^{0,\gamma-\delta}_p(\dz,E),$$
given by the action of $A$. 
Both $A_{T,\min}$ and $\wt A_{T,\min}$ are Fredholm operators according to Theorem 
\text{\rm\ref{real514.5}}.
Moreover, their indices coincide, provided $\delta$ is small enough. 
\end{proposition}
\begin{proof}
By Theorem \ref{real514}, 
$\text{\rm ker}\,A_{T,\min}\subset\text{\rm ker}\,\wt A_{T,\min}$. 
If $u\in\text{\rm ker}\,\wt A_{T,\min}$, then 
 $$u\in\calD(A_{T,\max})\cap\calH^{\mu,\gamma+\mu-\delta}_p(\dz,E)=\calD(A_{T,\min})$$ 
for small $\delta$ according to Lemma \ref{real63}. Thus 
$\text{\rm ker}\,A_{T,\min}=\text{\rm ker}\,\wt A_{T,\min}$. 
Taking adjoints of the minimal extensions leads to $A^t_{\wt T,\max}$ and 
$\wt A^t_{\wt T,\max}$ acting as unbounded operators in 
$\calH^{0,-\gamma}_{p^\prime}(\dz,E)$ and $\calH^{0,-\gamma+\delta}_{p^\prime}(\dz,E)$, 
respectively, by Theorem \ref{real58}. 
Obviously $\text{\rm ker}\,\wt A^t_{\wt T,\max}\subset\text{\rm ker}\,A^t_{\wt T,\max}$. 
By Corollary \ref{real510} the reverse inclusion also holds for small $\delta$. 
Hence 
 $\text{\rm ind}\,\wt A_{T,\min}=\text{\rm dim}\,\text{\rm ker}\,\wt A_{T,\min}-
   \text{\rm dim}\,\text{\rm ker}\,\wt A^*_{T,\min}=\text{\rm ind}\,A_{T,\min}.$
\end{proof}

\begin{theorem}\label{real65}
 $\text{\rm ind}\,(A_0)_{T_0,\min}=\text{\rm ind}\,A_{T,\min}$. 
\end{theorem}

\begin{proof}
In view of Proposition \ref{real64} we may assume ellipticity with respect to $\gamma+\mu$. 
By Theorem \ref{real81} of the Appendix, the statement is equivalent to 
$\text{\rm ind}\,\begin{pmatrix}A\\T\end{pmatrix}=
\text{\rm ind}\,\begin{pmatrix}A_0\\T_0\end{pmatrix}$. However, 
$\begin{pmatrix}A\\T\end{pmatrix}$ and $\begin{pmatrix}A_0\\T_0\end{pmatrix}$ 
can be connected by a homotopy through elliptic elements, hence have the same index: 
In fact, if $A=t^{-\mu}\,\op_M(h)$ near $t=0$, let $h_\eps(t,z)=h(\eps t,z)$ and set 
 $$A_\eps=\omega\,t^{-\mu}\,\op_M(h_\eps)+(1-\omega)A,\qquad 0\le\eps\le 1,$$
with $\omega$ supported sufficiently close to 0. Similarly we define $T_\eps$. 
\end{proof}

\begin{corollary}\label{real66}
The domain of the maximal extension differs from the domain of the minimal 
extension by a finite-dimensional space. 
More precisely,
 $$\text{\rm dim}\,\calD(A_{T,\max})\Big/\calD(A_{T,\min})=\text{\rm dim}\,
   \calD((A_0)_{T_0,\max})\Big/\calD((A_0)_{T_0,\min})=
   \smsum_{-\gamma-\mu<\re z-\frac{n+1}{2}<-\gamma}\,M(r_0,z).$$
\end{corollary}

\begin{proof}
Let $\iota:\calD(A_{T,\min})\hookrightarrow\calD(A_{T,\max})$ denote the inclusion, so that  
$A_{T,\min}=A_{T,\max}\iota$. 
By Theorem \ref{real514.5}, $\iota$ is a Fredholm operator, and
 $$\text{\rm dim}\,\calD(A_{T,\max})\Big/\calD(A_{T,\min})=-\text{\rm ind}\,\iota=
   \text{\rm ind}\,A_{T,\max}-\text{\rm ind}\,A_{T,\min}=
   -\text{\rm ind}\,(A^t)_{\wt T,\min}-\text{\rm ind}\,A_{T,\min},$$
where the last identity follows from Theorem \ref{real510.5}. Next we observe that the adjoint 
boundary condition  $(T_0)\,\wt{}$ can be constructed in such a way that it 
coincides with $(\wt T)_0$ near $t=0$. From Theorem \ref{real65} we obtain 
\begin{align*}
   \text{\rm dim}\,\calD(A_{T,\max})\Big/\calD(A_{T,\min})&=
   -\text{\rm ind}\,(A_0^t)_{\wt T_0,\min}-\text{\rm ind}\,(A_0)_{T_0,\min}\\
   &=\text{\rm dim}\,\calD((A_0)_{T_0,\max})\Big/\calD((A_0)_{T_0,\min}),
\end{align*}
and the assertion follows from Proposition \ref{real61}. 
\end{proof}

\subsection{Lower order conormal symbols}\label{section6.3a}
In order to describe the maximal domain for operators with $t$-dependent coefficients we have 
to take into 
account the first $\mu$ Taylor coefficients of the Mellin symbols of $A$ and $T$ at $t=0$. 
Writing $A$ as in \eqref{real4A.5}, we let 
 $$f_j=\frac1{j!}\sum_{k=0}^\mu \frac{d^ja_k(0)}{dt^j}z^k\in M_O^{\mu,0}(X,E_0), 
   \quad j=0,\ldots,\mu-1,$$
so that, near $t=0$ ,
\begin{equation}\label{real6D} 
 A=t^{-\mu}\smsum_{j=0}^{\mu-1}\,t^{j}\,\op_M(f_j)+t^\mu\,A_\mu
\end{equation}
with a Fuchs type differential operator $A_\mu$ of order $\mu$.
Similarly, writing the boundary condition $T=S\varrho$ near $t=0$ in the form 
$T=\diag(1,t^{-1},\ldots,t^{-\mu+1})\,\op_M(s)\varrho$, we let 
 $$s_j=\frac1{j!}\, \frac{d^js(0)}{dt^j}\varrho,\quad j=0,\ldots,\mu-1,$$
so that, near $t=0$,
\begin{equation}\label{real6E} 
 T=\diag(1,t^{-1},\ldots,t^{-\mu+1})\,
   \smsum_{j=0}^{\mu-1}\,t^j\,\op_M(s_j)+t^\mu\,T_\mu.
\end{equation}
with a resulting boundary condition $T_\mu$. We then call 
 $$\sigma_M^{\mu-j}\left(\begin{pmatrix}A\\T\end{pmatrix}\right)=
   \begin{pmatrix}f_j\\s_j\end{pmatrix},\quad j=0,\ldots,\mu-1,$$
the conormal symbol of order $\mu-j$.

\subsection{Description of the maximal domain}\label{section6.4}
Following the method developed in \cite{ScSe2} for the boundaryless case,  
we define recursively 
$\begin{pmatrix}r_j&k_j\end{pmatrix}$, $0\le j<\mu$, by 
\begin{align*}
\begin{pmatrix}r_0&k_0\end{pmatrix}&=\begin{pmatrix}f_0\\s_0\end{pmatrix}^{-1}\\
\begin{pmatrix}r_j&k_j\end{pmatrix}&=-\begin{pmatrix}T^{-j}r_0&T^{-j}k_0\end{pmatrix}
  \smsum_{l=0}^{j-1}\,\begin{pmatrix}T^{-l}f_{j-l}\\T^{-l}s_{j-l}\end{pmatrix}
  \begin{pmatrix}r_l&k_l\end{pmatrix},\qquad j=1,\ldots,\mu-1.
\end{align*}
Recall that $T^k$ denotes the operator of shifting the argument by $k$, i.e.\ 
$(T^kg)(z)=g(z+k)$. These equations are equivalent to 
 $$\smsum_{j=0}^l\,\begin{pmatrix}T^{-j}f_{l-j}\\T^{-j}s_{l-j}\end{pmatrix}
   \begin{pmatrix}r_j&k_j\end{pmatrix}=\delta_{0l}\begin{pmatrix}1&0\\0&1\end{pmatrix}\qquad
   \forall\;0\le l<\mu,$$
with the Kronecker symbol $\delta_{jk}$. Written componentwise, we particularly obtain
\begin{equation}\label{real6F} 
 \smsum_{j=0}^l\,T^{-j}f_{l-j}\,r_j=\delta_{0l},\qquad
 \smsum_{j=0}^l\,T^{-j}s_{l-j}\,r_j=0\qquad \forall\;0\le l<\mu.
\end{equation} 

\begin{definition}\label{real67}
We define the operators 
$G_k=\sum\limits_{l=0}^kG_{kl}:
\cicomp(\rz_+\times X,E_0)\longrightarrow\calC^{\infty,\gamma}(\dz,E)$, $0\le k<\mu$,  
by
 $$G_{00}=\omega\left(\op_M^{\gamma+\mu+\varepsilon-1-\frac{n}{2}}(r_0)-
   \op_M^{\gamma+\mu-\varepsilon-\frac{n}{2}}(r_0)\right),$$
and, if $1\le k\le\mu-1$ and $0\le l\le k$, 
 $$G_{kl}=\omega\,t^l
   \left(\op_M^{\gamma+\mu+\varepsilon-k-1-\frac{n}{2}}(r_l)-
   \op_M^{\gamma+\mu+\varepsilon-k-\frac{n}{2}}(r_l)\right).$$
Here, $\omega$ is an arbitrary fixed cut-off function. Moreover, $\eps>0$ is fixed and so 
small that for each pole $p$ of one of the symbols $r_0,\ldots,r_\mu-1$ all the distances 
$|\frac{n+1}{2}-\gamma-\mu+j-\re p|$, $j=0,\ldots,\mu-1$, are either zero or $>\eps$.
\end{definition} 
 
Note that, by construction,  $G_{kl}$  maps into 
$\calC^{\infty,\gamma+\mu-k+l-1+\eps}(\dz,E)\hookrightarrow
\calC^{\infty,\gamma+\eps}(\dz,E)$. 

\begin{lemma}\label{real68}
For $k=0\,\ldots,\mu-1$, 
 $$A:\text{\rm im}\,G_k\longrightarrow\calC^{\infty,\gamma+\eps}(\dz,E),\qquad
   T:\text{\rm im}\,G_k\longrightarrow
   \mathop{\oplus}\limits_{j=0}^{\mu-1}\calC^{\infty,\gamma+\mu+\eps-j-\frac{1}{2}}(\bz,F_j).$$
\end{lemma}

\begin{proof}
Let $u\in\cicomp(\rz_+\times X,E_0)$. Writing  $A$ as in \eqref{real6D} and using the mapping 
properties of 
the $G_{kl}$, we see that $AG_ku\in\calC^{\infty,\gamma+\eps}(\dz,E)$ if and only if
 $$\widetilde{\omega}\,\smsum_{j=0}^k\smsum_{l=0}^{k-j}
   t^j\,\op_M(f_j)\,G_{kl}\,u\;\in\;\calC^{\infty,\gamma+\mu+\eps}(\dz,E)$$
for every cut-off function $\widetilde{\omega}$. Choosing $\widetilde{\omega}$ with
$\wt{\omega}\prec\omega$, inserting the definition of $G_{kl}$ and 
rearranging the order of summation, this is equivalent to  
 $$\widetilde{\omega}\,\smsum_{j=0}^kt^j\smsum_{l=0}^j\left(
   \op_M^{\gamma+\mu+\varepsilon-k-1-\frac{n}{2}}((T^{-l}f_{j-l})r_l)-
   \op_M^{\gamma+\mu+\varepsilon-k-\frac{n}{2}}((T^{-l}f_{j-l})r_l)\right)u
   \in\calC^{\infty,\gamma+\mu+\eps}(\dz,E).$$
However, this expression actually equals zero by \eqref{real6F}. For $T$ we argue in the same 
way, using the representation \eqref{real6E}. 
\end{proof}

\begin{corollary}\label{real69}
$\text{\rm Im}\,PG_k\subset\calD(A_{T,\max})$ for  $ k=0,\ldots,\mu-1$. 
\end{corollary}
\begin{proof}
The mapping property of $T$ in Lemma \ref{real68} in connection with Lemma \ref{real48}.iii)
implies that $PG_ku\equiv G_ku$ modulo $\calC^{\infty,\gamma+\mu+\eps}(\dz,E)$ 
for $u\in\cicomp(\rz_+\times X,E_0)$. 
Therefore $APG_ku\in\calC^{\infty,\gamma+\eps}(\dz,E)$ 
according to Lemma \ref{real68}. 
\end{proof}

\begin{definition}\label{real610}
We set $\calE:=\text{\rm im}\,PG_{0}+\ldots+\text{\rm im}\,PG_{\mu-1}\subset\calD(A_{T,\max})$.
\end{definition}

 $\calE$ is a finite-dimensional subspace of $\calC^{\infty,\gamma+\varepsilon}(\dz,E)_T$. 
As in the proof of Corollary \ref{real69} we see that
\begin{equation}\label{real6F.5}
 \calE+\calC^{\infty,\gamma+\mu+\varepsilon}(\dz,E)=
 \text{\rm im}\,G_{0}+\ldots+\text{\rm im}\,G_{\mu-1}+
 \calC^{\infty,\gamma+\mu+\varepsilon}(\dz,E).
\end{equation}
In particular,  there exist complex numbers $q_j$, determined from the meromorphic structure 
of the symbols $r_l$, with 
\begin{equation}\label{real6G}
 \frac{n+1}{2}-\gamma-\mu\le\re q_j<\frac{n+1}{2}-\gamma 
\end{equation}
such that any element $u$ of $\calE$ is of the form 
\begin{equation}\label{real6H}
 u(t,x)\equiv\omega(t)\,\smsum_{j=0}^N\smsum_{k=0}^{l_j} 
        u_{jk}(x)\,t^{-q_j}\,\log^k t\qquad\text{\rm mod }
        \calC^{\infty,\gamma+\mu+\eps}(\dz,E)
\end{equation}
with smooth functions $u_{jk}\in\ci(X,E_0)$. 

In case $A=A_0$ and $T=T_0$ have $t$-independent coefficients near $t=0$, 
all $G_{kl}$, $l\ge1$, vanish and   
 $$\calE=\calE_0=
   \mbox{\rm im}\,P_0\,\omega\left(
   \op_M^{\gamma+\mu-\varepsilon-\frac{n}{2}}(r_0)-
   \op_M^{\gamma+\varepsilon-\frac{n}{2}}(r_0)\right),$$
cf.\ Proposition \ref{real61}. In particular, we have twice strict inequality `$<$' in
\eqref{real6G}. In general, equality in \eqref{real6G} is possible. 

\begin{proposition}\label{real611}
The dimension of $\calE$ can be estimated by 
 $$\dim\,\calE\ge\dim\,\calE_0=\dim\,\text{\rm im}\,\omega
     \left(\op_M^{\gamma+\varepsilon-\frac{n}{2}}(r_0)-
     \op_M^{\gamma+\mu-\varepsilon-\frac{n}{2}}(r_0)\right).$$
Moreover, there exists a subspace of $\calE$ which has the same dimension as $\calE_0$ and has 
trivial intersection with $\calD(A_{T,\min})$. 
\end{proposition}
\begin{proof}
Proposition \ref{real61} shows the identity for $\dim\,\calE_0$. 
For $0\le k<\mu$ pick $u_{k1},\ldots,u_{kn_k}\in\cicomp(\rz_+\times X,E_0)$ 
in such a way that $\{G_{k0}u_{kj}\st 1\le j\le n_k\}$ is a basis of $\mbox{\rm im}\,G_{k0}$. 
Then $n_1+\ldots+n_{\mu-1}=\dim\,\calE_0$. We shall show that 
 $$\{PG_ku_{kj}\st 0\le k<\mu,\, 1\le j\le n_k\}\subset\calE$$
is a set of linearly independent functions with 
 $$\mbox{\rm span}\{PG_ku_{kj}\st 0\le k<\mu,\, 1\le j\le n_k\}\cap
   \calD(A_{T,\min})=\{0\}:$$
Let $\alpha_{jk}\in\cz$ and 
$\sum\limits_{k=0}^{\mu-1}\sum\limits_{j=1}^{n_k}\alpha_{jk}PG_ku_{kj}\in\calD(A_{T,\min})$. 
Since $PG_ku_{kj}\equiv G_ku_{kj}$ modulo $\calC^{\infty,\gamma+\mu}(\dz,E)$, cf.\ the proof of  
Corollary \ref{real69}, we obtain that
 $$\smsum_{k=0}^{\mu-1}\smsum_{j=1}^{n_k}\alpha_{jk}G_ku_{kj}=u\;\in\;
   \calH^{0,\gamma+\mu-\varepsilon}_p(\dz,E).$$
Setting $l=\mu-1$, we conclude that 
 $$\smsum_{j=1}^{n_l}\alpha_{jl}G_{l0}u_{lj}=u-
   \smsum_{k=0}^{l-1}\smsum_{j=1}^{n_k}\alpha_{jk}G_ku_{kj}-
   \smsum_{j=1}^{n_l}\alpha_{kl}(G_l-G_{l0})u_{lj}.$$
Now the right-hand side belongs to $\calH^{0,\gamma+\mu-\varepsilon}_p(\dz,E)+
\calH^{0,\gamma+1+\varepsilon}_p(\dz,E)$ which has trivial intersection with 
$\mbox{\rm im}\,G_{l0}$. Hence the left hand side is zero, and 
$\alpha_{jl}=0$ for all $1\le j\le n_l$, since the $G_{l0}u_{lj}$ are linearly independent by 
assumption. Taking $l=\mu-2, \ldots,0$, we see that all $\alpha_{jk}$ must equal zero. 
\end{proof}

\begin{theorem}\label{real612}
With $\calE$ from Definition \text{\rm\ref{real610}}, the domain of the maximal extension is 
 $$\calD(A_{T,\max})=\calD(A_{T,\min})+\calE.$$
The sum is direct at least in the cases where $A$ has $t$-independent coefficients near 
$t=0$ or $r_0$ has no pole on the line $\re z=\frac{n+1}{2}-\gamma-\mu$. 
In any case, 
 $$\calD(A_{T,\min})\cap\calE\subset
   \mbox{\rm im}\,P\,\omega
   \left(\op_M^{\gamma+\mu-\varepsilon-\frac{n}{2}}(r_0)-
     \op_M^{\gamma+\mu+\varepsilon-\frac{n}{2}}(r_0)\right).$$
\end{theorem}
\begin{proof}
By Corollary \ref{real66}, $\dim\,\calD(A_{T,\max})\big/\calD(A_{T,\min})=\dim\,\calE_0$. Thus 
$\calD(A_{T,\max})=\calD(A_{T,\min})+\calE$ by Proposition \ref{real611}. 
If $A$ has $t$-independent coefficients near $t=0$, the intersection of $\calD(A_{\min})$ and
$\calE=\calE_0$ is zero by Proposition \ref{real61}a). 

In order to see that the sum is also direct if $r_0$ has no pole on the line in question, 
suppose that $u\in\calD(A_{T,\min})$ and $u=Pv$ with 
$v\in\textrm{im}\,G_0+\ldots+\textrm{im}\,G_{\mu-1}$, cf.\ Definition \ref{real610}. 
We can write  
$v=\omega(t)\,\sum\limits_{j=0}^N\sum\limits_{k=0}^{l_j}u_{jk}(x)\,t^{-q_j}\,\log^k t$ 
with $q_j$ satisfying \eqref{real6G}. 
Since $u-v\in\calC^{\infty,\gamma+\mu+\eps}(\dz,E)$, cf.\ the proof of Corollary \ref{real69},
 we have $v\in  
\calH^{\mu,\gamma+\mu-\delta}_p(\dz,E)$ for all $\delta>0$ and 
$Av\in\calH^{0,\gamma}_p(\dz,E)$. 
The first property of $v$ implies that $\re q_j=(n+1)/{2}-\gamma-\mu$ for all $q_j$. 
Now $A$ maps $v$ to a function with a similar structure, where $t$ has exponents 
$-q_j+l$, $l\in\nz_0$. Noting that $\omega \,t^{-p_j}\log^kt\in
\calC^{\infty,\gamma}$ if and only if $\re p_j< (n+1)/2-\gamma$, 
we conclude from the second property that 
$Av\in\calC^{\infty,\gamma+\delta}(\dz,E)$ for all $\delta<1$.  
Moreover, $P_0v\equiv P_0u=u+(P_0-P)u\equiv v$ modulo 
$\calH^{\mu,\gamma+\mu+\eps}_p(\dz,E)$ by \eqref{44a}. 
We next use the facts that we may assume $P_0v$ to be supported near the singularity 
by Remark \ref{real44.5}
and that $r_0f_0+k_0s_0=1$ and $\op_M(s_0)P_0v=T_0P_0v=0$.
This implies that
\begin{equation}\label{real6I} 
 v\equiv P_0v=\wt\omega\op_M^{\gamma+\mu-\eps-\frac{n}{2}}(r_0f_0+k_0s_0)P_0v=
 \wt\omega\op_M^{\gamma+\mu-\eps-\frac{n}{2}}(r_0)\op_M(f_0)P_0v
\end{equation}
modulo $\calH^{\mu,\gamma+\mu+\eps}_p(\dz,E)$ with a suitable cut-off function $\wt\omega$. 
Now $\wt u:=\op_M(f_0)P_0v\equiv\op_M(f_0)v=t^\mu A_0v$ 
modulo $\calH^{0,\gamma+\mu+\eps}_p(\dz,E)$. 
Moreover, $A_0v\in\calC^{\infty,\gamma+\delta}(\dz,E)$, 
since $Av\in\calC^{\infty,\gamma+\delta}(\dz,E)$, $\delta<1$, 
and $A-A_0$ gains a factor $t$. Thus 
$\wt u\in\calH^{0,\gamma+\mu+\eps}_p(\dz,E)$ and with \eqref{real6I} 
 $$v-\omega\left(\op_M^{\gamma+\mu-\varepsilon-\frac{n}{2}}(r_0)-
   \op_M^{\gamma+\mu+\varepsilon-\frac{n}{2}}(r_0)\right)\wt u\equiv
   \wt\omega\op_M^{\gamma+\mu+\eps-\frac{n}{2}}(r_0)\wt u\equiv0$$
modulo $\calH^{\mu,\gamma+\mu+\eps}_p(\dz,E)$ (note that changing $\omega$ to $\wt\omega$ 
induces only a remainder in $\cicomp(\dz,E)$). By the structure of the left-hand side, this 
is only possible if the left-hand side equals 0. This finishes the proof.  
\end{proof}

\begin{remark}\label{real613}
If $\wt\calE=\mbox{\rm span}\{PG_ku_{kj}\st 0\le k<\mu,\, 1\le j\le n_k\}$, as constructed in 
the proof of Proposition {\rm \ref{real611}}, then 
$\calD(A_{T,\max})=\calD(A_{T,\min})\oplus\wt\calE$.
\end{remark}


\section{Ellipticity of the adjoint problem}\label{section7}

Let $\wt\calA=\begin{pmatrix}A^t\\ \wt T\end{pmatrix}$ denote the adjoint problem to  
$\calA=\begin{pmatrix}A\\T\end{pmatrix}$, cf.\ Definition \ref{real55}. 

\begin{proposition}\label{real71}
If $\calA$ is $\dz$-elliptic then so is $\wt\calA$. 
\end{proposition}
\begin{proof}
Clearly with $\sigma_\psi^\mu(\calA)$ and $\wt\sigma_\psi^\mu(\calA)$ also 
$\sigma_\psi^\mu(\wt\calA)$ and $\wt\sigma_\psi^\mu(\wt\calA)$ are invertible (
note that these symbols only involve the differential operators $A$ and $A^t$, and not the 
boundary conditions). 
The fact that invertibility of $\sigma_\partial^\mu(\calA)$ implies that of 
$\sigma_\partial^\mu(\wt\calA)$ is proven in \cite{Grub1}, Theorem 1.6.9. 
This proof passes over to the rescaled boundary symbols. 
\end{proof}

Our next goal is to show that ellipticity of  $\calA$ with respect to $\gamma\in\rz$ 
in the sense of Definition \ref{real41} implies ellipticity of the adjoint problem with 
respect to $-\gamma+\mu$. 
According to Proposition \ref{real71} and Proposition \ref{real42} 
reduces to a question on the invertibility of the conormal symbol. We shall show: 

\begin{proposition}\label{real72}
The conormal symbol of the adjoint problem satisfies 
 $$\sigma_M(\wt\calA)(z)=
   \begin{pmatrix}\sigma_M^\mu(A^t)\\ \sigma_M^\mu(\wt T)\end{pmatrix}(z)=
   \begin{pmatrix}\sigma_M^\mu(A)^t\\ \sigma_M^\mu(T)\,\wt{}\end{pmatrix}
   (n+1-\mu-\overline{z}),\qquad z\in\cz.$$
\end{proposition}

The statement contains some notation we shall now explain: 
For each fixed $z$, 
$\sigma_M^\mu(A)(z)$ is a differential operator on $X$; 
$\sigma_M^\mu(A)^t(z)$ denotes its formal adjoint. 
Moreover 
$\begin{pmatrix}\sigma_M^\mu(A)(z)\\ \sigma_M^\mu(T)(z)\end{pmatrix}$ 
is a boundary value problem on the smooth manifold $X$; 
$\sigma_M^\mu(T)\,\wt{}\,(z)$ denotes the adjoint boundary condition in the 
sense of \ \cite{Grub1}. 
Note that $\sigma_M^\mu(T)(z)$ is a normal boundary condition: 
The entries $\underline{S}_{kk}$ in \eqref{real4C} are Fuchs type operators
of order zero so that the Mellin symbol -- just like the operator -- is 
simply multiplication by the surjective morphism $\underline{S}_{kk}$.

If one goes through the construction (Proposition 1.3.2, Lemma 1.6.1, (1.6.26) in \cite{Grub1} 
and Proposition \ref{real47}, Theorem \ref{real54}, Definition \ref{real55} in this paper) 
one finds that 
 $$\sigma_M^\mu(T)\,\wt{}=I^\times(\sigma_M(\mathfrak{A})\sigma_M(\ulC^\prime))^t,$$
where 
 $$\sigma_M(\ulC^\prime)=(T^k\sigma_M^{j-k}(\ulC^\prime_{jk}))_{0\le j,k<\mu},\qquad 
   \sigma_M(\mathfrak{A})=(T^k\sigma_M^{\mu-j-k-1}(\mathfrak{A}_{jk}))_{0\le j,k<\mu}.$$ 
Now the proof is a purely algebraic calculation, using the rules for computing 
conormal symbols of formally adjoint operators, cf.\ \cite{Schu2}. Here are the details: 

For abbreviation let us write $f^{[t]}(z):=f(n-\overline{z})^t$ and 
$f^{(t)}(z):=f(n+1-\overline{z})^t$. We have 
 $$\wt T=I^\times(\mathfrak{A}C^\prime)^t=(\wt\ulS_{jk})_{j,k}\,
   \diag(1,t^{-1},\ldots,t^{-\mu+1}),\qquad \wt\ulS_{jk}=
   \smsum_l\ulC^{\prime t}_{l,\mu-j-1}\underline{\mathfrak{A}}^t_{kl}.$$
Therefore, $\sigma_M(\wt T)=\big(T^k\sigma_M^{j-k}(\wt\ulS_{jk})\big)_{jk}$ with 
\begin{align*}
 T^k\sigma_M^{j-k}(\wt\ulS_{jk})&=
   \smsum_lT^{\mu-l-1}\sigma_M^{l-\mu+j+1}(\ulC^{\prime t}_{l,\mu-j-1})\,
   T^k\sigma_M^{\mu-k-l-1}(\underline{\mathfrak{A}}^t_{kl})\\
 &=\smsum_lT^j\sigma_M^{l-\mu+j+1}(\ulC^{\prime}_{l,\mu-j-1})^{[t]}\,
   T^{\mu-l-1}\sigma_M^{\mu-k-l-1}(\underline{\mathfrak{A}}_{kl})^{[t]}.
\end{align*}
On the other hand, using $f^{(t)}=T^{-1}f^{[t]}$, we obtain 
\begin{align*}
 I^\times(\sigma_M(\mathfrak{A})\sigma_M(\ulC^\prime))^{(t)}&=
   \Big(\smsum_lT^{j-\mu+1}\sigma_M^{l-\mu+j+1}(\ulC^{\prime}_{l,\mu-j-1})^{(t)}\,
   T^{-l}\sigma_M^{\mu-k-l-1}(\underline{\mathfrak{A}}_{kl})^{(t)}\Big)_{j,k}\\
 &=T^{-\mu}\Big(\smsum_lT^{j}\sigma_M^{l-\mu+j+1}(\ulC^{\prime}_{l,\mu-j-1})^{[t]}\,
   T^{\mu-l-1}\sigma_M^{\mu-k-l-1}(\underline{\mathfrak{A}}_{kl})^{[t]}\Big)_{j,k}\\
 &=T^{-\mu}(T^k\sigma_M^{j-k}(\wt\ulS_{jk}))_{j,k}.
\end{align*}
Hence $\sigma_M(T)\,\wt{}\,(n+1-\overline{z})=\sigma_M(\wt T)(z-\mu)$, and this proves 
Proposition \ref{real72}. 

\begin{corollary}\label{real73}
If $\calA$ is elliptic with respect to $\gamma$ then the adjoint problem is elliptic with 
respect to $-\gamma+\mu$. 
\end{corollary}
\begin{proof}
$\dz$-ellipticity of $\wt\calA$ is shown in Proposition \ref{real71}. By assumption, 
$\sigma_M(\calA)(z)$ is invertible whenever $\re z=\frac{n+1}{2}-\gamma$. This is equivalent, 
by Corollary \ref{real82}, to the invertibility of 
 $$\sigma_M(A)(z):H^\mu_p(X,E_0)_{\sigma_M(T)(z)}\subset L_p(X,E_0)\to L_p(X,E_0).$$
Thus also the adjoint of this operator is invertible. But this adjoint is just given by 
 $$\sigma_M(A)(z)^t:H^\mu_{p^\prime}(X,E_0)_{\sigma_M(T)(z)\,\wt{}}
   \subset L_{p^\prime}(X,E_0)\to L_{p^\prime}(X,E_0).$$
Equivalently, $\begin{pmatrix}\sigma_M(A)^t\\ \sigma_M(T)\,\wt{}\end{pmatrix}(z)$ is invertible 
whenever $\re z=\frac{n+1}{2}-\gamma$. Then, due to Proposition \ref{real72}, also 
$\sigma_M(\wt\calA)(n+1-\mu-\overline{z})$ is invertible. This gives the assertion, since 
$\re(n+1-\mu-\overline{z})=\frac{n+1}{2}+\gamma-\mu$. 
\end{proof}


\section{Appendix: A theorem on Fredholm operators}\label{section8}

Let $E,F,H$ be vector spaces and 
$\calA=\begin{pmatrix}A\\T\end{pmatrix}:H\longrightarrow
\begin{matrix}E\\ \oplus\\F\end{matrix}$
be a linear map. We shall recall here a few known results relating the invertibility and 
Fredholm property of $\calA$ to that of the induced operator 
$A_T:\text{\rm ker}\,T\to E$ defined by the action of $A$. 
Being a Fredholm operator in this context just means that the dimension of the kernel and the 
codimension of the image are finite. Clearly, 
\begin{equation}\label{real8A}
 \text{\rm ker}\,\calA=\text{\rm ker}\,A\,\cap\,\text{\rm ker}\,T=\text{\rm ker}\,A_T.
\end{equation}

\begin{lemma}\label{real80}
$\calA$ is surjective if and only if $A_T$ is surjective and $T:H\to F$ is surjective. 
\end{lemma}
\begin{proof}
Surjectivity of $\calA$ implies that of $A_T$ and $T$, since $E\oplus0$ and $0\oplus F$ are 
contained in the image of $\calA$.   
Vice versa, if $e\oplus f$ is given, there exists a $u_1$ with $Tu_1=f$, and a $u_2$ with 
$A_T u_2=e-A u_1$. Then $\calA u=e\oplus f$ for $u:=u_1+u_2$.  
\end{proof}

Together with \eqref{real8A} we at once obtain the following corollary: 

\begin{corollary}\label{real82}
$\calA$ is invertible if and only if $A_T$ is invertible and $T:H\to F$ is surjective. 
\end{corollary}

\begin{theorem}\label{real81}
$\calA$ is a Fredholm operator if and only if $A_T$ is a Fredholm operator 
and the image of $T$ has finite codimension in $F$. 
In this case, 
 $$\text{\rm ind}\,\calA=
   \text{\rm ind}\,A_T-\text{\rm codim}\,\text{\rm im}\,T.$$
\end{theorem}
\begin{proof}
i) Without loss of generality, we may assume that $T$ is surjective. 
In fact, if 
$F=\text{\rm im}\,T\oplus F_0$ we pass to $\calA_0$ given by the action of $\calA$ on $H$ but 
mapping to $E\oplus\text{\rm im}\,T$. Then obviously $\calA$ is a Fredholm operator if and only 
if $F_0$ is finite dimensional and $\calA_0$ is Fredholm. Then 
$\text{\rm ind}\,\calA_0=\text{\rm ind}\,\calA+\text{\rm dim}\,F_0$, 
and it suffices to show that $\text{\rm ind}\,\calA_0=\text{\rm ind}\,A_T$. 

ii) Let us assume that $\calA$ is a Fredholm operator. 
Let $e_1\oplus f_1,\ldots,e_k\oplus f_k$ be a basis of a
complement to $\text{\rm im}\,\calA$. 
If we define 
 $$K:\cz^k\longrightarrow E,\,c\mapsto\smsum_{j=1}^kc_j\,e_j,\qquad 
   Q:\cz^k\longrightarrow E,\,c\mapsto\smsum_{j=1}^kc_j\,f_j,$$
and set $\wt H=H\oplus\cz^k$, 
$\wt A=\begin{pmatrix}A&K\end{pmatrix}$, and $\wt T=\begin{pmatrix}T&Q\end{pmatrix}$, 
then 
 $\wt\calA:=\begin{pmatrix}\wt A\\ \wt T\end{pmatrix}:\wt H\longrightarrow 
   \begin{matrix}E\\ \oplus\\F\end{matrix}$  
is surjective. According to Lemma \ref{real80}, also 
${\wt A}_{\wt T}:\text{\rm ker}\,\wt T\to E$ is surjective. If $\text{\rm ker}\,\wt T=
\text{\rm ker}\,{\wt A}_{\wt T}\oplus V$ then $\wt A_V:V\to E$, 
defined by the action of $\wt A$, is bijective. 
Let us define the map $S$ by 
 $$S=\begin{pmatrix}S_1\\S_2\end{pmatrix}:E\xrightarrow{{\wt A}_V^{-1}}V\hookrightarrow 
    \wt H=\begin{array}{c}H\\ \oplus\\ \cz^k\end{array}.$$
For an $e\in\text{\rm ker}\,S_2$ we then obtain 
 $$\begin{pmatrix}A\\ T\end{pmatrix}S_1e=
   \begin{pmatrix}A&K\\T&Q\end{pmatrix}\begin{pmatrix}S_1e\\S_2e\end{pmatrix}=
   \begin{pmatrix}\wt A\\ \wt T\end{pmatrix}Se=\begin{pmatrix}e\\0\end{pmatrix},$$
since $S$ maps into the kernel of $\wt{T}$. Thus $\text{\rm ker}\,S_2\subset\text{\rm im}\,A_T$ 
and therefore 
 $$\text{\rm codim im}\,A_T\le\text{\rm codim ker}\,S_2=
   \text{\rm dim im}\,S_2\le k=\text{\rm codim}\,\text{\rm im}\,\calA.$$
Together with \eqref{real8A} this shows that $A_T$ is a Fredholm operator with 
$\text{\rm ind}\,A_T\ge\text{\rm ind}\,\calA$. 

iii) Now let $A_T$ be a Fredholm operator (and $T$ surjective, cf.\ i)). 
Let us write $E=\text{\rm im}\,A_T\oplus E_0$ with a finite dimensional $E_0$, 
and $H=\text{\rm ker}\,T\oplus V$. 
Then $T_V:V\to F$, defined by the action of $T$, is bijective. 
We define the map $S$ by $S:F\xrightarrow{T_V^{-1}}V\hookrightarrow H$ and set 
 $$W=\Big\{\begin{pmatrix}e+ASf\\f\end{pmatrix}\st e\in\text{\rm im}\,A_T,\,f\in F\Big\}.$$
Then $W\subset\text{\rm im}\,\calA$, since for $e=A_T u\in\text{\rm im}\,A_T$ and $f\in F$
 $$\calA(u+Sf)=\begin{pmatrix}A\\T\end{pmatrix}(u+Sf)=
   \begin{pmatrix}e+ASf\\TSf\end{pmatrix}=\begin{pmatrix}e+ASf\\f\end{pmatrix}.$$
The proof if finished if we can show that the codimension of $W$ in $E\oplus F$ equals 
$\text{\rm dim}\,E_0$, for then 
 $$\text{\rm codim im}\,\calA\le\text{\rm codim}\,W=\text{\rm dim}\,E_0=
   \text{\rm codim im}\,A_T;$$
this, together with \eqref{real8A}, yields that $\calA$ is a Fredholm operator with 
$\text{\rm ind}\,\calA\ge\text{\rm ind}\,A_T$. We define 
 $$R=\begin{pmatrix}1&\pi AS&0\\0&1&0\\0&(1-\pi)AS&1\end{pmatrix}:
   \begin{matrix}\text{\rm im}\,A_T\\ \oplus\\F\\ \oplus\\E_0\end{matrix}
   \longrightarrow 
   \begin{matrix}\text{\rm im}\,A_T\\ \oplus\\F\\ \oplus\\E_0\end{matrix},$$
where $\pi$ denotes the projection in $E$ onto $\text{\rm im}\,A_T$ along $E_0$. 
Then $R$ is Fredholm with index 0, since it differs from an invertible operator 
by a finite rank operator (the latter is the $3\times3$-matix with $(1-\pi) AS$ 
as the only non-zero entry). 
Even more is true: $R$ is  invertible, since its kernel is trivial. If 
$\alpha:E\oplus F
=\text{\rm im}\,A_T\oplus E_0\oplus F\to\text{\rm im}\,A_T\oplus F\oplus E_0$ 
is the isomorphism given by exchanging the last two components, 
then $W=(\alpha^{-1}R\alpha)(\text{\rm im}\,A_T\oplus 0\oplus F)$. 
Hence 
$\text{\rm codim}\,W=\text{\rm dim}\,E_0$ as desired. 
\end{proof}


\section{Appendix: The parameter-dependent Boutet de Monvel 
         algebra on smooth manifolds}\label{section2}

We recall some basic facts about Boutet de Monvel's calculus on smooth manifolds \cite{Bout}. 
For details we refer to \cite{Grub1}, also to \cite {ReSc}, \cite{Schr}, \cite{ScSc1} for an 
approach based on operator-valued symbols. In the sequel, $X$ is a smooth compact 
$n$-dimensional manifold with boundary, $s$ is a real number, and $p$ 
is a real number with $1<p<\infty$. 

\subsection{Function spaces}\label{section2.1}

The Sobolev space 
 $$H^s_{p}(\rz^n)= \{ u\in \calS^\prime(\rz^n)\st\spk{D}^s u \in L_p(\rz^n)\},$$ 
carries the norm $\|u\|_{H^s_{p}(\rz^n)} = \| \spk{D}^s u \|_{L_{p}(\rz^n)}$. 
A closed subspace is 
 $$\kringel{H}^s_{p}(\rpbar^n)=\{u\in H^s_{p}(\rz^n)\st
   \text{\rm supp}\,u\subset\rpbar^n\}.$$
Restriction of $H^s_p(\rz^n)$ to the half-space $\rz^n_+$ leads to 
 $$H^s_{p}(\rpbar^n)=r^{+}H^s_{p}(\rz^n),$$ 
endowed with the quotient norm 
$\|u\|_{H^s_{p}(\rpbar^n)}=\inf\{\|v\|_{H^s_{p}(\rz^n)}\st r^+ v=u\}$. 

We write  $B^s_{p}(\rz^n)$ for the Besov space $B^s_{p,p}(\rz^n)$;  
$B^s_{p,q}(\rz^n)$ consists of all tempered distributions $u$ with finite norm 
\begin{equation*}
 \|u\|_{B^s_{p,q}(\rz^n)} =\Big(\|\varphi_{0}(D)u\|^q_{L_{p}(\rz^n)}+
 \smsum_{k=1}^\infty 2^{skq}\| \varphi(2^{-k}D)u\|_{L_{p}(\rz^n)}^q\Big)^{1/q}. 
\end{equation*} 
Here, $\varphi\in\cicomp(\rz^n)$ is chosen such that (i)
$\varphi(\xi) > 0$ for $2^{-1} < |\xi| < 2$ and $\varphi(\xi)=0$ otherwise, (ii)
$\sum\limits_{k=-\infty}^\infty \varphi(2^{-k} \xi) = 1$ when $\xi \not= 0$, and $\varphi_0$ is 
given by $\varphi_{0}(\xi) = 1 - \sum\limits_{k=1}^\infty \varphi(2^{-k} \xi)$. 
Note that $B^s_2(\rz^n)=H^s_2(\rz^n)$. 

\begin{remark}
The $L_2$-scalar product allows an identification of 
$\kringel{H}^{-s}_{p^\prime}(\rpbar^n)$ with the dual of 
$H^s_{p}(\rpbar^n)$, and of ${B}^{-s}_{p^\prime}(\rz^n)$ with the dual of 
$B^s_{p}(\rz^n)$. Here, $p^\prime$ is the dual number to $p$, i.e.\ 
$\frac{1}{p}+\frac{1}{p^\prime}=1$. 
\end{remark}  

\begin{remark}\label{trace}
The trace operator $\gamma_j$, first defined on $\cicomp(\rpbar^n)$ by 
$\gamma_{j}u(x^\prime)=(D^j_{x_{n}}u)(x^\prime,0)$, extends by continuity to 
 $$\gamma_{j}: H^s_{p}(\rpbar^n)\longrightarrow B^{s-j-\frac{1}{p}}_{p}(\rz^{n-1}),\qquad 
   s>j+\frac{1}{p}.$$
\end{remark}

Using a partition of unity we then define  
$H^s_{p}(X)$, $\kringel{H}^s_{p}(X)$, and $B^s_p(X)$ for a manifold $X$ (with or without 
boundary), as well as  $H^s_{p}(X,E)$, $\kringel{H}^s_{p}(X,E)$, and  $B^s_p(X,E)$ for a vector 
bundle $E$ over $X$. 

\subsection{The transmission condition}\label{section2.2}
 
Let us  denote by $2X$ the double of $X$ (or any smooth manifold without boundary which 
contains $X$ as a submanifold). With a classical parameter-dependent pseudodifferential operator 
$P(\tau)$ on $2X$ we can form  $P_{+}(\tau)=r^{+}P(\tau)e^{+}$, where $r^{+}$ and $e^{+}$ are 
the operators of restriction to $X$ and extension-by-zero to $2X$, respectively. In general, 
$P_{+}(\tau)$ does not map the space $\calC^\infty(X)$ to itself: when this is true, $P(\tau)$ 
is said to satisfy the transmission property (cf. \cite{Horm}). We will assume the so-called 
(parameter-dependent) \textit{transmission condition}, described in Definition \ref{deftp} 
below in terms of properties of the local symbols of $P$. 

\begin{definition}\label{deftp} 
A classical parameter-dependent symbol 
$p(x,\xi,\tau)\in S^\mu_{\mathrm{cl}}(\rz^n\times\rz^n_\xi;\rz^l_\tau)$ with integer order 
$\mu\in\gz$, such that $p\sim\sum\limits_{j\in\nz_0}p_{(\mu-j)}$ with $p_{(\mu-j)}$ being 
positively homogeneous of degree $\mu-j$ in $(\xi,\tau)$ for $(\xi,\tau)\not=0$, has the 
transmission condition at $x_{n}=0$ if, 
\begin{equation}\label{trp}
 D^k_{x_{n}}D^\alpha_{(\xi,\tau)}p_{(\mu-j)}(x^\prime,0,0,0,1,0)=
 (-1)^{\mu-j-|\alpha|}
 D^k_{x_{n}}D^\alpha_{(\xi,\tau)}p_{(\mu-j)}(x^\prime,0,0,-1,0)\qquad\forall\;k,\alpha. 
\end{equation}
Here, we used the splittings $x=(x^\prime,x_n)$ and $\xi=(\xi^\prime,\xi_n)$. 
\end{definition}  

The difference between the transmission condition and the  transmission property is that the 
first ensures the transmission property for $P$ on both sides of $\partial X$. For a detailed 
description of the transmission condition see, e.g., \cite{Grub1}, \cite{GrHo}, \cite{Schr}.

\subsection{Parameter-dependent Boutet de Monvel's algebra}\label{section2.3}

The parameter-dependent version of Boutet de Monvel's calculus deals with families of 
block-matrix operators $\calA=\calA(\tau)$, $\tau \in \rz^l$, of the form  
\begin{equation}\label{bdma}
 \calA = \begin{pmatrix}P_{+}(\tau)+G(\tau) & K(\tau)\\T(\tau)&Q(\tau)\end{pmatrix}:
 \begin{matrix}\ci(X,E_{0})\\ \oplus\\ \ci(\partial X, F_{0})\end{matrix} 
 \longrightarrow 
 \begin{matrix}\ci(X,E_{1})\\ \oplus\\ \ci(\partial X, F_{1})\end{matrix}, 
\end{equation} 
where $E_{0},E_{1}$ are vector bundles over $X$ and $F_{0},F_{1}$ are  vector bundles over 
$\partial X$ (possibly zero-dimensional). As described above, $P_{+}(\tau)=r^{+}P(\tau)e^{+}$, 
satisfies the transmission condition. $G(\tau)$ is an operator family on $X$ called 
(\textit{parameter-dependent}) \textit{singular Green operator}, arising, in particular, in 
the composition of two block-matrix operators by 
$(PP^\prime)_{+}(\tau) -  P_{+}(\tau)P^\prime_{+}(\tau)$. Moreover, $T(\tau)$ is a 
(\textit{parameter-dependent})  \textit{trace operator},  $K(\tau)$ a 
(\textit{parameter-dependent}) \textit{Poisson operator}, and $Q(\tau)$ is a 
parameter-dependent pseudodifferential operator on $\partial X$. 

We now describe the various entries. For simplicity, we switch to the half-space case 
$X=\rpbar^n$ (local situation) and take all bundles trivial one dimensional.  

\textbf{Trace operators.}  A parameter-dependent  trace operator of order $\mu\in\rz$ and 
type $d\in\nz$ is a family of operators $T(\tau):\cicomp(\rpbar^n)\to\ci(\rz^{n-1})$, 
$\tau\in\rz^l$, defined as 
 $$T(\tau) = \smsum_{j=0}^{d-1} S_{j}(\tau) \gamma_{j} + \wt T(\tau), $$ 
where $\gamma_{j}$ is as in Remark \ref{trace}, $S_{j}(\tau)$ is a parameter-dependent 
pseudodifferential operator  on $\rz^{n-1}$ of order $\mu-j$ and $\wt T(\tau)$ is an operator 
family of the form 
 $$(\wt T(\tau) u)(x^\prime) = \int_{\rz^{n-1}}\int_{0}^\infty 
   e^{i x^\prime \xi^\prime}\tilde{t}(x^\prime, x_{n}, \xi^\prime,\tau) \,
  (\calF_{x^\prime\to\xi^\prime}u)(\xi^\prime, x_{n})  \, dx_{n} \dbar\xi^\prime, $$ 
where $\tilde{t}$ is in $\calS(\rpbar)$ as a function of $x_{n}$ and satisfies 
 $$\tilde{t} (x^\prime,x_{n},\xi^\prime,\tau)=[\xi^\prime,\tau]^{\frac{1}{2}}
   t(x^\prime,\xi^\prime,\tau;[\xi^\prime,\tau]x_n)$$
with 
 $$t(x^\prime,\xi^\prime,\tau;x_n)\in
   S_{\textrm{cl}}^{\mu+\frac{1}{2}}(\rz^{n-1}_{x^\prime}\times\rz^{n-1}_{\xi^\prime};
   \rz^l_{\tau})\pit\calS(\rz_{+}).$$
Here, $[\,\cdot\,]$ denotes a smooth positive function with $[\,\cdot\,]=|\cdot|$ outside a 
neighborhood of the origin. Note that this structure implies the estimates 
\begin{equation*}
 \| x_{n}^l D^{l^\prime}_{x_{n}} D^\beta_{x^\prime} D^\alpha_{(\xi^\prime,\tau)} \,
 \tilde{t}(x^\prime, x_{n},\xi^\prime,\tau) \|_{L_2(\rz_{+,x_n})} \le 
 c\,\spk{\xi^\prime,\tau}^{\mu+\frac{1}{2}-l+l^\prime-|\alpha|} 
\end{equation*}
for every choice of the indices $l,l^\prime,\alpha,\beta$ with a resulting constant $c$. 
The \textit{principal boundary symbol} of $T$, 
 $$\sigma_\partial^\mu(T)(x^\prime,\xi^\prime,\tau):\calS(\rpbar)\longrightarrow\cz,$$
defined for $(\xi^\prime,\tau)\not=0$, is given by 
 $$u\mapsto 
   \smsum_{j=0}^d\sigma_\psi^{\mu-j}(S_j)(x^\prime,\xi^\prime,\tau)\gamma_j u+
   |(\xi^\prime,\tau)|^{\frac{1}{2}}\int_0^\infty\sigma_\psi^{\mu+\frac{1}{2}}(t)
   (x^\prime,\xi^\prime,\tau;|(\xi^\prime,\tau)|x_n)u(x_n)\,dx_n.$$
Here, $\sigma_\psi^m$ denotes the usual homogeneous principal symbol of order $m$ and 
$\gamma_j$ acts on functions on $\rpbar$. 
With these choices, it turns out that any operator of the form  
$T(\tau)=\gamma_{0} P_{+}(\tau)$ is a parameter-dependent trace operator if $P$ satisfies 
the transmission condition.  

\textbf{Poisson operators.} A parameter-dependent Poisson operator of order $\mu\in\rz$ is of 
the form 
\begin{equation*}
 (K(\tau)v)(x^{\prime},x_{n})=\int_{\rz^{n-1}} e^{i x^\prime \xi^\prime}
 \tilde{k}(x^\prime,x_{n},\xi^\prime,\tau)\hat{v}(\xi^{\prime})\dbar \xi^\prime, 
\end{equation*} 
where $\tilde{k}$ has the same structure as $\tilde{t}$ above, but with $\mu-1$ in place of 
$\mu$. The {\em principal boundary symbol}, defined for $(\xi^\prime,\tau)\not=0$, is 
 $$\sigma_\partial^\mu(K)(x^\prime,\xi^\prime,\tau):\cz\longrightarrow\calS(\rpbar),\quad
   a\mapsto|(\xi^\prime,\tau)|^{\frac{1}{2}}
   \sigma_\psi^{\mu-\frac{1}{2}}(k)
   (x^\prime,\xi^\prime,\tau,|(\xi^\prime,\tau)|\,{\cdot}\,)\,a.$$

\textbf{Singular Green operators.} A parameter-dependent singular Green operator of order 
$\mu\in\rz$ and type $d\in\nz_0$ is a family  
 $$G(\tau)=\smsum_{j=0}^{d-1}K_{j}(\tau)\gamma_{j} + \wt G(\tau) $$ 
with parameter-dependent Poisson operators $K_{j}(\tau)$ of order $\mu-j$, the standard trace 
operators $\gamma_{j}$ and an operator $\wt G(\tau)$ of the form 
 $$(\wt G(\tau) u)(x) = \int_{\rz^{n-1}}\int_{0}^\infty e^{ix^\prime\xi^\prime}
   \tilde{g}(x^\prime,x_{n},y_{n},\xi^\prime,\tau)(\calF_{x^\prime\to\xi^\prime}u)
   (\xi^\prime, y_{n}) dy_{n} \dbar \xi^\prime.$$ 
The symbol-kernel $\tilde{g}$ is of the form 
 $$\tilde{g}(x^\prime,x_n,y_n,\xi^\prime,\tau)=
   [\xi^\prime,\tau]\,g(x^\prime,\xi^\prime,\tau,[\xi^\prime,\tau]x_n,[\xi^\prime,\tau]y_n)$$
with 
 $$g(x^\prime,\xi^\prime,\tau,x_n,y_n)\in 
   S^\mu_{\textrm{cl}}(\rz^{n-1}_{x^\prime}\times\rz^{n-1}_{\xi^\prime};\rz^l_\tau)\pit
   \calS(\rpbar\times\rpbar).$$
In particular, $\tilde{g}$ satisfies, for all indices, the estimates 
 $$\| x_{n}^k D^{k^\prime}_{x_{n}} y_{n}^m D^{m^\prime}_{y_{n}}
   D^\beta_{x^\prime}D^\alpha_{(\xi^\prime,\tau)}
   \tilde{g}(x^\prime,x_{n},y_{n},\xi^\prime,\tau)\|_{L_2(\rz_{+,x_n}\times\rz_{+,y_n})}
   \le c \spk{\xi^\prime,\tau}^{\mu-k+k^\prime-m+m^\prime-|\alpha|}.$$
The {\em principal boundary symbol} of $G$, 
 $$\sigma_\partial^\mu(G)(x^\prime,\xi^\prime,\tau):
   \calS(\rpbar)\longrightarrow\calS(\rpbar),$$
defined for $(\xi^\prime,\tau)\not=0$, is given by 
 $$u\mapsto \smsum_{j=0}^{d-1}\sigma_\partial^{\mu-j}(K_j)(x^\prime,\xi^\prime,\tau)\gamma_ju
   \,+\,|(\xi^\prime,\tau)| \int_0^\infty
   \sigma_\psi^\mu(g)(x^\prime,\xi^\prime,\tau,|(\xi^\prime,\tau)|\,\cdot\,,
   |(\xi^\prime,\tau)|y_n)u(y_n)\,dy_n.$$

The {\em principal boundary symbols} of the remaining elements are 
 $$\sigma_\partial^\mu(P_+)(x^\prime,\xi^\prime,\tau)=
   r^+\,\sigma_\psi^\mu(P)(x^\prime,D_{x_n},\xi^\prime,\tau)\,e^+: 
   \calS(\rpbar)\longrightarrow\calS(\rpbar),$$
and $\sigma_\partial^\mu(Q)(x^\prime,\xi^\prime,\tau)=
\sigma_\psi^\mu(Q)(x^\prime,\xi^\prime,\tau)$ is the usual homogeneous principal symbol. 

The definitions made above are all invariant under smooth change of coordinates which 
preserve the boundary, so operators can be defined on the manifold $X$ via a partition of unity. 

\begin{definition}\label{defbdm}
$B^{\mu,d}(X;\rz^l;E_0,F_0;E_1,F_1)$, $\mu\in\gz$ and $d\in\nz_0$, consists of all 
operator-families $\calA=\calA(\tau)$ of the form \eqref{bdma} with each entry being of order 
$\mu$ and type $d$.  
\end{definition}  

\subsection{Algebra and mapping properties, ellipticity}\label{section2.4}

The next result is due to \cite{GrKo}.  

\begin{theorem}\label{mapping}
Let $\calA \in B^{\mu,d}(X;\rz^l;E_0,F_0;E_1,F_1)$. 
Then, for each $\tau$,  $\calA(\tau)$ extends to continuous operators
\begin{equation}\label{cont} 
 \calA(\tau):
 \begin{matrix}H^s_{p}(X,E_0)\\ \oplus\\ B^{s-\frac{1}{p}}_{p}(\partial X,F_0)\end{matrix}
 \longrightarrow
 \begin{matrix}H^{s-\mu}_{p}(X,E_1)\\ \oplus\\ B^{s-\mu-\frac{1}{p}}_{p}(\partial X,F_1)
 \end{matrix},\qquad s > d-1+\frac{1}{p}.
\end{equation}
\end{theorem} 

The notion algebra refers to the following behaviour under composition: 

\begin{theorem}\label{bdmcomp}
Let $(E_{0},F_{0})$, $(E_{1},F_{1})$ and $(E_{2},F_{2})$ be pairs of vector bundles (on $X$ and 
$\partial X$, respectively) such that composition of operators 
 $$\calA_{0}:
   \begin{matrix}\ci(X,E_1)\\ \oplus 	\\ \ci(\partial X, F_1)\end{matrix}
   \longrightarrow
   \begin{matrix}\ci(X,E_2)\\ \oplus 	\\ \ci(\partial X, F_2)\end{matrix},\qquad
   \calA_{1}:
   \begin{matrix}\ci(X,E_0)\\ \oplus 	\\ \ci(\partial X, F_0)\end{matrix}
   \longrightarrow
   \begin{matrix}\ci(X,E_1)\\ \oplus 	\\ \ci(\partial X, F_1)\end{matrix}$$
makes sense. The pointwise composition 
$(\calA_{0}(\tau),\calA_{1}(\tau))\mapsto \calA_{0}(\tau)\,\calA_{1}(\tau)$ then induces a map 
 $$B^{\mu_{0},d_{0}}(X;\rz^l;E_1,F_1;E_2,F_2)\times B^{\mu_{1},d_{1}}(X;\rz^l;E_0,F_0;E_1,F_1)
   \longrightarrow B^{\mu,d}(X;\rz^l;E_0,F_0;E_2,F_2),$$ 
with $\mu=\mu_{0}+\mu_{1}$ and $d = \max\{\mu_{1}+d_{0}, d_{1}\}$. 
\end{theorem}  

Assuming all bundles to be equipped with a hermitian metric (compatible with the product 
structure near the boundary), we have scalar products in $L_2(X,E_j)$ and 
$H^{-\frac{1}{2}}_2(\partial X,F_j)$. If $\calA$ is an operator as in \eqref{cont} of order 
and type 0, we can define its adjoint $\calA^*$ by 
 $$\skp{\calA u_1}{u_2}_{L_2(X,E_1)\oplus H^{-1/2}_2(\partial X,F_1)}=
   \skp{u_1}{\calA^* u_2}_{L_2(X,E_0)\oplus H^{-1/2}_2(\partial X,F_0)},$$
where $u_j\in\ci(X,E_j)\oplus\ci(\partial X,F_j)$ are arbitrary. Note that we take the scalar 
product of $H^{-\frac{1}{2}}_2(\partial X,F_j)=B^{-\frac{1}{2}}_2(\partial X,F_j)$, since by 
this the duality of $B^{s-\frac{1}{p}}_p(\partial X,F_j)$ and 
$B^{-s-\frac{1}{p^\prime}}_{p^\prime}(\partial X,F_j)$ is induced. As the following theorem 
states, the subalgebra of such operators is closed under taking adjoints. 

\begin{theorem}
Let $\mu\le0$. Taking pointwise the formal adjoint induces a map
 $$\calA\mapsto \calA^*:B^{\mu,0}(X;\rz^l;E_0,F_0;E_1,F_1)\longrightarrow 
   B^{\mu,0}(X;\rz^l;E_1,F_1;E_0,F_0).$$
\end{theorem}

Finally, we give a brief account to ellipticity. It is determined by the invertibility of two 
principal symbols. If  $\calA\in B^{\mu,d}(X;\rz^l;E_0,F_0;E_1,F_1)$, its 
{\em homogeneous principal symbol} is 
 $$\sigma_\psi(\calA)=\sigma_\psi^\mu(P):\pi^*_X E_0\longrightarrow\pi^*_X E_1,$$
where $\pi^*_X$ denotes the bundle pull-back under the canonical projection 
$\pi_X:(T^*X\times\rz^l)\setminus0\to X$. The {\em principal boundary symbol}
 $$\sigma_\partial^\mu(\calA):
   \pi^*_{\partial X}\begin{pmatrix}E_0|_{\partial X}\otimes\calS(\rpbar)\\ \oplus\\F_0
   \end{pmatrix}\longrightarrow
   \pi^*_{\partial X}\begin{pmatrix}E_1|_{\partial X}\otimes\calS(\rpbar)\\ \oplus\\F_1
   \end{pmatrix},$$
$\pi_{\partial X}:(T^*\partial X\times\rz^l)\setminus0\to\partial X$ denoting the canonical 
projection, is defined by 
 $$\sigma_\partial^\mu(\calA)(x^\prime,\xi^\prime,\tau)=
   \begin{pmatrix}\sigma_\partial^\mu(P_++G)&\sigma_\partial^\mu(K)\\
    \sigma_\partial^\mu(T)&\sigma_\partial^\mu(Q)\end{pmatrix}(x^\prime,\xi^\prime,\tau)$$
with the single components as introduced in Section \ref{section2.3}. 

\begin{definition}
$\calA\in B^{\mu,d}(X;\rz^l;E_0,F_0;E_1,F_1)$ is parameter-dependent elliptic if both 
$\sigma_\psi(\calA)$ and $\sigma_\partial^\mu(\calA)$ are invertible. 
\end{definition}

As for the usual pseudodifferential calculus, ellipticity implies (in fact, is equivalent to) 
the existence of a parametrix within the calculus. 

\begin{theorem}
Let $\calA\in B^{\mu,d}(X;\rz^l;E_0,F_0;E_1,F_1)$ with $d=\max(\mu,0)$ be elliptic. Then there 
is an operator $\calB\in B^{-\mu,d^\prime}(X;\rz^l;E_1,F_1;E_0,F_0)$, $d^\prime=\max(-\mu,0)$, 
such that 
 $$\calA\calB-I \in B^{-\infty,d^\prime}(X;\rz^l;E_1,F_1;E_1,F_1),\qquad 
   \calB\calA-I \in B^{-\infty,d}(X;\rz^l;E_0,F_0;E_0,F_0).$$
In particular, $\calA(\tau)$ acting as in \textrm{\eqref{cont}} is a Fredholm operator of 
index 0 for each $\tau$ and is invertible for $|\tau|$ sufficiently large. 
\end{theorem}


\begin{small}
\bibliographystyle{amsalpha}

\end{small}


\end{document}